\documentclass{article}  
\usepackage{amssymb}
\usepackage{amsthm} 
\usepackage{amsmath}
\usepackage{graphics} 
\usepackage{eepic}    
\usepackage[all,arc,poly]{xy}

\newcommand{\nc}{\newcommand}
\nc{\look}{\marginpar{$\bullet$}} 

\nc{\Section}{\section}
\nc{\SubSection}{\subsection}

\newtheorem{theo}{Theorem}[section]   
\newtheorem{ddef}[theo]{Definition}
\newtheorem{llem}[theo]{Lemma} 
\newtheorem{oobs}[theo]{Observation} 
\newtheorem{rrem}[theo]{Remark} 
\newtheorem{prop}[theo]{Proposition} 
\newtheorem{ccor}[theo]{Corollary}  
\newtheorem{claim}[theo]{Claim} 
\newtheorem{qquest}[theo]{Question} 
\newtheorem{fact}[theo]{Fact} 
\newtheorem{pprov}[theo]{Proviso}
\newtheorem{eexam}[theo]{Example} 

\nc{\bT}{\begin{theo}} 
\nc{\eT}{\end{theo}}
\nc{\bD}{\begin{ddef} \rm }
\nc{\eD}{\end{ddef}}
\nc{\bC}{\begin{ccor}}
\nc{\eC}{\end{ccor}}
\nc{\bCl}{\begin{claim}}
\nc{\eCl}{\end{claim}}
\nc{\bQ}{\begin{qquest}}
\nc{\eQ}{\end{qquest}}
\nc{\bL}{\begin{llem}}
\nc{\eL}{\end{llem}}
\nc{\bP}{\begin{prop}}
\nc{\eP}{\end{prop}}
\nc{\bR}{\begin{rrem}}
\nc{\eR}{\end{rrem}}
\nc{\bO}{\begin{oobs}}
\nc{\eO}{\end{oobs}}
\nc{\bF}{\begin{fact}}
\nc{\eF}{\end{fact}}
\nc{\bProv}{\begin{pprov}}
\nc{\eProv}{\end{pprov}}
\nc{\bE}{\begin{eexam} \rm }
\nc{\eE}{\end{eexam}}

\nc{\prf}{\begin{proof}}
\nc{\eprf}{\end{proof}}

\renewcommand{\phi}{\varphi}
\renewcommand{\geq}{\geqslant}
\renewcommand{\leq}{\leqslant}

\newcommand{\strictsubset}{\varsubsetneq}
\renewcommand{\subset}{\subseteq}

\renewcommand{\supset}{\supseteq}

\newenvironment{romanenumerate}%
{\begin{list}{(\roman{enumi})}{\usecounter{enumi}
\setlength{\labelwidth}{2cm}
\setlength{\itemindent}{0pt}
\setlength{\itemsep}{0.5\itemsep}
\setlength{\topsep}{\itemsep}
\setlength{\parsep}{0pt}
}}{\end{list}}
\nc{\bre}{\begin{romanenumerate}}
\nc{\ere}{\end{romanenumerate}}
\newenvironment{alphaenumerate}%
{\begin{list}{(\alph{enumii})}{\usecounter{enumii}
\setlength{\labelwidth}{2cm}
\setlength{\itemindent}{0pt}
\setlength{\itemsep}{0.5\itemsep}
\setlength{\topsep}{\itemsep}
\setlength{\parsep}{0pt}
}}{\end{list}}
\nc{\bae}{\begin{alphaenumerate}}
\nc{\eae}{\end{alphaenumerate}}
\newenvironment{numenumerate}%
{\begin{list}{(\arabic{enumiii})}{\usecounter{enumiii}
\setlength{\labelwidth}{2cm}
\setlength{\itemindent}{0pt}
\setlength{\itemsep}{0.5\itemsep}
\setlength{\topsep}{\itemsep}
\setlength{\parsep}{0pt}
}}{\end{list}}
\nc{\bne}{\begin{numenumerate}}
\nc{\ene}{\end{numenumerate}}

\nc{\ins}[1]{\bigskip\noindent
\framebox{\begin{minipage}{.95\textwidth} \sloppy \noindent \em #1 \end{minipage}}\bigskip}

\nc{\str}[1]{{\mathcal{#1}}}
\nc{\brck}[1]{[\![ #1 ]\!]}
\nc{\restr}{\!\restriction\!}
\nc{\A}{\mathbb{A}}
\nc{\B}{\mathbb{B}}
\nc{\G}{\mathbb{G}}
\renewcommand{\H}{\mathbb{H}}
\nc{\I}{\mathbb{I}}
\nc{\U}{\mathbb{U}}
\nc{\V}{\mathbb{V}}
\nc{\abar}{\mathbf{a}}
\nc{\bbar}{\mathbf{b}}
\nc{\cbar}{\mathbf{c}}
\nc{\xbar}{\mathbf{x}}
\nc{\ybar}{\mathbf{y}}
\nc{\zbar}{\mathbf{z}}

\nc{\barr}{\begin{array}}
\nc{\earr}{\end{array}}
\nc{\btab}{\begin{tabular}}
\nc{\etab}{\end{tabular}}

\nc{\nothing}{\rule{0em}{1ex}}
\nc{\highnothing}{\rule{0em}{3ex}}
\nc{\hnt}{\highnothing}
\nc{\nt}{\nothing}
\nc{\nnt}{\rule{.1pt}{0pt}}

\renewcommand{\sc}{\scriptstyle}
\nc{\ssc}{\scriptscriptstyle}
\nc{\N}{{\mathbb N}}
\nc{\Z}{{\mathbb Z}}
\nc{\cym}{\mathbb{G}}

\renewcommand{\epsilon}{\varepsilon}

\begin{document}

\title{Amalgamation and Symmetry: From Local to Global Consistency in The Finite}
\author{Martin Otto\thanks{Research partially supported 
by DFG grant OT~147/6-1: \emph{Constructions and Analysis in
  Hypergraphs of Controlled Acyclicity.} I also gratefully acknowledge participation
in a programme on \emph{Logical Structure in Computation} at the
Simons Institute in Berkeley in 2016, which helped to refine some of
the ideas in this work.}%
\\
Department of Mathematics\\
Technische Universit\"at Darmstadt}
\date{July~2024}

\maketitle

\begin{abstract}
\noindent
We present a generic construction of finite realisations of 
amalgamation patterns. 
An amalgamation pattern is specified by a finite collection of finite
template structures together with a collection of partial isomorphisms 
between them. A realisation is a globally consistent solution to the locally
consistent specification of this amalgamation problem:
this is a single structure equipped with an atlas of
distinguished substructures associated with the template structures, 
their overlaps realising precisely the identifications 
induced by the given partial isomorphisms.
Our construction is based on natural reduced
products with suitable groupoids. The resulting realisations are generic 
in the sense that they can be made to preserve all symmetries of
the specification. They can also be made to be universal w.r.t.\ to local 
homomorphisms up to any specified size. As key
applications of the main construction we discuss finite hypergraph 
coverings of specified levels of acyclicity and a new route to the
lifting of local symmetries to global automorphisms in finite
structures.
\end{abstract}

\vspace{.25cm}

\vspace{.5cm}\hnt\quad
\btab[t]{@{}r@{}} mathematical subject classification: 
primary 20L05, 05C65, 05E18;\\
secondary 20B25, 20F05, 57M12 \etab

\pagebreak   

\tableofcontents 

\pagebreak 

\section{Introduction}

At the centre of this investigation are generic algebraic-com\-bi\-na\-to\-rial 
constructions of finite amalgams of finite structures. These
amalgams are based on gluing instructions specified as free amalgamations 
between pairs of given templates. The input data for this 
problem, which we call \emph{amalgamation patterns},
consist of a finite family of finite
structures as templates, equipped with a set of partial isomorphism
between pairs of these. An amalgamation pattern serves as a local 
specification of the overlap pattern between designated parts of 
a desired global solution or realisation.
Such a \emph{realisation} consists of a single finite 
structure made up of a union of designated substructures; 
each designated substructure is isomorphically related to a member of the 
given family of templates; and the overlaps between these designated
substructures is in accordance with the given family of partial 
isomorphisms --- as prescriptions for the overlap in pairwise
amalgamations. 
Some aspects of the combinatorial core of the matter can be illustrated in
the most fundamental instance of such a problem: here the input data 
consists of a family of just disjoint sets together with designated
partial bijections between pairs of these sets. In this scenario,
a realisation is just a hypergraph: its hyperedges 
represent the co-ordinate domains in an atlas of charts into the given
family of sets, so that changes of co-ordinates in their overlaps 
are induced by the given partial bijections. In the general case, we provide generic
constructions of finite realisations that allow us to respect all intrinsic
symmetries of the specification and to achieve any degree of local acyclicity
for the global intersection pattern between the constituent substructures.
The main constructions are based on reduced products with suitable finite 
groupoids. The existence of such finite groupoids, which need to have 
strong acyclicity properties in order to support global consistency, 
is established in~\cite{Cayleynew20}. 
That new approach replaces constructions proposed
in~\cite{lics2013}, which attempted to generalise more directly
corresponding constructions from~\cite{OttoJACM} 
but contained a serious flaw.%
\footnote{\newline
I am very grateful to Julian~Bitterlich, 
who prominently used these results in~\cite{Bitterlich}, for  
discovering a serious mistake in the proposed construction of
$N$-acyclic groupoids from~\cite{lics2013}, which eventually
led to the substantially new approach in~\cite{Cayleynew20}.}
The upshot of the construction
is summarised in Theorem~\ref{Nacycgroipoidthm}. 
Our main theorem here, Theorem~\ref{realisethm}, states that any
finite amalgamation pattern possesses highly symmetric realisations, 
obtained as reduced products with suitable groupoids, with additional
local acyclicity and universality properties.
As one key application we obtain finite branched hypergraph coverings 
of any desired local degree of hypergraph acyclicity. 
As a further application 
we obtain a novel approach to extension properties for partial
isomorphisms, which provide liftings of local symmetries 
to global symmetries (i.e., from partial to full automorphisms)
in finite extensions. Our generic realisations of amalgamation
patterns induced by such extension tasks (EPPA tasks in the sense of 
Herwig and Lascar~\cite{HerwigLascar}) yield not just a new route to
such extensions in the style of the powerful Herwig--Lascar theorem, 
but apparently more generic solutions. Indeed, our formulations
in Theorem~\ref{eppathm} and Corollary~\ref{eppacorr}
are more specific w.r.t.\ the symmetries involved, 
w.r.t.\ to the local-to-global relationship between the 
parts and the whole, and w.r.t.\ their universality properties. 

\medskip
At the algebraic and combinatorial level, our constructions 
illustrate new uses of finite groupoids and associated graph
and hypergraph structures. Analogous to the established use of 
Cayley graphs as a combinatorial representation of groups, 
we can use a Cayley graph representation of groupoids --- and a dual 
picture in a hypergraph of cosets --- to capture the specific acyclicity
properties that are essential for our constructions.
In controlling cyclic configurations in the dual hypergraph, i.e.\
cycles formed by overlapping cosets, the relevant acyclicity 
requirements exceed the scope of classical Cayley graph constructions 
as discussed e.g.\ in~\cite{Alon95}. Classical notions of acyclicity 
in Cayley graphs concern the graph-theoretic girth and 
measure the length of generator cycles. Here we need to control
\emph{coset cycles} and measure the length of
cycles formed by overlapping cosets in the Cayley graph. The technical
challenge lies in the fact that the building blocks of the relevant
cyclic configurations are elements of the dual hypergraph, and thus 
second-order objects in the Cayley graph itself, and that their size 
cannot a priori be bounded. Correspondingly, these acyclicity criteria are 
tailored to combinatorial contexts where local acyclicity for
decompositions of hypergraph-like rather than graph-like structures
matters. Here we use reduced products with the Cayley graphs of 
suitable finite groupoids to obtain natural realisations. 
Our first application then shows how this approach naturally lends itself
to the construction of finite hypergraph coverings of controlled
acyclicity. These coverings
provide interesting classes of highly homogeneous and highly acyclic finite 
hypergraphs. As very adaptable synthetic constructions of locally acyclic
hypergraphs these can play a r\^ole analogous to that played by known
constructions of Cayley graphs of large girth in the setting of
graphs. They are a source of generic examples of interest in relation to
structural decomposition techniques in combinatorics and algorithmic 
model theory. Indeed, hypergraph acyclicity has long been
recognised as an important criterion in combinatorial, algorithmic and
logical contexts~\cite{Berge,BeeriFaginetal}. In the well understood
setting of graphs and graph-like structures, local acyclicity can be achieved in finite
bisimilar unfoldings or coverings based on Cayley groups of large
girth~\cite{OttoAPAL04}. We here see that Cayley groupoids and the 
notion of coset acyclicity support the
adequate generalisations in the hypergraph setting.
Applications of closely related
structural transformations to questions in logic, and especially in finite
model theory, have already been explored, e.g.,
in~\cite{HO,OttoJACM,lics2013,GO14,CO17} 
w.r.t.\ expressive completeness results as well as w.r.t.\ 
algorithmic issues~\cite{BGO,BtCO}. 
The genericity of our local-to-global
constructions as exemplified in their
applications to extension problems for local automorphisms 
may have further applications in the study of
automorphisms of countable structures built from finite
substructures and of amalgamation classes that arise in the 
model-theoretic and algebraic analysis of homogeneous 
structures \cite{MacPherson}. 

\medskip
At the more conceptual level, the groupoidal constructions presented here
suggest interesting discrete and even finite analoga of classical
concepts like branched coverings~\cite{Fox}, notions of path
independence and contractability, and the study of local 
symmetries~\cite{Lawson}, which invite further exploration.
The alge\-bra\-ic-com\-bi\-na\-to\-ri\-al approach to phenomena of local 
versus global consistency in finite relational
structures may also point towards potential applications in 
relational models for quantum information theory as proposed in~\cite{Abramsky}. 

An amalgamation pattern 
serves as a local specification of \emph{pairwise} local amalgamation steps.
Correspondingly, the pattern itself is indexed by an underlying 
graph-like structure of sites and links, the \emph{incidence patterns} of
Definition~\ref{incpatdef}. This stands to the amalgamation task 
embodied in the amalgamation pattern as the intersection graph of a
hypergraph stands to the actual hypergraph. And just as branched 
hypergraph coverings can be associated with bisimilar coverings 
at the level of the underlying intersection graphs~\cite{OttoJACM}, 
our realisations of amalgamation patterns involve bisimilar 
local unfoldings at the level of the underlying incidence patterns. 
Notions of bisimilarity and local unfoldings thus also play a 
crucial r\^ole in the context of local versus global symmetries and 
of local versus global consistency in general relational structures.

\section{Amalgamation patterns and their realisations} 

\paragraph*{Basic conventions.}
Dealing with relational structures we always appeal to a 
fixed finite relational signature $\sigma$ and typically write
$\str{A} = (A, (R^{\str{A}})_{R \in \sigma})$ for a $\sigma$-structure over
universe $A$ with interpretations $R^{\str{A}} \subset A^r$ for
a relation symbol $R \in \sigma$ of arity $r = \mathrm{ar}(R)$.   
By \emph{substructures} $\str{A}_0 \subset \str{A}$ we 
always mean induced substructures $\str{A}_0 = \str{A}\restr A_0$ for
which $R^{\str{A}_0} = R^{\str{A}} \cap A_0^r$. This should be contrasted 
with \emph{weak substructures} (sometimes considered to be 
the standard notion, e.g.\ for subgraphs) where  just $R^{\str{A}_0} \subset
R^{\str{A}} \cap A_0^r$ is required. 
We shall often have occasion to consider more complex structural scenarios
that involve e.g.\ some $\sigma$-structure together with a family of
substructures, or families of partial isomorphisms within some
$\str{A}$ or between the members of a family of $\sigma$-structures, et cetera.
In such situations we shall often adopt a multi-sorted 
formalisation, with distinct sorts for different kinds of objects and 
an explicit encoding of relationships between objects, e.g.\
by means of families of functions between different sorts. 
As we shall see, even the labelling of such families
should sometimes not be regarded as static but as a structural feature
that may for instance be subject to permutations if we want to account 
for all relevant symmetries. The idea is illustrated in the
following section for the key notion of \emph{amalgamation patterns}, where
we start from an ad-hoc preliminary formalisation in order to develop 
the more symmetry-aware formalisation. A \emph{partial isomorphism}
between $\sigma$-structure $\str{A}$ and $\str{A}'$ is a partial function $p$
from $A$ to $A'$ that is an isomorphism between the induced
substructures on its domain and image: $p \colon \str{A}\restr
\mathrm{dom}(p) \simeq \str{A}' \restr \mathrm{image}(p)$; we write 
$\mathrm{part}(\str{A},\str{A}')$ for the set of all partial
isomorphisms between $\str{A}$ and $\str{A}'$. For partial maps like
these, suggestive notation like $p \subset q$ is used to say that
$q$ is an extension of $p$.

The underlying structures that we think
of as given data for amalgamation tasks will crucially be finite, but in
order to discuss the intended solutions, which are also finite, we
occasionally refer to related infinite structures. For many basic
notions finiteness is not essential and we choose formulations that
could equally be applied in infinite settings. Specific conventions will be discussed in
context, but overall we just assume the standard terminology of
basic universal algebra.

\subsection{Amalgamation patterns}

Amalgamation patterns specify intended overlaps between substructures that
are modelled on templates. This specification is locally tight, and guaranteed to be
locally consistent in terms of local one-to-one overlaps, but leaves
under-specified the global overlap structure implicitly required for its
realisations, and may not in itself be globally consistent in any
straightforward sense. We give a preliminary definition to discuss
first examples at a more naive level before turning to the proper
definitions. Those
more formal definitions of amalgamation patterns and their
realisations will be slightly more involved in order to support 
the intended sensitivity to internal symmetries.

\bD [preliminary version]
\label{amalgpatpreldef}
An \emph{amalgamation pattern}  over some relational signature
$\sigma$ is a structure of type 
$\H = \bigl( \str{H}, (\str{A}_s)_{s \in S}, (\rho_e)_{e \in E}\bigr)$
describing \emph{links} between \emph{sites}. 
$\H$ consists of  a $\sigma$-structure $\str{H}$ that is the disjoint
union $\str{H} = \dot{\bigcup}_{s \in S} \str{A}_s$
of $\sigma$-structures $\str{A}_s = \str{H}\restr A_s$ as
\emph{sites}, based on a partition $H =\dot{\bigcup}_{s \in S} A_s$ into
non-empty $A_s$ for $s \in S$.
These substructures are related by a collection of \emph{links} $\rho_e \in
\mathrm{part}(\str{A}_s,\str{A}_{s'})$, which are 
partial isomorphisms between $\str{A}_s$ and $\str{A}_{s'}$.
The index set $E$ of links is partitioned according to 
$E = \dot{\bigcup}_{s,s' \in S} E[s,s']$ such that $E[s,s']$ is
the set of labels of links directed from $\str{A}_s$ to $\str{A}_{s'}$
(possibly empty).
\eD

\bE
\label{explodedviewex}
The \emph{exploded view} of a relational structure $\str{A}$ w.r.t.\ a
suitable family
of distinguished substructures $(\str{A}_s)_{s \in S}$ that cover
$\str{A}$ is a special
case of an amalgamation pattern. The $\str{A}_s =
\str{A}\restr s$, for a collection of subsets $s \in S
\subset \mathcal{P}(A)$, provide an \emph{atlas} (in a sense
to be defined below) if $\str{A} = \bigcup_{s \in S} \str{A}_s$. In this
case, we let $\str{H}$ be the disjoint union of $s$-tagged copies $(\str{A}\restr
s) \times \{ s \}$
of the
$\str{A}_s$:
\[
\textstyle
\str{H} := \dot{\bigcup}_{s \in S} (\str{A}\restr s \times \{ s \}).
\]
Then the natural underlying indexing by a graph structure $(S,E)$, 
where $E = \{ (s,s') \colon s \not= s', s \cap s' \not=
\emptyset \}$, together with  partial isomorphisms 
\[
\barr{rcl}
\rho_e \colon \str{A}_s \times \{ s \} &\stackrel{\ssc \mathrm{part}}{\longrightarrow}&
\str{A}_{s'} \times \{ s' \}
\\
(a,s) &\longmapsto&(a,s') \; \mbox{ for } a \in s \cap s', 
\earr
\] 
record the actual identifications between
elements of $\str{A}_s$ and $\str{A}_{s'}$ in $\str{A}$, for all
non-trivial combinations $s,s'$.
This very simple example holds some of the essential
intuition of how an amalgamation pattern arises as an overlap
specification, and also points us to what it should mean to \emph{realise}
an amalgamation pattern: clearly $\str{A}$
realises its exploded view specification.    
Due to its great simplicity, this class of examples trivially displays many 
of the special features that we shall want to guarantee through
specific pre-processing in the general case. 
\eE

\bE
\label{explodedviewhypex}
An even more special case ensues if we completely abstract away from 
relational content, i.e.\ for $\sigma = \emptyset$. Then we are
dealing we the \emph{exploded view of a hypergraph} $(A,S)$, which 
represents the hyperedges as disjoint sets together with an explicit
specification of overlaps between pairs of hyperedges.
\eE

\bD
\label{atlasdef}
An \emph{atlas} $\A$ for a $\sigma$-structure $\str{A}$ 
augments $\str{A}$ according to 
\[
\A = (\str{A},U,(U_s)_{s \in S}, (\pi_{u,s})_{u\in U_s})
\]
by a collection of \emph{charts}, i.e., isomorphisms
\[
\pi_{u,s} \colon \str{A}\restr u \simeq \str{A}_s
\]
between induced
substructures $\str{A}\restr u$ and members $\str{A}_s$ 
of some collection of external \emph{co-ordinate structures} $(\str{A}_s)_{s \in S}$. 
The \emph{co-ordinate domains} $u \subset A$ of the charts
form a superimposed hypergraph structure $(A,U)$ on $A$, 
with hyperedge set $U \subset \mathcal{P}(A)$ consisting of the
domains $u$ of the charts, which are required to cover $\str{A}$ in the sense 
that 
\[
\textstyle
\str{A} = \bigcup_{u \in U} \str{A}\restr u, 
\]
i.e., not just every element of $A$, but every tuple in the
interpretation of any relation over $\str{A}$ must be fully contained in
one of the hyperedges $u$. 
\eD

We note that there may be multiple charts (into distinct, but
necessarily isomorphic) co-ordinate structures $\str{A}_s$ on the same 
co-ordinate domain $u \subset A$. The collection $U$ of co-ordinate
domains thus becomes the union, but not necessarily a disjoint union, 
of subsets $U_s$ of those co-ordinate domains that are associated with 
the co-ordinate structure $\str{A}_s$ by some $\pi_{u,s}$, for $s \in S$.

One could also associate a collection of
\emph{changes of co-ordinates} with a given atlas, which would be a
collection of partial isomorphisms between co-ordinate structures $\str{A}_s$
as induced by pairs of charts $\pi_{u,s}\colon \str{A}\restr u \simeq \str{A}_s$ and  
$\pi_{u',s'} \colon \str{A}\restr u' \simeq \str{A}_{s'}$ 
with non-trivial intersection
$u \cap u'$ of co-ordinate domains.

\paragraph*{Extensive multi-sorted formalisations.}
In the following we shall want to avoid fixed labellings in the formalisation
of complex structures like amalgamation patterns or atlases. This is essential in order to remain sensitive to all relevant
  symmetries, including symmetries that permute the index structure 
 used for labelling purposes. For instance, in the harmless 
example of an atlas $\A$ for
 $\str{A}$: an automorphism of the relational structure $\str{A}$ 
could map each $u \in U_s$ to some $u' \in U_{s'}$ which may be
matched by a permutation on the index set $S$ and an isomorphism 
between the co-ordinate structures $\str{A}_s$ and $\str{A}_{s'}$ that 
may also be compatible with the chosen co-ordinate maps 
$\pi_{u,s}$ and $\pi_{u',s'}$ and possibly with changes of co-ordinates
(best  depicted in a commuting diagram). 
Such a symmetry would not be apparent if the labelling of objects 
(such as co-ordinate domains, co-ordinate structures and co-ordinate
maps and changes) were taken as static. It therefore makes sense for our concerns
to adopt a multi-sorted formalisation, which allows for various 
functional and relational links between and within individual sorts
and takes into account the sorting and typing of various composite objects
so as to trace their behaviour under the most general kinds of 
overall symmetries. The compressed presentation of a composite structure 
like an atlas as a labelled tuple of objects will still be convenient,
but we shall only regard it as a shorthand for its extended
multi-sorted formalisation, which in the case of an atlas would
naturally involve disjoint sorts for elements of $\str{A}$, for 
elements of the $\str{A}_s$, of $U \subset \mathcal{P}(A)$,  
of  a set of co-ordinate maps $P = \{ \pi_{u,s} \colon u \in U_s,
s \in S \}$, possibly also of corresponding changes of
co-ordinates, and a sort for the elements of the index set $S$ itself.

\medskip
Following this idea, we firstly make explicit the underlying sites-and-links
structure of an amalgamation pattern (and other related structures) 
by formalising it as a multigraph $\I = (S,E)$ according 
to the following definition. 
We shall then regard the amalgamation pattern $\H$ --- in the official
definition to be given in Definition~\ref{amalgpatdef} below --- as 
an \emph{amalgamation pattern $\H$ over $\I$}.

\bD
\label{incpatdef}
An \emph{incidence pattern} $\I$ is a finite two-sorted structure
\[
\I  = \bigl( S,E, \iota \colon E \rightarrow S^2, \cdot^{-1} \colon E
\rightarrow E\bigr)
\] 
with sorts $S$ (vertices, to be
viewed as types of \emph{sites} for amalgamation) and 
$E$ (edges, to be viewed as types of \emph{links}) 
connected by a pair of functions 
$\iota = (\iota_1,\iota_2) \colon E \rightarrow S^2$ for the
allocation of source vertex $\iota_1(e) \in S$ and target vertex
$\iota_2(e) \in S$ to every edge $e \in E$, and with an involutive 
operation of edge reversal on $E$, $e \mapsto e^{-1}$ subject to the
requirements that $(e^{-1})^{-1} = e$ and 
$\iota_1(e^{-1}) = \iota_2(e)$.
\eD

The incidence pattern $\I$ underlying the amalgamation pattern $\H$ according to
Definition~\ref{amalgpatpreldef} is $\I = (S,E,\iota,\cdot^{-1})$
where $S$ and $E$ are the index sets for sites and links in $\H$, 
with $\iota_i(e) = s_i$ determined by the requirement that 
$\rho_e \in \mathrm{part}(\str{A}_{s_1},\str{A}_{s_2})$
and $e^{-1}$ determined by the condition that
$\rho_{e^{-1}}= \rho_e^{-1}$. 
While we identify $\rho_{e^{-1}}$ with the inverse
$\rho_e^{-1}$ in the context of a given amalgamation pattern $\H$ over $\I$,
we keep in mind that at the level of the overlap
pattern $\I$ itself, the operation $\cdot^{-1}$ is just edge reversal
in a directed multigraph.

The following definition of an amalgamation pattern \emph{over} 
an incidence pattern in an explicitly multi-sorted format is guided by 
the above considerations concerning an adequate account of symmetries.

\bD
\label{amalgpatdef}
An \emph{amalgamation pattern} $\H$ over an
incidence pattern $\I$, denoted $\H/\I$ to make the reference to $\I$ explicit, 
is a multi-sorted structure 
\[
\bigl(\I; \str{H}, \delta \colon H \rightarrow S;
P, \eta \colon P \rightarrow E\bigr)
\]
where 
\bre
\item[--]
$\I  = (S,E,\iota,\cdot^{-1})$ is an incidence pattern (with sorts $S$ and $E$ as above);
\item[--]
$\str{H}$
is a $\sigma$-structure whose universe $H$ is $S$-partitioned 
by $\delta$ such that 
\[
\str{H} = \dot{\bigcup}_{s \in S} \str{A}_s
\]
is a disjoint union of non-empty $\sigma$-structures $\str{A}_s = \str{H} \restr \delta^{-1}(s)$;
\item[--]
$P$ is a collection of 
partial $\sigma$-isomorphisms between pairs of these $\str{A}_s$,
$E$-partitioned into singleton sets (i.e.~$E$-labelled) by 
$\eta$ such that 
\[
P = \{ \rho_e \colon e \in E \},
\]
where $\rho_e$ is the unique element of $\eta^{-1}(e)$ 
and such that $\rho_e \in \mathrm{part}(\str{A}_s,\str{A}_{s'})$ if 
$\iota(e) = (s,s')$,
and $\rho_{e^{-1}} = \rho_e^{-1}$.
\ere
\eD

We also think of an atlas as a
multi-sorted structure with sorts $A$, $\dot{\bigcup}_{s \in S} A_s$, 
$U = \bigcup_{s \in S}U_s$, 
$P= \{ \pi_{u,s} \colon u \in U, s \in S \}$ and the index sort $S$,
with the natural encoding of the labellings involved. We shall only 
formalise this explicitly for the r\^ole of an atlas in the 
context of a \emph{realisation}
of an amalgamation pattern, in Definition~\ref{realdef} below.

\paragraph*{Symmetries.} 
The explicitly multi-sorted formalisation of incidence and 
amalgamation patterns is meant to support a sufficiently broad notion
of symmetries. In particular, we do want to consider symmetries
of $\H/\I$ that permute the labels $s \in S$ of sites and $e \in E$ of links in
as much as such permutations respect the structure $\I$, i.e.\ induce a
symmetry of the underlying incidence pattern. This is
naturally reflected in automorphisms of the multi-sorted formalisation.

\bD
\label{amalgautdef}
A \emph{symmetry} of the incidence pattern  $\I = (S, E, \iota,\cdot^{-1})$ is an automorphism of this $2$-sorted
structure consisting of a pair of permutations $\pi = (\pi_S, \pi_E)$ 
of the two sorts $S$ and $E$ that are compatible with $\iota$ and $\cdot^{-1}$.

A \emph{symmetry} of an amalgamation pattern
$\H/\I = \bigl(\I, \str{H}, \delta, P, \eta\bigr)$
over the incidence pattern $\I$ 
is an automorphism $\pi$ of this multi-sorted structure. 
It is specified in terms of permutations $\pi = (\pi_S,\pi_E,\pi_H)$ 
of the three sorts $S$, $E$ and $H$ such that
$\pi_H$ is an automorphism of the relational structure $\str{H}$, 
the pair $(\pi_S, \pi_E)$ forms an automorphism of $\I$ and the triple 
$(\pi_S,\pi_E,\pi_H)$ is compatible with $\delta$ and $\eta$.
A symmetry of $\H/\I$ is \emph{$\I$-rigid} if its operation on $\I$ is
trivial, i.e.\ if it leaves the labelling of sorts and links fixed.
\eD

In the context of $\H/\I$, our notational convention for compositions
of partial bijections adheres to the format of multiplication/operation from
the right. For links $e_1,e_2 \in E$ that match in the sense that
$\iota_2(e_1) = \iota_1(e_2)$ we write $\rho_{e_1} \cdot \rho_{e_2}$ 
or just $\rho_{e_1}\rho_{e_2}$ for the natural (partial) composition
$\rho_{e_2} \circ \rho{e_1}$, which is a partial isomorphism 
from $\str{A}_s$ to $\str{A}_{s'}$ for 
$s := \iota_1(e_1)$ and $s' := \iota_2(e_2)$ (possibly empty).  
This convention naturally extends to \emph{walks} 
$w= e_1 \cdots e_n$ in $\I$, where the property of a walk implies  
that $\iota_2(e_i) = \iota_1(e_{i+1})$ for all $i < n$, so that 
$\rho_w := \Pi_{i=1}^n \rho_{e_i}$ stands for the composition 
\[
\rho_w = \Pi_{i=1}^n \rho_{e_i} = 
\rho_{e_n} \circ  \cdots \circ  \rho_{e_1},
\]
which is a partial isomorphism from $\str{A}_s$ to $\str{A}_{s'}$ 
where $s = \iota_1(e_1)$ and $s' = \iota_2(e_n)$. In the same vein, we
associate $\mathrm{id}_s:=\mathrm{id}_{A_s}$, the identity function on 
the domain $A_s$ of $\str{A}_s$, with 
the unique walk $\lambda_s$ of length $0$ at $s$.

\medskip
In Section~\ref{groupoidmonsec} we shall appeal to the inverse
semigroup structure generated by the $(\rho_e)_{e \in E}$, which may
be regarded as an inverse sub-semigroup of the symmetric inverse 
semigroup $I(X)$ over the set $X = H = \dot{\bigcup}_{s \in S} A_s$, 
see e.g.~\cite{Lawson}. 
Notions of \emph{coherence}, which are to be discussed in relation to global
consistency for an amalgamation pattern $\H$
in the next section, can be cast in terms of this 
inverse semigroup 
in relation to its generators $(\rho_e)_{e \in E}$. 
Its elements are the compositions of the $\rho_e$
with (partial) composition 
as the semigroup operation (cf.\ Definition~\ref {symminvsemdef}).
Note that all non-trivial compositions arise as compositions along
walks in $\I$: products (partial compositions) can only have non-empty 
results if at least sites match at the interface. This is going to be 
converted to a groupoidal formalisation in Section~\ref{groupoidmonsec}.
Unlike the situation in the groupoidal setting, compositions
of the $\rho_e$ in the inverse semigroup setting over $\H$ are typically
partial in restriction to their sites so that, for instance for $e
\in E$ with $\iota(e) = (s,s')$,  we need to distinguish between 
$\rho_{ee^{-1}} = \rho_{e^{-1}} \circ \rho_e = \rho_e^{-1} \circ \rho_e
  = \mathrm{id}_{\mathrm{dom}(\rho_e)}$ and $\rho_{\lambda_s} = \mathrm{id}_{s}= \mathrm{id}_{A_s}$.

\subsection{Measures of global consistency} 

\bD
\label{coherencedef}
An amalgamation pattern $\H/\I$ is called
\bre
\item
\emph{coherent} if for any two walks $w_1,w_2$ from $s$ to $s'$, 
the compositions $\rho_{w_1}$ and $\rho_{w_2}$
agree as maps from $A_s$ to $A_{s'}$ on the intersection of the domains.
\item
\emph{simple} if for each individual link $e \in E$, 
$\iota(e)=(s,s')$, the partial isomorphism $\rho_e \in
\mathrm{part}(\str{A}_s,\str{A}_{s'})$ extends 
every composition $\rho_w$ along walks $w$ from $s$ to $s'$ in $\I$,
in the sense that $\rho_e \supset \rho_w$.
\item
\emph{strongly coherent} if for any two walks $w_1, w_2$ 
from $s$ to $s'$ there is a walk 
$w$, also from $s$ to $s'$, such that $\rho_w$ is a common extension of
$\rho_{w_i}$ for $i=1,2$.
\ere
\eD

We note that coherence as defined in~(i) is 
equivalent to the condition that any composition $\rho_w$ along
a walk $w$ that loops at site $s$ is a restriction of the
identity $\rho_{\lambda_s} = \mathrm{id}_{A_s}$ on $A_s$. It is also
equivalent to the condition that the union of the 
bijections induced along any two walks between the same sites 
is again a bijection (but, unlike the case of strong coherence, not necessarily a partial isomorphism).

Strong coherence may be viewed as a confluence property for
compositions along different walks linking the same sites;
it in particular implies coherence. Simplicity will be of technical
interest later. It is clear that these notions impose non-trivial 
structural constraints on amalgamation patterns.

\bE
\label{IoverIex}
We may regard $\I$ itself as an amalgamation pattern $\I/\I$ (the
minimal one over $\I$) 
in the natural manner, with $S$ partitioned into singleton sets 
$\{ s\}$ and singleton maps 
$\rho_e \colon \iota_1(e) \mapsto \iota_2(e)$ for  $e \in E$.  
For any walk $w = e_1\cdots e_n$ from $s = \iota_1(e_1)$ to $s' = \iota_2(e_n)$
in $\I$, the composition $\rho_w$ precisely maps $s$ to $s'$. 
This amalgamation pattern satisfies coherence, simplicity and strong
coherence.

It is also instructive to check that all instances of the other most 
basic class of examples of amalgamation patterns, viz.\ exploded views according to 
Example~\ref{explodedviewex} and~\ref{explodedviewhypex},  
trivially satisfy all three conditions.
\eE

With an amalgamation pattern $\H/\I$ we associate the equivalence 
relation $\approx$ on $H$ that is induced by the $\rho_e$ if we regard 
them as identifications (in the sense of a prescribed overlap). I.e., we
let $\approx$ be the reflexive transitive closure of the
union of the graphs of the $\rho_e$ for $e \in E$. Then for $a \in
A_s$ and $a' \in A_{s'}$, 
\[
a \approx a' \; \mbox{ iff } \; a' = \rho_w(a) 
\mbox{ for some walk $w$ from $s$ to $s'$.} 
\] 

In the following we write $[a]$ 
for the equivalence class of $a \in H$ w.r.t.\
$\approx$:
\[
[a]  := \{ \rho_w(a) \colon \mbox{ $w$ a walk in $\I$, } 
a \in \mathrm{dom}(\rho_w) \}.
\] 

\bL
\label{cohapproxlem}
Coherence of $\H/\I$ guarantees that $\approx$ is trivial (coincides
with equality) in restriction to each $\str{A}_s$. 
Simplicity guarantees that, for any two sites $\str{A}_s$ and
$\str{A}_{s'}$ that are directly linked by some $e \in E$ 
with $\iota(e) = (s,s')$, $\approx$ identifies just the 
pairs linked by $\rho_e$. 
Strong coherence moreover implies that $\approx$ is a congruence
w.r.t.\ the interpretation of relations in $\str{H}$.
\eL

\prf
The first claim is an immediate consequence of the definition of
coherence, Definition~\ref{coherencedef},
and the characterisation of $\approx$ in terms of the action
of the $\rho_w$ on $H$ as given above. Similarly, the claim regarding simplicity
is obvious from the definition. For the third claim we observe
that, as a composition of partial isomorphisms $\rho_e$,  every 
$\rho_w$ is a partial isomorphism. While unions of two or more partial
isomorphisms may not be partial isomorphisms
(cf.~Example~\ref{moebiusex}), strong coherence in 
Definition~\ref{coherencedef} is designed to guarantee that every
union of partial isomorphisms $\rho_{w_i}$ along different walks 
that link the same pair $\str{A}_s$ and $\str{A}_{s'}$ admits a
common extension $\rho_{w}$. It follows that the union $\rho$ of all 
those partial isomorphisms linking $\str{A}_s$ and
$\str{A}_{s'}$ is itself a partial isomorphism of this kind, say
$\rho = \rho_w$. Any tuple $\abar$ in $\str{A}_s$ that is (component-wise)
$\approx$-equivalent with $\abar'$ in $\str{A}_{s'}$ thus is in the
domain of this maximal $\rho$, and compatibility of $\abar \approx
\abar'$ with the interpretations of relations $R$ in $\str{A}$ and
$\str{A}_{s'}$ follows since $\rho \in \mathrm{part}(\str{A}_s,\str{A}_{s'})$.
\eprf

The proof indicates that coherence deals with global consistency at the level of
elements, while strong coherence enforces global consistency at the level of
tuples, and therefore at the relational level. The straightforward
definition of global consistency, as an outright property of a given
amalgamation pattern $\H$ is the following.

\bD
\label{gobalconsdef}
The amalgamation pattern $\H/\I$ is \emph{globally consistent}
if 
\[
\str{A}:= \bigl(\dot{\bigcup}_{s \in S} \str{A}_s\bigr)/{\approx}
\] 
is a well-defined $\sigma$-structure $\str{A}$ that admits an atlas 
with co-ordinate structures $\str{A}_s$ and co-ordinate domains
$\str{A}_s/{\approx}$ with the natural charts 
$\pi_s \colon \str{A}_s/{\approx} \simeq \str{A}_s$ that associate
the $\approx$-equivalence class $[a]$ of $a \in A_s$ with $a$.
\eD

The following is clear from Lemma~\ref{cohapproxlem}.

\bR
\label{conscoherentrm}
Strong coherence implies global consistency,  and global consistency implies coherence.
\eR

\bE
\label{moebiusex}
Consider the following amalgamation pattern: $\str{A}_1$ and
$\str{A}_2$ are isomorphic copies of a two-element structure
consisting of a single directed edge ($\sigma = \{ R \}$, $R$ binary)
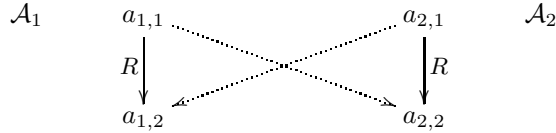
\[
\xymatrix{
\str{A}_1 & a_{1,1} \ar@{->}[d]_{\makebox(5,10){\small $R$}} \ar@{.>}[rrrd]
&&& a_{2,1} \ar@{->}[d]^{\makebox(5,10){\small $R$}}
\ar@{.>} [llld]
&\str{A}_2 
\\
& a_{1,2} &&& a_{2,2} &
}
\]

Two links identify the endpoints of these single edges in a cross-wise
fashion: $\rho_1$ matches the source of the $R$-edge of $\str{A}_1$ with the
target of the $R$-edge of $\str{A}_2$, $\rho_2$ likewise matches  
the source of the $R$-edge of $\str{A}_2$ with the
target of the $R$-edge of $\str{A}_1$ (dotted arrows in the diagram):
$\rho_1 \colon a_{1,1} \mapsto a_{2,2}$ 
and  $\rho_2 \colon a_{2,1} \mapsto a_{1,2}$.

Clearly, this amalgamation pattern is 
coherent, in fact trivially so, since the $\rho_i$
have no non-trivial compositions; it is neither strongly coherent, nor
is it globally consistent: $\rho_1 \cup \rho_2^{-1}$ is not a partial isomorphism.
The twist in this overlap pattern is reminiscent of the M\"obius strip:
it turns out that it, too, requires an at least two-fold covering to
overcome the global inconsistency in the overlap specification.

With links that are straight instead of twisted, $\rho_1' \colon a_{1,1} \mapsto a_{2,1}$ 
and  $\rho_2' \colon a_{2,2} \mapsto a_{1,2}$, the resulting amalgamation
pattern would be globally consistent (with $\str{A} = (\str{A}_1 \dot{\cup}
\str{A}_2)/{\approx}
\;\simeq
\str{A}_i$), though still not strongly coherent, since the combination of
the $\rho_i'$, albeit an isomorphism, is not realised as a composition of specified links.
\eE

\subsection{Realisations of amalgamation patterns} 

The idea of a \emph{realisation} 
of an amalgamation pattern strives for the converse of 
the passage from a structure to its exploded view. This
becomes interesting in the case where the pieces in the amalgamation
pattern do not already in themselves combine to form a single structure
--- i.e.\ when the local consistency built into the amalgamation
pattern fails to trivially add up to the global consistency of an
atlas as a global master plan according to Definition~\ref{gobalconsdef}. 
E.g.\ in the example reminiscent of a M\"obius strip, a
two-fold covering resolves the global inconsistency.

\bD
\label{realdef}
A \emph{realisation} of an amalgamation pattern $\H$ 
over $\I$ is a $\sigma$-structure $\str{A}$
together with an atlas 
\[
\A/\H = (\H/\I,\str{A},U,(U_s)_{s \in S}, (\pi_{u,s})_{u\in U_s}),
\]
such that $U_s \not= \emptyset$ for all $s \in S$ and,
for $u \in U_s$, $\pi_{u,s}$ is a chart from $\str{A}\restr u$ 
onto $\str{A}_s$, and the link structure of $\H/\I$ is reflected  
in the following tight manner: 
\bre
\item all links of $\H/\I$ are realised as overlaps locally: 
\\
for every $u \in U_s$ and $e \in E$ with $\iota(e) =(s,s')$
there is some $u' \in U_{s'}$ such that
$\mathrm{dom}(\rho_e) = \pi_{u,s}(u \cap u')$ and 
$\rho_e =  \pi_{u',s'} \circ \pi_{u,s}^{-1}$. 
\item there are no incidental overlaps, globally: 
\\
for any $u \in U_s, u' \in U_{s'}$ with $u \cap u' \not=\emptyset$, 
there is some walk $w$ in $\I$ such that 
$\mathrm{dom}(\rho_w) = \pi_{u,s}(u \cap u')$ and
$\rho_w = \pi_{u',s'} \circ \pi_{u,s}^{-1}$. 
\ere
Formally we think of $\A/\H$ as a multi-sorted structure with 
sorts for the elements of $\str{A}$, 
the sort $U = \bigcup_{s\in S} U_s \subset
\mathcal{P}(A)$ for the co-ordinate domains of its atlas, 
a sort $P  = \{ \pi_{u,s} \colon u \in U, s \in S \}$ for the
charts, augmenting the multi-sorted structure $\H/\I$ with
the underlying incidence pattern $\I$, whose first sort $S$ serves as
an index sort for the atlas in the above presentation as a labelled family.
\eD

\bD
\label{realsymmdef}
A \emph{symmetry} of a realisation $\A/\H$ is an automorphism of 
the multi-sorted structure that links the relational structure
$\str{A}$ on sort $A$ 
to the amalgamation pattern $\H/\I$ through its atlas based on the
sorts $U \subset \mathcal{P}(A)$ and $P = \bigl\{ \pi_{u,s} \colon u
\in U_s, s \in S \}$. Its \emph{$\H$-rigid symmetries} are those 
automorphisms that fix $\str{H}$ (and thus all of $\H/\I$ and in
particular $\I$) pointwise.
A realisation $\A/\H$ of the amalgamation pattern
$\H/\I$ is \emph{fully symmetric} over  $\H/\I$  if every symmetry of
$\H/\I$ extends to a symmetry of $\A/\H$, and if the 
automorphism group of its $\H$-rigid symmetries acts transitively on
$U_s$, for every $s \in S$.
\eD

In terms of overlaps between isomorphic copies of the distinguished
template structures $\str{A}_s$ of $\H$, condition~(i) in
Definition~\ref{realdef} is a richness 
condition. It guarantees an extension property: 
all specified overlaps occur wherever applicable. 
Condition~(ii) on the other hand guarantees a 
minimality property: just the specified overlaps, and trivially induced ones, do occur. 

Condition~(i), as an \emph{extension property}, 
says that (for suitable choices of $u' \in U_{s'}$ for
given $u \in U_s$ and $e \in E$) the following diagram
of partial maps commutes (including the information about these maps
as partial isomorphisms that are bijective w.r.t.\ the indicated 
domains and ranges, which implies that the relationship between 
$\str{A}\restr u \simeq \str{A}_s$ and $\str{A}\restr u' \simeq
\str{A}_{s'}$ is that of a disjoint amalgam): 
\[
\xymatrix{
& \str{A}\restr (u \cap u')
\ar[ld]_{\pi_{u,s}} \ar[rd]^{\pi_{u'\!\!,s'}}&
\\
\str{A}_s\restr \mathrm{dom}(\rho_e) \ar[rr]^{\rho_e} && \str{A}_{s'}
\restr \mathrm{image}(\rho_e)
}
\]  

Similarly, condition~(ii) says that (for suitable choices of $w = e_1
\cdots e_n$ for given $u \in U_s$ and $u' \in U_{s'}$)
the following diagram
of partial maps commutes (including the information about these maps
as partial isomorphisms that are bijective w.r.t.\ the indicated 
domains and ranges):
\[
\xymatrix{
& \str{A}\restr (u \cap u')
\ar[ld]_{\pi_{u,s}} \ar[rd]^{\pi_{u'\!\!,s'}}&
\\
\str{A}_s\restr \pi_{u,s}(u \cap u') \ar[rr]^{\rho_{w}} && \str{A}_{s'}
\restr \pi_{u',s'}(u \cap u')
}
\]  

Given the overlaps that have to be realised
according to condition~(i), also any chain of
compositions of $\rho_{e_i}$ along a walk $w = e_1\cdots e_n$ from 
$s$ to $s'$ in $\I$ must be realised as an actual sequence of overlaps
that induces some intersection of the form $u \cap u'$ provided the
composition $\rho_w$ is non-empty in $\H$.  
Condition~(ii), as a \emph{minimality} requirement, requires that,
conversely, any non-trivial
intersection $u \cap u'$ arises in precisely this manner.

Condition~(ii) goes beyond the requirement of global consistency, 
similar to the manner in which strong coherence goes beyond global
consistency (cf.\ Definition~\ref{gobalconsdef}
and discussion in Example~\ref{moebiusex}). Correspondingly, 
strong coherence together with simplicity guarantees that the quotient $\str{A} :=
(\dot{\bigcup}_{s \in S} A_s)/{\approx}$ induces a realisation, as stated
in Observation~\ref{strongcoherelqutientlem} below. 
Also compare Remark~\ref{caninfrealrem} below for a discussion of 
a canonical infinite realisation for any amalgamation pattern, which can intuitively
be obtained by unfolding $\H/\I$ and $\I$ itself 
in a tree-like fashion and gluing disjoint copies of the $\str{A}_s$ in $e$-related
locations with overlaps as prescribed by $\rho_e$.

\bO
\label{strongcoherelqutientlem}
If $\bigl( \str{H}, (\str{A}_s)_{s \in S}, (\rho_e)_{e \in   E}\bigr)$
is strongly conherent and simple, then its natural
quotient w.r.t.\ $\approx$, with equivalence classes $[a] \in H/{\approx}$, 
induces a realisation $\H/{\approx} $
based on the relational structure 
$\str{A} = (\dot{\bigcup}_{s \in S} \str{A}_s)/{\approx}$ 
and the atlas of charts with co-ordinate domains 
$u_s = \{ [a] \colon a \in A_s \}$ onto the $\str{A}_s = \str{H}\restr
A_s$.  
\eO

\subsection{Homomorphisms and coverings}

A \emph{homomorphism} 
\[
\pi \colon \H^{\ssc (2)} \!
\longrightarrow \,\H^{\ssc (1)} 
\]
between amalgamation patterns
$\H^{\ssc (i)} = \bigl( \I^{\ssc (i)}, \str{H}^{\ssc (i)}, (\str{A}^{\ssc (i)}_s)_{s \in
  S^{\ssc (i)}}, (\rho^{\ssc (i)}_e)_{e \in E^{\ssc (i)}}\bigr)$ 
over incidence patterns 
$\I^{\ssc (i)} =(S^{\ssc (i)},E^{\ssc (i)})$, for $i =2,1$, 
consists of compatible projection maps $\pi$ between the various sorts of
the $\H^{\ssc (i)}/\I^{\ssc (i)}$, whose incarnation over the domain $H^{\ssc (2)}$ restricts to
isomorphisms between the templates $\str{A}^{\ssc (2)}_s$ of $\H^{\ssc (2)}$ and 
their images $\str{A}^{\ssc (1)}_{\pi(s)}$ in $H^{\ssc (1)}$, with
commuting diagrams
\[
\nt\hspace{-2cm}
\xymatrix{
\H^{\ssc (2)}\!/\I^{\ssc (2)}
\ar[dd]|{\makebox(10,10){$\pi$}}
&
*++{H^{\ssc (2)}}
\ar[dd]|{\makebox(10,10){$\pi$}}
& *+++{\str{A}^{\ssc (2)}_s}
\ar[rr]^{\rho ^{\ssc (2)}_e} 
\ar[dd]|{\rotatebox{270}{$\sc \,\simeq\,$}}
&&
*+++{\str{A}^{\ssc (2)}_{s'}}
\ar[dd]|{\rotatebox{270}{$\sc \,\simeq\,$}}
&
*+++{\I^{\ssc (2)}}
\ar[dd]|{\makebox(10,10){$\pi$}}
&
*++++{E^{\ssc (2)}}
\ar[dd]
\ar[rr]^{\iota^{\ssc (2)}_i}
&&
*++++{S^{\ssc (2)}}
\ar[dd]
\\
\\
\H^{\ssc (1)}\!/\I^{\ssc (1)}
&
*+++{H^{\ssc (1)}}
& 
*+++{\str{A}^{\ssc (1)}_{\pi(s)}}
\ar[rr]^{\rho ^{\ssc (1)}_{\pi(e)}} 
&&
*+++{\str{A}^{\ssc (1)}_{\pi(s')}}
&
*+++{\I^{\ssc (1)}}
&
*++++{E^{\ssc (1)}}
\ar[rr]^{\iota^{\ssc (1)}_i}
&&
*++++{S^{\ssc (1)}}
}
\]

Note that this notion of a homomorphism is rather strict at the local
level, in that it requires local bijectivity.
As the restrictions of $\pi$ to the $\str{A}_s^{\ssc (2)}$ are
bijective, and as $\rho_{\pi(e)}^{\ssc (1)} =
(\pi\restr \str{A}_{s}^{\ssc (2)})^{-1} \circ 
\rho_{e}^{\ssc (2)} \circ \pi\restr \str{A}_{s'}^{\ssc (2)}$, 
the compositions $\rho_e^{\ssc (2)}$ along
walks of $\I^{\ssc (2)}$ translate directly into compositions of the 
$\rho_e^{\ssc (1)}$ along the image walks in $\I^{\ssc (1)}$, and
equalities like  $\rho_{w_1}^{\ssc (2)} = \rho_{w_2}^{\ssc (2)}$,
$\rho_{w}^{\ssc (2)}=\emptyset$, or inclusions like $\rho_{w}^{\ssc (2)} \subset
\mathrm{id}_s$ at the level of $\H^{\ssc (2)}/\I^{\ssc (2)}$ 
imply the analogous 
$\rho_{\pi(w_1)}^{\ssc (1)} = \rho_{\pi(w_2)}^{\ssc (1)}$, 
$\rho_{\pi(w)}^{\ssc (1)}=\emptyset$, or $\rho_{\pi(w)}^{\ssc (1)} \subset
\mathrm{id}_{\pi(s)}$ at the level of $\H^{\ssc (1)}/\I^{\ssc (1)}$, 
but not vice versa. 
One essential way in which $\H^{\ssc (2)}$ can deviate
from its homomorphic image in $\H^{\ssc (1)}$ is in terms of 
a potentially thinner link structure, which could involve 
walks in $\I^{\ssc (1)}$ that are not $\pi$-images 
of walks of $\I^{\ssc (2)}$, as well as in terms of 
a potentially richer branching in the link structure so that 
walks in $\I^{\ssc (1)}$ could arise as $\pi$-images of 
walks of $\I^{\ssc (2)}$ in multiple ways. The first of these
deviations is ruled out by the following definition of a covering, 
which requires the homomorphism to provide lifts of walks. The second 
kind of variance is unaffected since lifts need not be unique. In
fact, a unique lifting condition would lead to a notion of unbranched
covering, while our notion crucially is one of \emph{branched
  covering}  (cf.\ Section~\ref{hypcovsec} for a discussion of
how unbranched coverings are too restrictive for our purpose).%
\footnote{E.g., the exploded view of the 
$3$-uniform hypergraph corresponding to an apex over an $n$-cycle 
(a triangulation of an $n$-gon from its centre)
only admits trivial unbranched coverings by disjoint copies of itself.}

\bD
\label{amalgcovdef}
A \emph{covering} of an amalgamation pattern
$\H^{\ssc (1)}$over $\I^{\ssc (1)}$ by another amalgamation pattern
$\H^{\ssc (2)}$ over $\I^{\ssc (2)}$
is a surjective homomorphism 
\[
\pi \colon \H^{\ssc (2)} \!\longrightarrow\,
\H^{\ssc (1)},
\]
surjective at the level of every one of the sorts involved, 
with the natural lifting property w.r.t.\ links: 
for every $s \in S^{\ssc (2)}$ and every $e \in E^{\ssc (1)}$ 
s.t.\ $ \pi(s) = \iota^{\ssc (1)}(e)$ there exists some 
$e' \in E^{\ssc (2)}$ s.t.\ $\pi(e') = e$. 
\eD

Note that the lifting property in this definition is 
expressed entirely in terms of the action of $\pi$ on the incidence
pattern, $\pi \colon \I^{\ssc (2)} \rightarrow \I^{\ssc (1)}$ where it
stipulates that $\pi$ induces a \emph{bisimulation} (rather than just a 
surjective homomorphism). 

\bD
\label{covsymmdef}
A \emph{symmetry} of a covering $\pi \colon \H^{\ssc (2)} \longrightarrow
\H^{\ssc (1)}$ is an automorphism of the multi-sorted structure 
consisting of the $\H^{\ssc (i)}/\I^{\ssc (i)}$ and the maps
induced by $\pi$ between corresponding sorts;  it is 
\emph{$\H$-rigid} w.r.t.\ $\H = \H^{\ssc (1)}$
if it fixes the amalgamation pattern $\H^{\ssc (1)}/\I^{\ssc (1)}$
  pointwise.
The covering $\pi \colon \H^{\ssc (2)} \longrightarrow
\H^{\ssc (1)}$ is \emph{fully symmetric} (over its base $\H^{\ssc (1)}/\I^{\ssc (1)}$) 
if every symmetry of $\H^{\ssc (1)}/\I^{\ssc (1)}$ extends to a symmetry
of the covering and if the group of its $\H^{\ssc (1)}$-rigid
symmetries acts transitively on each of the sets  
$\pi^{-1}(\{s\}) \subset S^{\ssc (2)}$, for $s\in S^{\ssc (1)}$.
\eD

Note that the transitivity condition for $\H^{\ssc (1)}$-rigid
symmetries implies that these also act transitively on each of 
the fibres $\pi^{-1}(a)$ above individual
elements $a \in A^{\ssc (1)}_s \subset H^{\ssc (1)}$.
The notion of a realisation in Definition~\ref{realdef} specifies required and 
admitted overlaps between charts in an atlas for the desired global
structure just up to bisimulation. This fact is borne out in precise terms
in the following.

\bR
\label{coverealrem}
If $\pi \colon \H^{\ssc (2)}\! \longrightarrow\,
\H^{\ssc (1)}$ is a covering, then any realisation 
$\A^{\ssc (2)}/\H^{\ssc (2)}$ induces a realisation 
$\A^{\ssc (1)}/\H^{\ssc (1)}$, which is obtained 
by the natural composition of charts of $\A^{\ssc (2)}$ with the 
local bijections induced by $\pi$.
\eR

\prf
Let $\A^{\ssc (2)} = (\H^{\ssc (2)}/\I^{\ssc (2)},\str{A}, U, (U^{\ssc (2)}_s)_{s \in S^{\ssc (2)}},
(\pi_{u,s}^{\ssc (2)})_{u \in U_s^{\ssc (2)}})$ be a 
realisation of $\H^{\ssc (2)}$ and put 
$\A^{\ssc (1)} = (\H^{\ssc (2)}/\I^{\ssc (1)},\str{A}, U, (U^{\ssc (1)}_s)_{s \in S^{\ssc (1)}},
(\pi_{u,s}^{\ssc (1)})_{u \in U_s^{\ssc (1)}})$ where
\[
\pi_{u,s}^{\ssc (1)} := \pi \circ \pi_{u,\hat{s}}^{\ssc (2)} \; \mbox{
  for } \; u \in 
U^{\ssc (1)}_s:= \bigcup \{ U^{\ssc
  (2)}_{\hat{s}}  \colon \pi(\hat{s}) = s \}.
\]
Clearly the new charts $\pi_{u,s}^{\ssc (1)}$ locally map
$\str{A}\restr u$ isomorphically onto $\str{A}_s$ for $u \in U^{\ssc (1)}_s$.
One checks that conditions~(i) and~(ii) for realisations carry over as required.
\eprf

\subsection{Realisations from coverings} 

Remark~\ref{coverealrem} shows that the notion of a covering fits well with
the intuition that an amalgamation pattern specifies the
combinatorial structure of an atlas of the desired realisations just
up to bisimulation. In fact, it 
also offers a useful route to 
special realisations, which will guide our further investigation.

\bL
\label{coverreallem}
For an amalgamation pattern 
$\H = \bigl(\I,\str{H}, (\str{A}_s)_{s \in S}, (\rho_e)_{e \in E}\bigr)$, 
any covering $\pi \colon \hat{\H} \rightarrow \H$ by 
an amalgamation
pattern $\hat{\H} = \bigl( \hat{\I}, \hat{\str{H}}, 
(\str{A}_{\hat{s}}) _{\hat{s} \in \hat{S}}, 
(\rho_{\hat{e}})_{\hat{e} \in \hat{E}}\bigr)$ 
that is simple and strongly coherent, induces a realisation of $\H$.
This realisation $\A/\H$ is based on the quotient of $\hat{\str{H}}$ 
w.r.t.\ $\approx$, 
$\str{A} = \bigl(\dot{\bigcup}_{\hat{s} \in \hat{S}} \str{A}_{\hat{s}}\bigr)/{\approx}$.
\eL

\prf
Let $\pi \colon \hat{\H} \rightarrow \H$ be a covering,
$\hat{\H} = \bigl( \hat{H}, 
(\str{A}_{\hat{s}})_{\hat{s} \in \hat{S}}, 
(\rho_{\hat{e}})_{\hat{e} \in \hat{E}}\bigr)$ 
strongly coherent. By Remark~\ref{conscoherentrm}, $\hat{\H}$ is in particular globally
consistent so that 
$\str{A} = \dot{\bigcup} \str{A}_{\hat{s}}/{\approx}$ is well-defined
as a $\sigma$-structure on the universe $A = (\dot{\bigcup} 
A_{\hat{s}})/{\approx} = \{ [a] \colon a \in \bigcup
A_{\hat{s}} \}$. Let $U \subset \mathcal{P}(A)$ consist of the 
subsets $u_{\hat{s}} 
= \{ [a] \colon a \in A_{\hat{s}} \}$ for 
$\hat{s} \in \hat{S}$. We put $U_s := \{ u_{\hat{s}} \colon
\pi(\hat{s}) = s\}$ (noting that we may have $u_{\hat{s}_1} = u_{\hat{s}_2}$
for $\hat{s}_1 \not= \hat{s}_2$, due to $\approx$-identifications). 
By construction $\str{A}$ admits an atlas of charts 
$\pi_{\hat{s}}\colon \str{A}_{\hat{s}}/{\approx} \,\simeq\,
\str{A}_{\hat{s}}$, whose co-ordinate domains are these subsets
$u_{\hat{s}} 
\in U$. 
By composition of these $\pi_{\hat{s}}$ with the
covering projection $\pi$, we obtain charts 
\[
\pi_{u_{\hat{s}},\pi(\hat{s})} = \pi \circ \pi_{\hat{s}}
\;\colon \str{A}_{\hat{s}}/{\approx} \;\simeq\;
\str{A}_{\pi(\hat{s})}
\] 
from co-ordinate domains $u_{\hat{s}} \in U$ onto the co-ordinate
structures $\str{A}_s$, $s= \pi(\hat{s})$, of $\H$. 
We check the conditions from Definition~\ref{realdef} on realisations:
condition~(i) uses the lifting property for the covering $\pi$, for 
a traversal of a single link $e \in E$, which yields a link $\hat{e}
\in \hat{E}$ that is implemented as the actual, full overlap 
between corresponding sites $\str{A}_{\hat{s}}/{\approx}$ and
$\str{A}_{\hat{s}'}/{\approx}$ due to simplicity of $\hat{\H}$.
Condition~(ii) for $\A/\H$ directly corresponds to the strong
coherence requirement on $\hat{\H}$.
\eprf

It is instructive to see how a realisation 
$\A/\H = (\H,\str{A}, U, (U_s)_{s \in S}, (\pi_{u,s})_{u \in U_s})$
of $\H$ conversely also gives rise to a covering $\pi \colon \hat{\H} \rightarrow
\H$. A natural induced covering is based on 
the disjoint union of $u$-tagged copies of $\str{A}\restr u
\simeq \str{A}_s$ for $u \in U_s$ with some choice of 
$\hat{\I} = (\hat{S},\hat{E})$ with $\hat{S} = \{
(u,s) \colon u \in U_s \}$. The covering link set $\hat{E}$ can be
based on a choice of pairs $\hat{e} = ((u,s),(u',s'))$ 
over $\hat{S}$, for each $e \in E$ with $\iota(e) = (s,s')$, such that
$(u,s)$ for $u \in U_s$ is paired with $(u',s')$ for one choice of a $u' \in U_{s'}$
such that the identity on $u \cap u'$ precisely corresponds to 
$\pi_{u',s'} \circ \rho_e \circ \pi_{u,s}$, witnessing condition~(i)
for the realisation $\A/\H$. This covering will, however,
typically fail to be strongly coherent (rather than just simple and globally
consistent). This is because the
project-and-lift relationship between walks in $\hat{\I}$ and their
projections to $\I$ preserves equalities between source and target
sites only in the projection, not in the lifting: if $\hat{w}_1$ and
$\hat{w}_2$ both link $\hat{s} = (u,s)$ to  $\hat{s}' = (u',s')$ in
$\hat{\I}$, then $\pi(\hat{w}_i) = w_i$ 
both link $s$ to $s'$; but some $\rho_w$ for which
$\rho_w \supset \rho_{w_i}$, which exists according to condition~(ii) for
realisations, may lift to $\rho_{\hat{w}}$ for some $\hat{w}$ from $\hat{s}$ whose target
site could be some other site  $(u'',s') \in \pi^{-1}(s')$ rather than the
desired $\hat{s}' = (u',s')$.  This is in contrast with the
situation of Remark~\ref{coverealrem}, where condition~(ii) for realisations
is available at the level of the covering $\H^{\ssc (2)}$  (here $\hat{\H}$, above) rather
than at the level of the base structure $\H^{\ssc (1)}$  (here $\H$, below).

\section{Towards generic realisations}

By \emph{genericity} we mean to emphasise that the desired 
solutions to the realisation problem are fully compatible with
isomorphisms between given instances. In particular our constructions 
do not involve any ad-hoc choices that could break internal symmetries 
of the given instance. Apart from greater
mathematical elegance, this consideration will be of crucial 
importance for some of the applications to be
discussed later, especially in Section~\ref{HLsec}.

\subsection{An algebraic-combinatorial approach}

The term monoid is traditionally used to relax the
requirements on groups w.r.t.\ existence of inverses while the term 
groupoid, in its intended meaning here, relaxes the requirement that
the fundamental binary operation be total or that the underlying carrier
set be single-sorted.
The combination of both relaxations is captured in the notion of 
a \emph{category}. 
To emphasise the salient points for our considerations, we 
choose to call the resulting structures \emph{groupoidal monoids},
regarded as  multi-partite algebraic structures, which
may not necessarily provide inverses for their partial,
sort-dependent composition operation. 
The groupoidal monoids to be considered here reflect, as algebraic structures,
features of the underlying link structure $\I$ of an amalgamation pattern. 
Towards realisations we shall later want to lift them to 
proper groupoid structure with inverses. Groupoids should be thought of
as multi-partite analogues of groups, as well as categories with
invertible morphisms or as close counterparts of inverse semigroups, 
cf.~\cite{Lawson}.

\subsection{Groupoidal monoids and groupoids}
\label{groupoidmonsec}

Recall from Definition~\ref{incpatdef}
the format of an incidence pattern $\I = (S,E)$ in multigraph notation, 
which stands as shorthand for the more adequate multi-sorted
formalisation $\I  = (S,E,\iota, \cdot^{-1})$. 
Here $\iota = (\iota_1,\iota_2)$ stands for the pair of functions
that associate source and target vertices in $S$ with every edge $e \in
E$ in a manner that is compatible with edge reversal $e \mapsto e^{-1}$. 
A \emph{walk} in $\I$ is a finite sequence or word $w = e_1 \cdots
e_n$ for some $n \in \N$ (the length of the walk $w$) such that 
$\iota_2(e_i) = \iota_1(e_{i+1})$ for $i < n$. We write $E^\ast$ for
the set of all walks in $\I$. Special walks are those of length $n=0$:
there is precisely one such at each $s \in S$, corresponding to 
the empty word (at $s$) denoted $\lambda_s$ in the following.
The functions $\iota = (\iota_1,\iota_2) \colon E \rightarrow S$
naturally extend to all of $E^\ast$, with 
$\iota_1(e_1\cdots e_n) = \iota_1(e_1)$ and  
$\iota_2(e_1\cdots e_n) = \iota_2(e_n)$. In particular 
$\iota(\lambda_s) = (s,s)$.

It is also convenient to partition $E^\ast$ into the
subsets of walks from $s$ to $s'$ in $\I$, for all pairs $s,s' \in S$,
\[
E^\ast[s,s'] = \{ w \in E^\ast \colon
\iota(w) = (s,s') \},
\]
and to define the partial concatenation operation on $E^\ast$ 
as the union over
\[
\barr{rcl}
\cdot \colon E^\ast[s,t] \times E^\ast[t,u] &\longrightarrow & E^\ast[s,u] 
\\
(w,w') &\longmapsto & w \cdot w' := ww',
\earr
\]
which is precisely defined on all pairs $(w,w')$ such that 
$\iota_2(w) = \iota_1(w')$.

\bD
From the incidence pattern $\I = (S,E,\iota,\cdot^{-1})$ we 
derive the \emph{free groupoidal monoid} $\I^\ast$ 
as a many-sorted structure 
\[
\I^\ast = (E^\ast, S, \iota, \cdot, \{ \lambda_s \colon s \in S\}),
\]
where $E^\ast$ is the set of all walks in $\I$, 
with the natural extensions of $\iota$ from $E$ to $E^\ast$ as given above,
with the (partial) concatenation operation $\cdot$ for walks as above
and with the (empty) walks $\lambda_s$ of length $0$ at $s$ as
a set of distinguished elements, which are the units w.r.t.\ concatenation.
\eD

Note that $\I^\ast$ is groupoidal in that it has a sort-restricted
concatenation operation, which is total in restriction to matching
sorts, and that it is monoidal in not having inverses for elements
other than the units.

For notational convenience we also use the shorthand 
\[
\I^\ast = \bigl( E^\ast,
(E^\ast[s,s'])_{s,s' \in S}, \cdot, (\lambda_s)_{s \in S}\bigr)
\] 
in close analogy with the shorthand 
$\I = (S,E)$ for $\I = (S,E,\iota,\cdot^{-1})$.
W.r.t.\ these shorthands we keep in mind that they do not 
reflect symmetries or automorphisms
appropriately in as far as they might wrongly suggest a rigid labelling of/by
elements of the sets $S$ and $E$.

The following idea of a \emph{reduced product} between an amalgamation
pattern $\H/\I$ and $\I^\ast$ yields a particularly natural, 
tree-like realisation of $\H/\I$. Albeit infinite, it illustrates and
motivates key features of our construction of finite
realisations. Indeed, the desired finite realisations can be pictured 
as suitable quotients of this canonical infinite realisation. 
Much of the remainder of this section
focuses on the combinatorics of 
finite algebraic derivatives of $\I$ and $\I^\ast$ in relation to
$\H/\I$ towards the construction of reduced products, ultimately with 
suitable finite groupoids over $\I$ in Section~\ref{redprodrealsec}. 

The \emph{reduced product} $\H\otimes \I^\ast$ between $\H/\I$ and $\I^\ast$
is based on the quotient of the relational structure consisting of the disjoint union 
\[
\dot{\bigcup}_{w \in E^\ast} \str{A}_{\iota_2(w)} \times \{w\},
\]
w.r.t.\ to the congruence relation $\approx$ induced as the reflexive and
symmetric transitive closure of all identifications according to 
\[
(a,w) \approx (\rho_e(a),we)
\]
for $\iota_2(w) = \iota_1(e)$.
It follows that $(a_1,w_1) \approx (a_2,w_2)$ if, and only if,  there are
some $w_0,w_1',w_2' \in E^\ast$ such that 
$w_i = w_0 w_i'$ for $i=1,2$ and $a_2 = \rho_{w_2'}
(\rho_{w_1'}^{-1} (a_1))$. 
Writing $[(a,w)]$ for the $\approx$-equivalence class of $(a,w)$,
the natural atlas is based on the co-ordinate domains 
\[
u[w] := \{ [(a,w)] \colon a \in A_{s} \}
\] 
for $s = \iota_2(w)$. In
$U_s := \{ u[w] \colon \iota_2(w) = s \}$ we then have 
natural charts onto $\str{A}_s$.
Choosing $w_0$ to be the maximal common prefix of $w_1$ and $w_2$ 
in the above representation $w_i = w_0 w_i'$ for $i=1,2$, we find 
the exact overlap between $u[w_1]$ and $u[w_2]$ 
accounted for by $\rho_{w_2'} \circ \rho_{w_1'}^{-1} = \rho_w$ for
$w = w_1'\nt^{-1} w_2'$. The reduced product 
$\H\otimes\I^\ast/{\approx}$
can be pictured as a forest of
trees of overlapping copies of the $\str{A}_s$. 
The following claim is then straightforward.

\bR
\label{caninfrealrem}
The reduced product $\H\otimes \I^\ast/{\approx}$ 
induces a realisation of $\H/\I$. 
\eR

We note that $\H\otimes\I^\ast/{\approx}$ is infinite if $\H$ has at
least one non-trivial link $\rho_e$, for which $\mathrm{image}(\rho_e)
\strictsubset A_{\iota_2(e)}$: for $a \in  A_{\iota_2(e)} \setminus
\mathrm{image}(\rho_e)$ and $w_i = (e^{-1}e)^i$ ($i$-fold
concatenation in $\I^\ast$), the elements
$(a,w_i)$ will be pairwise
inequivalent w.r.t.\ $\approx$.

\paragraph*{Inverse semigroup structure.}
To put the natural action of $\I^\ast$ on the amalgamation pattern
$\H/\I$ on a more formal footing
we regard the inverse semigroup generated by the $\rho_e$ as an inverse
sub-semigroup of the full inverse semigroup of all partial bijections
on the universe $H$ of $\H$. The latter is known as the 
\emph{symmetric inverse semigroup} on that set, in analogy with the 
\emph{symmetric group} of all global bijections, cf.~\cite{Lawson}.

\bD
\label{symminvsemdef}
The \emph{symmetric inverse semigroup} associated with a set $X$ is
the inverse semigroup $I(X) = \bigl( \mathrm{part}(X,X), \cdot \bigr)$
consisting of all partial bijections of $X$ with the binary operation
of (partial) composition: for $\rho_i \in \mathrm{part}(X,X)$, the
product $\rho_1 \cdot  \rho_2$ (typically also denoted $\rho_2 \circ
\rho_1$) is the composition whose domain is $\mathrm{dom}(\rho_1 \cdot \rho_2)
= \mathrm{dom}(\rho_1) \cap \rho_1^{-1} (\mathrm{dom}(\rho_2)) =
\rho_1^{-1} (\mathrm{image}(\rho_1) \cap \mathrm{dom}(\rho_2))$.  
\eD

Idempotents $1_u = \mathrm{id}_u$ (local identities) 
for $u \subset X$ serve as local left or right neutral elements w.r.t.\ 
composition with elements $\rho \in I(X)$ whose domain or range is
contained in $u$:
$ 1_u \cdot \rho = \rho$ if $\mathrm{dom}(\rho) \subset u$ and 
$\rho \cdot 1_u = \rho$ if $\mathrm{image}(\rho) \subset u$.
And conversely, the $1_u$ precisely arise as compositions 
$\rho \cdot \rho^{-1}$ and $\rho^{-1} \cdot \rho$ for those 
$\rho \in I(X)$ whose domain or range is precisely $u$, respectively.

In the case of $I(H)$ where $H$ is the universe of some amalgamation pattern 
$\H = (H,(\str{A}_s)_{s \in S}, (\rho_e)_{e \in E})$, we shall also 
write $\mathrm{id}_s$ instead of  $1_{A_s} = \mathrm{id}_{A_s}$ for the local 
identities related to the partition of $H$ into the universes 
$A_s$ of the $\str{A}_s$ for $s \in S$. In this context, we 
may consider the partial isomorphisms $\rho_e$ as
elements of $I(H)$, the full symmetric semigroup on the set $H$.
The inverse semigroup composition of elements $\rho_{e_1}$
and $\rho_{e_2}$ is defined for any such pair, regardless of the
sort-typing that is implied by the incidence pattern $\I$, or the
monoidal structure of $\I^\ast$. For $w \in E^\ast$, i.e.\ for
a walk $w = e_1\cdots e_n$ in $\I$, the associated 
inverse semigroup product is 
\[
\rho_w = \prod_{i=1}^n \rho_{e_i} = \rho_{e_1} \cdots \rho_{e_n}
= \rho_{e_n} \circ \cdots \circ \rho_{e_1}
\]
in the sense of partial composition of maps, in agreement with
the monoidal product of $\I^\ast$. In this case the 
$\iota$-values encode source and target sites 
that accordingly determine inclusions of the form $\mathrm{dom}(\rho_w) \subset
A_s$ for $s= \iota_1(w)$ and $\mathrm{image}(\rho_w) \subset A_{s'}$
for $s' = \iota_2(w)$. But the inverse semigroup structure of $I(H)$
allows arbitrary products $\prod_{i=1}^n \rho_{e_i}$ regardless of
$\iota$-values.  On one hand it partly reflects the above inclusion 
constraints in equations like $\rho_w = \mathrm{id}_s \cdot \rho_w$ 
and $\rho_w = \rho_w \cdot \mathrm{id}_{s'}$; on the other hand, this reflection
is just partial because all such equations trivialise for $\rho_w =
\emptyset$
(which may happen for products along walks in $\I$ that are good in
the sense of $\I^\ast$, and is always the case for products that are 
not aligned along walks in
$\I$ and hence not represented in $\I^\ast$ at all).
In particular $\prod_{i=1}^n \rho_{e_i} = \emptyset$  for 
$e_1 \cdots e_n \not\in E^\ast$.

\bD
\label{Hinvsemdef}
The inverse semigroup $I(\H) \subset I(H)$ is generated as a
sub-semigroup of $I(H)$ by 
the $(\rho_e)_{e \in E}$ of 
$\H = (H,(\str{A}_s)_{s \in S}, (\rho_e)_{e \in E})$ together with 
the local identities $\mathrm{id}_s := \mathrm{id}_{A_s}$.%
\footnote{Note that
local identities $\mathrm{id}_s$ will not necessarily arise as
non-trivial products of $\rho_e$ with their inverses, since all $\rho_e$ with $\iota_1(e) = s$
may have $\mathrm{dom}(\rho_e) \strictsubset A_s$.}
\eD

Since we regard the local identities $\mathrm{id}_s$ as induced by walks 
$\lambda_s$ of length~$0$ at $s \in S$, at least all the
\emph{non-empty} elements are partial isomorphisms 
of the form $\rho_w = \prod_{i=1}^n \rho_{e_i}$ for walks $w$ in $\I$, i.e.\
for $w \in E^\ast$. 
If we restrict to just those elements of
$I(\H)$ that are induced by compositions $\rho_w$ of $\rho_e$ along
walks $w$ in $\I$, we move away from the globally defined composition
operation of an inverse semigroup.
 
\paragraph*{Groupoid structure.}
We now want to look at groupoids whose partial and sort-sensitive operation of 
full rather than partial composition reflects the incidence pattern 
$\I$ as the underlying combinatorial schema of sorts. 
Groupoids can also be described as categories all of whose morphisms 
are isomorphisms. Fixing an underlying  incidence pattern
$\I$ with sorts $S$ and $E$ means that we
fix the object sort of the category to be the set $S$ (of sites or sorts)
and a set of generators $e \in E$ for the morphism sort that link
the objects of sort $s$ and $s'$ as prescribed by $\iota(e) = (s,s')$.
In these terms, the following definition puts the focus on the
algebraic composition structure for the morphisms.

\bD
\label{Igroupoiddedf}
A \emph{groupoid} $\G$ over the incidence pattern $\I$, denoted by the
shorthand $\G/\I$, is a multi-sorted structure of the form 
\[
\G/\I = \bigl(\I, G, \iota \colon G \rightarrow S^2, \cdot \bigr)
\]
whose set $G$ of groupoid elements is partitioned by $\iota$ into
the sets $G[s,s'] := \iota^{-1}(s,s') = \{ g \in G \colon \iota(g) =
(s,s') \}$ of groupoid elements of source/target sorts $s/s'$, with a 
partial composition operation defined on pairs of matching interface
sort according to 
\[
\barr{rcl}
\cdot \colon G[s,s'] \times G[s',s''] &\longrightarrow & G[s,s'']
\\
(g_1, g_2) &\longmapsto& g_1 \cdot  g_2 =: g_1 g_2
\earr
\]
that is associative and has, for all $s,s' \in S$,
\bre
\item
unique left and right neutral elements $1_s \in G[s,s]$, s.t.
\\
$1_s g = g$ for all $g \in G[s,s']$ and 
$g 1_s = g$ for all $g \in G[s',s]$;
\item
unique inverses $g^{-1} \in G[s',s]$ for $g \in G[s,s']$,
s.t.
\\
$g g^{-1} = 1_s$ and $g^{-1} g = 1_{s'}$.
\ere 

In addition to these basic algebraic conditions we require 
$G$ to be generated by $E$ (the
link sort of $\I$), via a map $e \mapsto e^\G$ that associates a
groupoid element $e^\G \in G[\iota_1(e),\iota_2(e)]$ with every $e \in
E$, such that $(e^{-1})^\G = (e^\G)^{-1}$, and such that
\bre
\addtocounter{enumi}{2}
\item
every $g \in G$ can be written 
as a product 
$g = \prod_{i=1}^n e_i^\G =: w^\G$
for some walk $w = e_1 \cdots e_n \in E^\ast$ in $\I$.
(\footnote{The length~$0$ walk $\lambda_s$ at $s$ is taken to
generate the (empty) product $1_s$ of sort $\iota(1_s) = (s,s)$.
Note that $w^\G \in G[s,s']$ if $w$ is a walk from $s$ to $s'$ in
$\I$, i.e.\ $\iota(w^\G) = \iota(w)$.})
\ere
\eD

It is convenient to use the enumeration of relevant ingredients as shorthand:
\[
\G = \bigl(
G, (G[s,s'])_{s,s' \in S}, \cdot , (1_s)_{s \in S}, (e^\G)_{e \in E}
\bigr).
\]

We also use the suggestive abbreviations
\[
G[\ast,s]:= \bigcup_{s' \in S} G[s',s]\; \mbox{ and } \; 
G[s,\ast] := \bigcup_{s' \in S} G[s,s'],
\]
so that, for instance, 
$\cdot$ is defined
precisely on $\bigcup_{s \in S} (G[\ast,s] \times G[s,\ast])$.

\bD
\label{simplegroupoiddef}
$\G/\I$ is called \emph{simple} if the generators $(e^\G)_{e \in E}$
are pairwise distinct and $e^\G \not= 1_s$ for all $e\in E, s \in S$.%
\footnote{$e^\G \not= 1_s$ is meaningful as a constraint just for reflexive 
links $e \in E[s,s]$ (loops at $s$).} 
\eD

For a simple groupoid $\G/\I$ we identify $e$ with $e^\G$ and regard 
the link sort $E$ of $\I$ as a subset of $G$.

\bE
\label{IoverIgroupoidex}
Looking at $\I/\I$ as an amalgamation pattern as in Example~\ref{IoverIex},
and at the action of $\I^\ast$ on this particular amalgamation pattern, 
we can isolate a natural groupoid $\G(\I)/\I$ within the semigroup 
$I(\I) \subset I(S)$: it consists of those $\rho_w \in I(\I)$ that are
induced by walks 
$w \in \I^\ast$, with $\iota(\rho_w) = \iota(w)$ and with 
the natural composition. 
Since these $\rho_w$ uniformly are just the singleton maps $\rho_w \colon 
\iota_1(w) \mapsto \iota_2(w)$,  their composition in matching
interface sorts is exact rather than partial, and induces a groupoid
operation with units $1_s = \mathrm{id}_{\{s\}} = \rho_{\lambda_s}$. 
This groupoid $\G(\I)/\I$ is simple if, and only if $\I = (S,E)$ is a simple graph
(rather than a multi-graph, possibly with loops).
\eE

\bE
\label{completemalgex}
Call an amalgamation pattern $\H/\I$ \emph{complete} if all the 
$\rho_e$ are bijections between the sites involved: $\rho_e \colon
A_{\iota_1(e)} \rightarrow A_{\iota_2(e)}$ is a bijection for all $e \in E$. 
In this case, too, the groupoidal action of $\I^\ast$ induces 
a natural groupoid $\G(\H)/\I$ within the semigroup 
$I(\H) \subset I(H)$. Again, $\G(\H)$ consists of the $\rho_w \in I(H)$ induced
by walks $w \in \I^\ast$, with the natural composition so that, due to completeness,
every $\rho_w$ is a bijection between
$\mathrm{dom}(\rho_w) = A_{\iota_1(w)}$ and 
$\mathrm{image}(\rho_w) = A_{\iota_2(w)}$.
\eE

\bR
\label{freegroupoidrem}
The free groupoidal monoid $\I^\ast$ induces a simple groupoid 
obtained as the quotient of $\I^\ast$ w.r.t.\ cancellation of all
factors $e e^{-1}$, so that 
$w w^{-1}$ is identified with $\lambda_{\iota_1(w)}$, and $w^{-1}$ 
inverse to $w$. We regard this groupoid, which is a
simple groupoid over $\I$ and infinite for non-trivial $\I$, as the 
\emph{free groupoid over $\I$}.   
The same groupoid up to isomorphism is obtained in
the manner of Example~\ref{IoverIgroupoidex} from the natural tree unfolding of the
multi-graph $\I$.
\eR

Towards the construction of realisations we shall be interested in 
groupoids over $\I$ that respect those algebraic identities that are
induced by the inverse semigroup action of $\I^\ast$ 
on a given amalgamation pattern over $\I$, as
expressed in $I(\H)$ (cf.\ Definition~\ref{Hinvsemdef}). The
corresponding notion of \emph{compatibility} is defined as follows.

\bD
\label{compatdef}
A groupoid $\G/\I$ is \emph{compatible} with the amalgamation pattern 
$\H/\I$, both over the same incidence pattern $\I$, if, for all $w \in
\I^\ast$,
\[
w^\G = 1_s 
\mbox{ (in $\G$) } \; \Rightarrow \; \rho_w \subset \mathrm{id}_s =
\mathrm{id}_{A_s}
\mbox{ (in $I(\H)$). } 
\]
\eD

The following shows how compatibility in this sense can be achieved 
as a first step towards the construction of groupoids related to a
given amalgamation task. 

\bE
\label{compatex}
Consider an amalgamation pattern $\H/\I$ whose template structures
$\str{A}_s$ all have the same size. (This is w.l.o.g.\ for our
considerations, since the universes of the $\str{A}_s$ can be blown up
to any common size by trivial extensions.) Let $\hat{\H}/\I$ be any
\emph{completion} $\hat{\H}$ obtained by extension of the partial bijections
$\rho_e\colon \str{A}_{\iota_1(e)} \rightarrow \str{A}_{\iota_2(e)}$
to full bijections $\hat{\rho}_e$; these $\hat{\rho}$ are no longer
required to be partial isomorphisms between the $\sigma$-structures
$\str{A}_s$ of $\H$, and we think of $\hat{\H}$ as an amalgamation
pattern of plain sets with $\hat{\sigma} = \emptyset$.
Then the $\I$-groupoid $\G(\hat{\H})$ generated by the natural
action of these $\hat{\rho}_e$ on the set $H=\bigcup_{s \in S} A_s$ as
in Example~\ref{completemalgex} is compatible with $\hat{\H}$ and hence with $\H$.
The choice of extensions
$\hat{\rho}_e \supset \rho_e$ may break symmetries of $\H/\I$, but
this is no longer the case if $\G$ is similarly based on the action on the 
disjoint union of \emph{all} such completions $\hat{\H}$. 
\eE

\bL
\label{cohercompatlem} 
For any amalgamation pattern $\H/\I$, the following conditions 
are equivalent 
\bre
\item
$\H/\I$ is coherent;
\item
\emph{every} groupoid $\G/\I$ over the same $\I$ is compatible with
$\H/\I$;
\item
the groupoid $\G(\I)/\I$ from Example~\ref{IoverIgroupoidex} is
compatible with $\H/\I$.
\ere
\eL

\prf
By Definition~\ref{coherencedef},  
$\H/\I$ is coherent if any $\rho_w$ induced by a walk in $\I$ 
that loops at site $s$ satisfies $\rho_w \subset \mathrm{id}_s$ in $\H$. 
If $w^\G = 1_s$ in some $\G/\I$, then $\iota(w) = (s,s)$ implies that
$w$ is a walk that loops at $s$,
whence coherence implies that $\rho_w \subset \mathrm{id}_s$ as required for 
compatibility. Conversely, compatibility of the groupoid 
$\G(\I)/\I$ 
with $\H/\I$ implies that $\H/\I$ is coherent: in this groupoid, 
$w^{\G(\I)} = 1_s$ 
precisely for those walks $w$ in $\I$ with $\iota(w) = (s,s)$; so
compatibility implies that $\rho_w \subset \mathrm{id}_s$ in $\H$ for all those, as
required for coherence. 
\eprf

\subsection{Cayley graphs of groupoids}
 \label{groupoidsec}

With a groupoid $\G/\I$ we associate its \emph{Cayley graph}. 
First and foremost, this is the natural generalisation to the setting
of groupoids of the notion of a Cayley graph for
groups w.r.t.\ to a given set of generators. Moreover, 
it provides us with a crucial link between amalgamation patterns $\H/\I$ 
and groupoids $\G/\I$: the Cayley graph of $\G$
is an amalgamation pattern, indeed a very special, 
homogeneous and complete 
amalgamation pattern over the same underlying incidence pattern $\I$.

A group $G$, with generator set $E \subset G$ that is closed under 
inverses, can be seen as a special
case of a groupoid, based on an incidence pattern $\I = (S,E)$ with a
singleton set $S = \{ 0 \}$ of sites, and reflexive links $e \in E$
with trivial $\iota(e) = (0,0)$. This format of $\I$ implies that
$\cdot$ is total and indeed $(G,\cdot)$ a group; conversely, any
group with generator set $E = E^{-1}$ can be cast in this format. In
this sense, the following definition subsumes the familiar notion of
the Cayley graph of a group as a special case.

\bD
\label{Cayleydef}
With a groupoid $\G/\I$ over the incidence pattern $\I$ associate its
\emph{Cayley graph}, which is the many-sorted
(multi-)graph 
\[
\G/\I =
\bigl(
\I, G, \delta \colon G \rightarrow S, R, \eta \colon R \rightarrow E
\bigr),
\]
with vertex set $G = \dot{\bigcup}_{s \in S} G[\ast,s]$,
vertex-labelled by $S$ w.r.t.\ this partition, which is 
formalised by the partition function $\delta$, 
and with $R = \{ R_e \colon e \in E \}$ as a collection of edge relations $R_e$,
$E$-partitioned by $\eta$, with 
$e$-labelled edge sets 
\[
\textstyle
R_e = \{ e[g] \colon g \in G[\ast,\iota_1(e)] \}
\; \mbox{ with induced } \;
\barr[t]{@{}rcl@{}}
\iota \colon \bigcup R
& \rightarrow & G^2
\\
e[g] &\mapsto&  (g,ge^\G).
\earr
\]
\eD

The edge-labelling in this multi-graph suggests, 
as one choice for a natural relational encoding, to put
$e[g] := (g,e,ge^\G)\in R_e$, while for simple groupoids $\G$ 
even just $e[g] := (g,ge^\G) \in R_e \subset G^2$ is
good enough to determine the label. 

\bO
\label{Cayleygrapamalgpatternobs}
The Cayley graph of $\G/\I$ has the format of an amalgamation
pattern, where $\rho_e$ is the partial bijection induced by right
multiplication with the 
generator $e^\G$, whose graph is $R_e$. Note that $\rho_e$
bijectively relates the associated sites $G[\ast,s]$ and $G[\ast,s']$ if
$\iota(e)=(s,s')$. This amalgamation pattern is \emph{complete} in the sense
of Example~\ref{completemalgex}. Indeed the groupoid $\G$ can be retrieved
as the groupoid induced by this amalgamation pattern in the manner discussed
in that example. 
\eO

The simplified, labelled representation of the
underlying amalgamation pattern is 
\[
\bigl( G, (G[\ast,s])_{s \in S}, (\rho_e)_{e \in E} \bigr).
\]

Note that inversion corresponds to edge reversal,
$R_{e^{-1}} = R_e^{-1}$, and that walks in $\I$ represent groupoid
multiplication by groupoid elements as represented as products of
generators.

In the following the notation $\G/\I$ can refer to $\G/\I$ in these
interchangeable contexts: the groupoid, its Cayley graph or the
associated amalgamation pattern. Differences arise especially 
w.r.t.\ symmetries, to be discussed next.

\paragraph{Symmetries.}
Similar to Definition~\ref{amalgautdef} for amalgamation patterns, 
we also use the term \emph{symmetry} in connection with groupoids and
their Cayley graphs 
to refer to the adequate notion of automorphisms of the multi-sorted 
structures that allows for the parallel permutation of all the sorts
involved, rather than fixing, say, the labelling of sites, links, or
generators. We recall that it is largely to this end that we have,
e.g., officially
made the underlying incidence pattern 
$\I$ an internal component of the structure of a groupoid \emph{over} $\I$.

\bD
\label{innersymdef}
A \emph{symmetry} of a groupoid
$\G/\I = \bigl(\I, G, \iota, \cdot \bigr)$
over the incidence pattern $\I = (S, E, \iota, \cdot^{-1})$
is an automorphism $\pi$ of this multi-sorted structure. 
It is specified in terms of permutations $\pi = (\pi_S,\pi_E,\pi_G)$ 
of the three sorts $S$, $E$ and $G$ that 
bijectively map each sort within itself such that 
the pair $(\pi_S,\pi_E)$ is an automorphism of $\I$ and the triple
$(\pi_S,\pi_E,\pi_G)$ is compatible with 
$\iota^\G = (\iota_1,\iota_2) \colon G \rightarrow S^2$ and 
with the groupoid operation $\cdot$ on $G$. 

The groupoid $\G/\I$ is \emph{fully symmetric over $\I$}
if every symmetry of $\I$ extends to a symmetry of $\G/\I$. 

Symmetries of the Cayley graph of $\G/\I$ are defined
in agreement with the notion of symmetries of amalgamation patterns 
in Definition~\ref{amalgautdef}. Correspondingly, 
a symmetry of the Cayley graph of $\G/\I$ is called 
\emph{$\I$-rigid} if it fixes $\I$ pointwise.
\eD

We may think of a symmetry of the groupoid $\G/\I$ either 
as an extension of an automorphism of the two-sorted structure $\I$,
or as an automorphism of the algebraic groupoid structure $\G$ that 
also permutes generators in a manner compatible with the structure of
$\I$, i.e. inducing an automorphism of $\I$. 
  
Despite the close relationship, there are 
crucial differences between the grou\-poid on one hand and 
its Cayley graph or amalgamation pattern on the other.
These differences are particularly apparent in relation to their symmetries. 
In particular, we see that the groupoid operation $\cdot$ cannot
be definable \emph{within} the Cayley graph, although $\G$ is isomorphic
to the groupoid generated from the groupoidal action of $\I^\ast$ on
this Cayley graph according to Example~\ref{completemalgex} --- which
means that the abstract algebraic structure of the groupoid is fully determined
by the Cayley graph.
Definability of the groupoid operation 
within the Cayley graph is ruled out by 
richer symmetries of the latter.%
\footnote{The relationship is similar to that between an
  affine space and its associated vector space: the latter is
  definable \emph{from} the former, as the linear space of translations; it
  fails to be definable \emph{within} the former because there is no 
canonical choice of an origin.} 
In fact, the group of 
automorphism of the Cayley graph acts transitively on each partition set $G[\ast,s]$:
any pair $g,g'$ from the same 
partition set $G[\ast,s]$ are related by an automorphism 
that is induced by left multiplication with 
$g' g^{-1}$ on $G[\iota_1(g),\ast]$, and trivial elsewhere: 
\[
\barr{rcll}
G &\longrightarrow& G
\\
h &\longmapsto &
\left\{ 
\barr{ll}
g' g^{-1} h & \mbox{for } h \in G[\iota_1(g),\ast]
\\
h & \mbox{else.} 
\earr\right.
\earr
\]

This bijectively maps $G[\iota_1(g),s']$ onto $G[\iota_1(g'),s']$,
for every $s'$, and in particular maps $g$ to $g'$. Since it
preserves both the vertex partition and the edge
relations $R_e$, it is an automorphism of the Cayley graph, 
in fact an $\I$-rigid one. 

\bR
\label{transsymmrem}
The group of $\I$-rigid symmetries of the Cayley graph of $\G/\I$ 
acts transitively on the subsets $G[\ast,s]$ for every $s \in S$. 
\eR

The $\I$-rigid automorphism group of the
groupoid, on the other hand, is trivial, since it preserves each
$G[s',s]$ as a set and in particular fixes $1_s \in G[s,s]$ and
all generator elements $e^\G$ for $e \in E$. 
All automorphisms of the
Cayley groupoid, not just the $\I$-rigid ones, must preserve 
the set $\{ 1_s \colon s \in S \}$ of units.  

\paragraph*{Measures of acyclicity.}
As an amalgamation pattern, the Cayley graph of a groupoid $\G/\I$
need in general neither be simple nor coherent, let alone
strongly coherent. Due to completeness, coherence and simplicity fall into one;
and both conditions precisely mean that the groupoid is degenerate
in the sense that $G[s,s] = \{ 1_s \}$ for all $s$.
Corresponding criteria for groupoids and products of amalgamation
patterns with groupoids will be studied below. These
considerations are closely related to the following notion of
qualified \emph{coset acyclicity} for groupoids, which generalises
corresponding notions for Cayley groups from~\cite{OttoAPAL04}. Coset  
acyclicity puts much more severe restrictions than the familiar
notion of \emph{large girth}: while girth concerns the length of
cycles formed by generators (i.e.\ of graph-theoretic cycles in the
Cayley graph), coset acyclicity concerns the length of cycles formed
by overlapping cosets w.r.t.\ subgroup(oid)s generated by subsets of
the set of generators (i.e.\ of hypergraph-theoretic cycles in a 
dual hypergraph associated with the Cayley graph).

\bD
\label{subgroupoiddef}
In a groupoid $\G/\I$, a subset $\alpha \subset E$ that is
closed under inversion ($\alpha = \alpha^{-1}$) induces a 
\emph{sub-groupoid $\G[\alpha] \subset \G$ generated by $\alpha$}. 
$\G[\alpha]/\I$ consists of groupoid elements in
$G[\alpha] := \{ w^\G \colon  w \in \I^\ast \restr \alpha \} \subset
G$, with groupoid operation induced by $\G$. 
We write $G[\alpha,s,s']$ for the set of elements
$g \in G[\alpha]$ with $\iota(g) = (s,s')$, which is the set of
those $g \in G[s,s']$ that admit a representation as a product of
generators in $\alpha$, including empty products of the form
$\lambda_s^\G = 1_s$ for all $s \in S$.%
\footnote{Due to the trivial nature of the extra components,
  $\G[\alpha]$ could also be cast as a groupoid over 
the restricted incidence pattern 
$\I\restr\alpha$ over $S\restr \alpha =\iota_1(\alpha) =
\iota_2(\alpha)$.}
The (left) \emph{$\alpha$-coset}  at $g \in G[\ast,s]$ 
is 
\[
g \G[\alpha] = \{ gh \colon
h \in \G[\alpha], \iota_1(h)= s \} =
\{ g h \colon h \in \G[\alpha,s,\ast] \}. 
\]
\eD

\bD
\label{cosetcycdef}
For $n \geq 2$, a \emph{coset cycle} of length $n$ in the groupoid $\G/\I$, is
a cyclically indexed tuple of \emph{pointed cosets} $\bigl( g_i,  g_i \G[\alpha_i]
\bigr)_{i \in \Z_n}$ such that, for all~$i$:
\bre
\item
$g_{i+1} \in g_i \G[\alpha_i]$;
\item
$g_i \G[\alpha_i\cap \alpha_{i-1}] \cap g_{i+1} \G[\alpha_i\cap
\alpha_{i+1}] = \emptyset$.
\ere
\eD

The first condition in this definition of a coset cycle says that consecutive cosets 
are linked in the sense that they overlap in the named representatives; 
the second condition says that there is no direct
shortcut inside any one of these cosets,
from its immediate predecessor to its immediate successor.

\bD
\label{Nacycdef}
For $N \geq 2$, a Cayley groupoid $\G/\I$ is called \emph{$N$-acyclic} if it does
not have any coset cycles of lengths $2 \leq n \leq N$. 
\eD

The weakest meaningful coset acyclicity condition, $2$-acyclicity, 
deserves special attention. If we think
of $\alpha$-cosets as $\alpha$-connected components, then the
following observation suggests to picture 
$2$-acyclicity as a notion of simple connectivity:
the intersection of two connected pieces is itself
connected.

\bO
\label{twoacycobs}
For any groupoid $\G/\I$: $\G$ is $2$-acyclic if, and only if, 
for all $\alpha_1,\alpha_2 \subset E$ such that 
$\alpha_i = \alpha_i^{-1}$,
\[
\G[\alpha_0] \cap \G[\alpha_1]
= \G[\alpha_0 \cap \alpha_1].
\] 
\eO

\prf
Generally, $\G[\alpha_0 \cap \alpha_1] \subset \G[\alpha_0] \cap
\G[\alpha_1]$.

A coset cycle of length $2$ consists of two pointed cosets, 
$g_0 \G[\alpha_0]$ at $g_0$ and $g_1 \G[\alpha_1]$ at $g_1$, such that 
$g_0, g_1 \in g_0 \G[\alpha_0] \cap g_1 \G[\alpha_1]$ while  
$g_0 \G[\alpha_0 \cap \alpha_1]$ is disjoint from $g_1 \G[\alpha_0
\cap \alpha_1]$. Therefore 
\[
g_0^{-1} g_1 \in (\G[\alpha_0] \cap \G[\alpha_1])
\setminus \G[\alpha_0 \cap \alpha_1]
\]
shows that $\G[\alpha_0] \cap \G[\alpha_1]
\not= \G[\alpha_0 \cap \alpha_1]$.
Conversely, any $h \in (\G[\alpha_0] \cap \G[\alpha_1])
\setminus \G[\alpha_0 \cap \alpha_1]$ gives rise to a coset cycle 
consisting of $\G[\alpha_0]$ at $1_s$ for $s = \iota_1(h)$,  
and $h \G[\alpha_1]$ at $h$.
\eprf

\bR
\label{freegroupoidacycrem}
The free groupoid over $\I$ of 
Remark~\ref{freegroupoidrem}, which is infinite in non-trivial
cases, is acyclic in the sense that it has no coset cycles. 
\eR

The following was claimed in~\cite{lics2013},
in a flawed construction that attempted to generalise rather directly
the corresponding, much simpler result for Cayley 
groups from~\cite{OttoJACM}. 
A corrected treatment can be found in~\cite{Cayleynew20}.

\bT
\label{Nacycgroipoidthm}
For every finite incidence pattern $\I$ and $N \geq 2$, there are
finite groupoids $\G/\I$ that are $N$-acyclic. Moreover, such $\G/\I$
can be chosen to be simple and compatible with any given 
amalgamation pattern $\H/\I$, as well as compatible with any symmetries of
the underlying amalgamation pattern $\H/\I$ in the sense that 
every symmetry of $\H/\I$ extends to a symmetry of $\G/\I$.
\eT

Compare Definition~\ref{compatdef} in connection with the compatibility assertion,
and Definitions~\ref{amalgautdef}
and~\ref{innersymdef} in connection with
symmetries of $\G/\I$ and $\H/\I$, as well as Example~\ref{compatex} 
for an initial groupoid $\G_0/\I$ that is symmetric over and
compatible with $\H/\I$, and which the desired $\G/\I$ then is to extend. 

We isolate the specific symmetry properties of the theorem 
in the following notion of being \emph{fully symmetric}. 
The above then says that finite $N$-acyclic
groupoids can be obtained that are fully symmetric over the given data. 

Recall from Definition~\ref{innersymdef} that  
a groupoid $\G/\I$ is \emph{fully symmetric} over $\I$ if every
symmetry of the incidence pattern $\I$ extends to a symmetry of the 
groupoid $\G$.

\bD
\label{fullysymmetricGdef} 
A groupoid $\G/\I$ 
is called \emph{fully symmetric over $\H/\I$} if every
symmetry of $\I$ that is induced by a symmetry of the amalgamation pattern $\H/\I$ 
extends to a symmetry of the groupoid $\G$.
\eD

Note that the definition only refers to symmetries of the groupoid
$\G/\I$;
see Remark~\ref{transsymmrem}
regarding the rich additional structure of $\I$-rigid automorphisms of
the Cayley graph of $\G/\I$ over and above these.

\subsection{Measures of global consistency in direct products}

Direct products of amalgamation patterns $\H/\I$ with suitable groupoids 
will allow us to obtain strongly coherent and simple coverings
$\hat{\H}/\hat{\I}$, and, through these by 
Lemma~\ref{coverreallem}, realisations. 
The ensuing realisations themselves thus
materialise as reduced products, 
viz.\ quotients of the coverings $\hat{\H}/\hat{\I}$ obtained
as direct products w.r.t.\ the congruence
relation $\approx$ induced by the identifications according to the 
$\rho_{\hat{e}}$ of $\hat{\H}$. 

The qualitative difference between coherence and strong coherence may
be interpreted as a fundamental difference between element-wise versus 
tuple-wise consistency.%
\footnote{There seems to be a vague analogy 
between these discrete combinatorial notions and classical geometric 
notions of path independence and contractibility in homotopy.}
At the technical level, this difference manifests itself as
a difference between graph-like and hypergraph-like phenomena. 
While graph-like structures admit (even unbranched) finite $N$-acyclic 
coverings through direct products with Cayley graphs of
groups~\cite{OttoAPAL04,DO}, finite $N$-acyclic hypergraph coverings
(in general necessarily branched coverings) naturally arise as
reduced products with Cayely graphs of groupoids. 

\medskip
We show that direct products of amalgamation patterns with 
Cayley graphs of compatible finite 
groupoids $\G/\I$ produce finite coverings that are coherent, while
further direct products with Cayley graphs of compatible and $2$-acyclic groupoids
guarantee simplicity and strong coherence. 
The coverings involved really arise at the level of the underlying incidence
pattern $\I$, where they can be pictured as \emph{bisimilar coverings} in the
usual graph setting (where an edge- and vertex-coloured graph encodes
a transition system or Kripke structure). 

In this sense also a groupoid $\G/\I$ induces
a covering $\hat{\I}$ of $\I$ by the Cayley graph of $\G$, as follows. 
In a groupoid $\G/\I$, the projection map
$\pi \colon g \longmapsto \iota_2(g)$ 
induces an unbranched bisimilar covering of the multi-graph structure 
$\I$ by the multi-graph structure of the Cayley graph of $\G$. By this
we mean a homomorphism, mapping edges $e[g]
\in R_e$
at $g \in G$ with $\iota_2(g) = \iota_1(e)$ to the edge $e$ in $\I$, which also 
satisfies this unique lifting property: 
at every $g \in G$ with $\iota_2(g) = \iota_1(e)$,  
$e$ uniquely lifts to $e[g]$. Putting 
\[
  \hat{E} := \{ e[g]
  \colon e \in E, g \in G, \iota_2(g) = \iota_1(e)\},
\]
we obtain an incidence pattern $\hat{\I}$, which is based on $G$ as
its sort of sites and $\hat{E}$ as its sort of links:
\[
\hat{\I} = \hat{\I}(\G) = \bigl( G, \hat{E}, \iota, \cdot^{-1} \bigr). 
\]

This covering $\pi \colon \hat{\I} \rightarrow \I$ is a
rendering of the Cayley graph of the groupoid $\G/\I$ ---
refined to become an incidence pattern with individually 
labelled edges in $\hat{E} = \dot{\bigcup}_{e \in E} R_e$.
With the following definition we lift this natural covering
relationship between $\hat{\I}(\G)$ and $\I$ to a covering 
$\pi \colon \hat{\H} \rightarrow \H$ between 
amalgamation patterns $\hat{\H} = \H \otimes \G$ over $\hat{\I}$
and $\H$ over $\I$.

\bD
\label{prodgroupdef}
The \emph{direct product} of an amalgamation pattern 
\[
\H/\I =
(\I,\str{H}, \delta \colon H \rightarrow S, P, \eta \colon P \rightarrow E)
\] 
over the incidence pattern $\I = (S,E,\iota,\cdot^{-1})$ with a groupoid 
$\G/\I$ is the amalgamation pattern 
\[
\H \otimes \G := \bigl( \hat{\I}, 
\str{H} \otimes \G , \hat{\delta} \colon H \otimes G  \rightarrow G, P \otimes
G, \hat{\eta} \colon P \otimes G \rightarrow \hat{E} \bigr)
\]
over the incidence pattern 
$\hat{\I} = (G, \hat{E}, \iota, \cdot^{-1} )$ induced by $\G/\I$. 
The universe $H \otimes G$ of $\str{H} \otimes \G$ is the set 
\[
\barr{r@{\;=\;}l}
H \otimes G &
\bigl\{ (a,g) \in H \times G \colon \delta(a) = \iota_2(g)
\bigr\}
\\
\hnt
& \bigl\{ (a,g) \in H \times G \colon 
a \in A_s, g \in G[\ast,s] \mbox{ for some } s \in S
\bigr\},
\earr
\]
partitioned by $\hat{\delta}$ which is the projection to the second 
component. $\str{H} \otimes \G$ as a $\sigma$-structure is the disjoint
union of the $\str{A}_g := \str{A}_s \times \{ g \}$ for $g \in G[\ast,s], s \in S$;
$P \otimes G$ correspondingly consists of all pairs 
$(p,g) \in P \times G$ where $\iota_1(\eta(p)) = \iota_2(g)$ regarded 
as a partial isomorphism between  $\str{A}_g$ and $\str{A}_{ge^\G}$ in 
the natural manner. 
\eD

Compare Definitions~\ref{amalgcovdef} and~\ref{fullysymmetricGdef}  
for the following.

\bL
\label{coverinprodsymmlem}
Let $\H/\I$ be an amalgamation pattern, $\H \otimes \G$ its direct
product with a groupoid $\G/\I$, and consider $\H \otimes \G$ as an
amalgamation pattern over the incidence pattern $\hat{\I}$ induced by
the Cayley graph of $\G$. In this situation
the natural projection from the direct product $\H \otimes\G$ to $\H$,
$\pi_1\colon (a,g) \mapsto a$, is a covering
of the amalgamation pattern $\H/\I$ by $\H \otimes \G/\hat{\I}$.
If $\G/\I$ is fully symmetric over $\H/\I$, then this covering is fully
symmetric over $\H/\I$. 
\eL

\prf
That $\pi_1$ satisfies the conditions for a covering
is straightforward. For the lifting property in particular, any
$\rho_e$ at the level of $\H$ lifts to $\rho_e \times g \colon 
(a,g) \mapsto (\rho_e(a),g e^\G)$ at every
$g \in \G[\ast,\iota_2(g)]$ for $\iota_2(g) = \iota_1(e)$; 
this is the partial isomorphism between $\str{A}_g$ and $\str{A}_{ge^\G}$ 
associated in $\H \otimes \G$ with the link
$e[g]
\in \hat{E}$.
We note that this further implies that 
the following diagram, in which the vertical isomorphisms 
are induced by $\pi_1$, commutes for any pair of walks $\hat{w}$
in $\hat{\I}$ and  $w$ in $\I$ that are in such a 
projection/lifting relationship, and where therefore  
$g' = g w^\G$:
\[
\xymatrix{
*++{\str{A}_g}
\ar[r]^{\mbox{$\rho_{\hat{w}}$}} 
\ar@{-}[d]|{\rotatebox{270}{$\,\simeq\,$}}
&
*++{\str{A}_{g'}}
\ar@{-}[d]|{\rotatebox{270}{$\,\simeq\,$}}
\\
*+{\str{A}_{\iota_2(g)}}
\ar[r]^{\mbox{$\rho_{w}$}} 
&
*+{\str{A}_{\iota_2(g')}}
}
\]

For the symmetry claim we first observe that any pair of symmetries
of the Cayley graph of $\G/\I$ and of $\H/\I$ that 
agree w.r.t.\ their action on $\I$ extends to a symmetry of $\H \otimes
\G$ and indeed of the covering relationship between $\H \otimes \G$
and $\H$. The salient point here is that $\H \otimes \G$ is defined
in terms of the Cayley graph of $\G/\I$ (as an amalgamation
pattern) rather than in terms of $\G$ as a groupoid. 
It is therefore clear that, if  $\G/\I$ is fully symmetric over $\H/\I$, 
every symmetry of $\H/\I$ extends to a symmetry of the covering. 
For the claim of full symmetry of $\H\otimes \G$ over $\H$, in this situation,
it remains to argue that the group of $\H$-rigid symmetries acts 
transitively on each of the sets  
$\pi^{-1}(\{s\}) = G[\ast,s] \subset G$ (regarded as sets of sites of
$\hat{\I}$), which is immediate from Remark~\ref{transsymmrem}: the 
$\I$-rigid automorphisms of the Cayley graph of $\G$ act transitively 
on the sets $G[\ast,s]$.
\eprf

\bL
\label{coherprodlem}
For the covering of $\H/\I$ by $\H\otimes\G/\hat{\I}$: 
\bre
\item
if $\G/\I$ is compatible with $\H/\I$, then 
$\H\otimes\G/\hat{\I}$ is coherent;
\item
if $\H/\I$ is coherent and if
$\G/\I$ is simple and $2$-acyclic, then 
$\H\otimes\G/\hat{\I}$ is simple and strongly coherent.
\ere
\eL

\prf
In~(i), compatibility of $\G/\I$ with $\H/\I$ means that, for every 
$w \in \I^\ast$ with $w^\G = 1_s$ in $\G$, we must have
$\rho_w \subset \mbox{id}_{A_s} = \mathrm{id}_s$ in $\H$, 
where $\rho_w$ is the composition of the $\rho_e$ 
along the walk from $s$ to $s$ described by $w$ in $\I$. 
Coherence of $\H \otimes \G$ on the other hand requires that for every
walk $\hat{w}$ in $\hat{\I}$ that loops back from some $g$ to the same
$g$, $\rho_{\hat{w}} \subset \mbox{id}_{A_g}$, where
$\rho_{\hat{w}}$ is the composition of the
$\rho_{\hat{e}}$ along the walk $\hat{w}$. But that $\hat{w}$ 
loops at $g$ in $\hat{\I}$, means that $w^\G = 1_s$ for 
the natural projection $w$ of $\hat{w}$ induced by $\pi \colon 
e[g]
\mapsto e$ and for $s = \iota_1(w)$. 
This projection also associates the action of
$\hat{\I}\nt^\ast$ on $\hat{\H} = \H \otimes \G/\hat{\I}$ with 
the action of $\I^\ast$ on
$\H = \H/\I$, as indicated in the commuting diagram above, so that 
\[
\rho_{\hat{w}}(a,g) = (\rho_w(a), g w^\G).
\]
So, if $w^\G = 1_s$ and therefore $\rho_w \subset \mbox{id}_{A_s}$
by compatibility of $\G$, it follows that  $\rho_{\hat{w}} \subset \mbox{id}_{A_g}$ as
required for coherence of $\hat{\H}$. 

For~(ii), consider now a direct product $\hat{\H} = \H \otimes \G$ where
$\G/\I$ is both compatible with $\H/\I$ and $2$-acyclic and
assume that $\H/\I$ is coherent: compatibility of $\G/\I$ with
$\H/\I$ is implied by coherence of $\H/\I$, see
Lemma~\ref{cohercompatlem}.  Coherence of $\H/\I$ means that 
the equivalence relation $\approx$, which is induced by identification 
of linked  elements $a \approx \rho_e(a)$ (for 
$a, e$ with $a \in \mathrm{dom}(\rho_e)$) never identifies any two elements 
in the same $\str{A}_s$ (cf.\ Lemma~\ref{cohapproxlem}). In other
words, $\approx$-equivalence classes $[a]$ of elements 
$a \in H$ intersect any one partition set $A_s \subset H$ 
in at most one element. So there is, for every
$a \in H$, a partial function $f_a \colon S
\rightarrow H$ such that $[a] \cap A_s = \{ f_a(s) \}$ 
if $f_a(s)$ is defined, and empty otherwise.%
\footnote{Think of $f_a$ as describing a partial section 
that traces the element $a$ through the desired realisation.}  
In particular, every $a \in H$ defines a subset $\alpha_a \subset E$ of
generators according to 
\[
\alpha_a := \bigl\{
e \in E \colon f_a(s) \in \mathrm{dom}(\rho_e)  \mbox{ for } s = \iota_1(e)
\bigr\},
\]
which can be seen as the set of links that carry $a$. The set 
$\alpha_a$ is closed under the converse operation of $\I$ (i.e.\ closed
under inverses as a generator set in $\G$): 
$f_a(s) \in \mathrm{dom}(\rho_e)$
(and $f_a(s)$ defined) for $s = \iota_1(e)$ if, and only if,
$f_a(s') \in \mathrm{dom}(\rho_{e^{-1}})$
for $s' = \iota_1(e^{-1}) = \iota_2(s)$  
(and $f_a(s')$ defined). 

Towards strong coherence of $\hat{\H}/\hat{\I}$,
 consider walks $\hat{w}_i$ in $\hat{\I}$ that link the same sites 
$g$ and $g' = gh = g w_i^\G$, where the $w_i$ are the projections of
$\hat{w}_i$ to $\I$, for $i=1,2$. The partial
isomorphisms $\rho_{\hat{w}_i}$ between $\str{A}_g$ and $\str{A}_{g'}$ 
are lifts of their projections $\rho_{w_i}$ between $\str{A}_s$ and
$\str{A}_{s'}$, for $s = \iota_2(g) = \iota_1(w_i)$ and $s' = \iota_2(g')= \iota_2(w_i)$.
Let $d_i =\mathrm{dom}(\rho_{w_i}) \subset A_s$, and put
\[
\textstyle
\alpha_i := \bigcap_{a \in d_i} \alpha_a,
\]
which is the set of those links of $\H$ that preserve the $\approx$-equivalent
copies of all elements of $d_i$. By $2$-acyclicity of $\G$, the
pointed cosets $(g, g\G[\alpha_1])$ and $(g',g'\G[\alpha_2])$ cannot form a
coset cycle. Since $\rho_{w_1}$ maps $d_1$ onto its image at $g'$,  
$w_1$ is composed of generators in $\alpha_1$, whence $g' \in
g\G[\alpha_1]$; similarly, $g \in g'\G[\alpha_2]$ since $\rho_{w_2}$ maps
$d_2$ onto its image at $g'$ (or its inverse maps this image onto $d_2$).
So the second condition of Definition~\ref{cosetcycdef} must be
violated, which for this potential $2$-cycle means that $g,g'$ are in
the same coset w.r.t.\ $\alpha := \alpha_1 \cap \alpha_2$, i.e.\ that 
\[
w_1^\G = w_2^\G = w^\G \mbox{ for some } w \in \G[\alpha]. 
\]

It follows that $d_1 \cup d_2 \subset \mathrm{dom}(\rho_w)$. For
the unique lift of this walk $w$ in $\I$ to a walk
$\hat{w}$ from $g$ to $g'$ therefore 
\[
\rho_{\hat{w}} \supset \rho_{\hat{w}_i}
\]
provides a common extension. This shows strong coherence of
$\hat{\H}/\hat{\I}$. 

Moreover, simplicity and 
$2$-acyclicity of $\G$ and strong coherence of $\hat{\H}/\hat{\I}$
together imply simplicity of $\hat{\H}/\hat{\I}$, as shown by a simple
indirect argument. (We write $e^\G = e$ and $e[g] = (g,ge)$
in this setting of a simple $\G$.)
By strong coherence there is a maximal 
element among all the $\rho_{\hat{w}}$ that link $\str{A}_g$ and
$\str{A}_{ge}$ (note that $ge \not= g$ since $\G$ is simple, which
implies that $\hat{\I}$ is
simple in that $\hat{E}$ has no loops or multiple edges). 
If this $\rho_{\hat{w}}$ were not the same as 
$\rho_{(g,ge)}$ as required for simplicity, then the
domain of $\rho_{\hat{w}}$ would be strictly larger than 
$\mathrm{dom}(\rho_{(g,ge)})$. In projection to $\H$, 
$ \mathrm{dom}(\rho_e) \strictsubset  \mathrm{dom}(\rho_w)$ implies 
that $w$ cannot use the generator $e$, i.e.\ $w \in
\G[\alpha]$ for some $\alpha \subset E \setminus \{ e, e^{-1} \}$. 
But then $(g,g\G[\alpha])$ and $(ge, ge \G[e,e^{-1}])$ would form a coset
cycle of length~$2$ in $\G$.
\eprf

Theorem~\ref{Nacycgroipoidthm} together with Lemma~\ref{coherprodlem}
underpins the following, as the crucial step towards
our main result, Theorem~\ref{realisethm}, in the following section. 
For symmetries of coverings compare Definition~\ref{covsymmdef}. 

\bC
\label{doublecovercor}
Any finite amalgamation pattern $\H/\I$ possesses a covering by 
a finite, simple and strongly coherent amalgamation pattern
$\tilde{\H}/\tilde{\I}$. Moreover, the covering 
$\pi \colon \tilde{\H}/\tilde{\I} \rightarrow \H/\I$
can be chosen to be fully symmetric over $\H/\I$. 
\eC

\prf
For $\tilde{\H}$ we take a two-fold direct product  
$(\H \otimes \G) \otimes \hat{\G}$ as obtained in
Lemma~\ref{coherprodlem}, which is a covering of (a covering of)
$\H/\I$ by Lemma~\ref{coverinprodsymmlem}, and simple and strongly
coherent for suitable choices of $\G$ and $\hat{\G}$, by the previous 
lemma. We may choose $\G$ compatible with $\H/\I$, fully symmetric
over $\H/\I$, and $\hat{\G}$ simple and fully symmetric over 
$(\H \otimes \G)/\hat{\I}$, so that $\tilde{\H}/\tilde{\I}$ is  
fully symmetric over $\H/\I$: any symmetry of $\H/\I$ first lifts 
to $(\H \otimes \G)/\hat{\I}$ and then further to
$\tilde{\H}/\tilde{\I}$. The required $\H$-rigid symmetries that
relate sites of the covering over the same $\str{A}_s$ of $\H/\I$
are induced by lifts of $\I$-rigid symmetries of the Cayley graph of 
$\G$ to symmetries of 
$(\H \otimes \G) \otimes \hat{\G}$ (recall that the group of $\I$-rigid
symmetries of the Cayley graph of $\G$ is transitive on the sets $G[\ast,s]$
by Remark~\ref{transsymmrem}). Any such symmetry extends to an 
$\H$-rigid symmetry of $\H \otimes \G$: 
$\I$-rigid symmetries of the Cayley graph of $\G$ induce 
$\H$-rigid symmetries of $\H \otimes \G$ by
Lemma~\ref{coverinprodsymmlem}. 
Since $\hat{\G}$ is fully symmetric over $(\H\otimes \G)/\hat{\I}$, 
it extends this symmetry of $\H \otimes \G$ further to a symmetry 
of $(\H\otimes \G) \otimes \hat{\G}$, which remains $\H$-rigid, again
by Lemma~\ref{coverinprodsymmlem}.
\eprf

The following analysis shows that 
a covering as in Corollary~\ref{doublecovercor} 
can even be obtained as a direct product 
$\H\otimes \tilde{\G}$ for a suitable finite groupoid
  $\tilde{\G}/\I$ in one step. In other words, we want to replace the 
twofold direct product construction, which first
guarantees coherence, then strong coherence at the level of a
partially unfolded incidence pattern ($\hat{\I}$ as an unbranched
covering of $\I$), by a single direct product.

Let us consider an iterated direct product $\tilde{\H} := (\H \otimes \G)
\otimes \hat{\G}$
as from the above construction, i.e.\ based on the following:
\bre
\item[--]
$\H/\I$ and $\G/\I$ are over $\I = (S,E)$, $\G$ is simple;
\item[--]
$\G$ is compatible with $\H$ (so that $\H \otimes \G$ is coherent by Lemma~\ref{coherprodlem});
\item[--]
$(\H \otimes \G)/\hat{\I} = \hat{\H}/\hat{\I}$ and $\hat{\G}$ 
are over $\hat{\I} = (G,\hat{E})$, the incidence pattern induced by
the Cayley graph of $\G$, as an unbranched bisimilar covering $\pi \colon \hat{\I}
\rightarrow \I$, where $\pi \colon g \mapsto \iota_2(g)$;
\item[--]
$\hat{\G}/\hat{\I}$ is simple and $2$-acyclic 
(and by Lemma~\ref{cohercompatlem} compatible with
$\hat{\H}/\hat{\I}$, so that $(\H \otimes \G) \otimes \hat{\G}$ is
simple and strongly coherent by Lemma~\ref{coherprodlem});
\item[--]
$\hat{\G}/\hat{\I}$ is fully symmetric over $\hat{\I}$ 
(cf.\ Definition~\ref{innersymdef}): all symmetries of 
$\hat{\I}$ (i.e.\ of the  
Cayley graph of $\G$) extend to symmetries of the groupoid $\hat{\G}$;
\item[--]
$\tilde{\H}/\tilde{\I} = 
((\H \otimes \G) \otimes \hat{\G})/\tilde{\I}$ over 
$\tilde{\I}$, the incidence pattern induced by the Cayley graph of $\hat{\G}$. 
\ere

\bP
\label{twoprodinoneprop}
In this situation, the Cayley graph of $\hat{\G}$ also carries the structure
of an amalgamation pattern over the coarser incidence pattern
$\I$. This amalgamation pattern is complete and induces 
a groupoid $\tilde{\G}/\I$ (by groupoidal action on 
$\hat{\G}$ over $\I$, in the sense of Example~\ref{completemalgex})
such that 
\bae
\item
$\tilde{\G}/\I$ is simple and compatible with $\H/\I$; 
\item
$\tilde{\G}/\I$ is fully symmetric over $\H/\I$ if 
$\G$ is fully symmetric over $\H/\I$;
\item
$\tilde{\G}/\I$ is $N$-acyclic if $\hat{\G}/\hat{\I}$ is;
\item
the direct product $\H \otimes \tilde{\G}$
of $\H/\I$ with this $\tilde{\G}/\I$ is 
a covering of $\H/\I$ by a simple and strongly coherent 
amalgamation pattern;
\item
the covering by $\H \otimes \tilde{\G}$ is fully symmetric
over $\H/\I$ if $\G$ is fully symmetric over $\H/\I$. 
\eae
\eP

In the r\^ole of a simple and strongly coherent 
covering of $\H/\I$, the direct product $\H \otimes \tilde{\G}$
of $\H/\I$ with this $\tilde{\G}/\I$ can thus 
replace the nested direct product $(\H \otimes \G) \otimes \hat{\G}$ 
considered above. By combination with Lemmas~\ref{coverreallem}, 
a (fully symmetric) realisation of 
$\H$ can therefore also be obtained as a reduced product of $\H$ 
with a suitable finite groupoid in a single step (cf.\ Theorem~\ref{realisethm} below).

\bC
\label{singlecovercor}
For any $N \geq 2$, 
any finite amalgamation pattern $\H/\I$ possesses a covering by 
a finite, simple and strongly coherent amalgamation pattern
$\tilde{\H}/\tilde{\I}$ that is fully symmetric 
over $\H/\I$, and which is obtained as a direct product of 
$\H/\I$ with a suitable $N$-acyclic finite groupoid $\tilde{\G}/\I$.
\eC

\prf[Proof of Proposition~\ref{twoprodinoneprop}]
We start with the construction of the groupoid
$\tilde{\G}/\I$, which is obtained by re-interpreting the
Cayley graph of $\hat{\G}$ as an amalgamation pattern over the coarser
incidence pattern $\I$ rather than the built-in $\hat{\I}$. 
Recall that $\hat{\I}$ is based on the Cayley graph of $\G$, and
related to $\I$ by an unbranched covering 
$\pi \colon \hat{\I} \rightarrow \I$ that maps $g$ to $\iota_2(g)$.  

As $\G$ is simple we may identify $e^\G$ with $e$ and think of $E$ as a
subset $E \subset G$; we also write $e[g] = (g,ge)$ 
without ambiguity in this case. To cast the Cayley graph of $\hat{\G}$ as an
induced amalgamation pattern over $\I$, we merge pieces of the finer partition
according to $\hat{\I}$ to fit the coarser pattern $\I$. This
coarsening is induced by $\pi$, which maps $g \in G$ (as a sort label of $\hat{\I}$) to  
$\iota_2(g) \in S$, $e[g] = (g,ge)$ (as a link label of
$\hat{\I}$) to $e \in E$, and correspondingly associates 
$\hat{g} \in \hat{G}[\ast,g]$ with sort $s = \iota_2(g) \in S$. 
The resulting amalgamation pattern (presented both
in proper multi-sorted format and as a labelled tuple) is
\[
\hat{\G}/\I := \bigl( \I, \hat{G}, \delta \colon \hat{G} \rightarrow S, P, \eta
\colon P \rightarrow E \bigr)
=
\; \bigl( \hat{G}, (\hat{G}_s)_{s \in S}, (\hat{\rho}_e)_{e \in E} \bigr)
\]
where $\delta$ partitions $\hat{G}$ into the sets 
\[
\hat{G}_s  = \bigcup \bigl\{ \hat{G}[\ast,g]
\colon \iota_2(g) = s \bigr\}
\]
and each $\hat{\rho}_e \in P$ 
is the disjoint union of the local bijections induced by right
multiplications with generators $\hat{e} = e[g] = (g,ge) \in \hat{E}$, 
for which $\pi(\hat{e}) = e$. Since $\pi \colon \hat{\I} \rightarrow
\I$ is an unbranched covering, there always is a unique such $\hat{e}$ 
applicable at every $\hat{g} \in \hat{G}_s$, dependent on just 
the partition set $\hat{G}_g = \hat{G}[\ast,g]$ to which $\hat{g}$
belongs in $\hat{\G}/\hat{\I}$. 
More specifically, the lift of $e \in E$ with $\iota(e)= (s,s')$ in
$\I$ to some $\hat{g} \in \hat{g} \in \hat{G}[\ast,g] \subset
\hat{G}_s$ with $\iota_2(g) = s$ is $e[g] = (g,ge)$. 
In terms of its operation by right multiplication on the Cayley graph
of $\hat{\G}$, 
\[
\hat{\rho}_e := \bigcup_{\pi(\hat{e}) = e} \rho_{\hat{e}}^{\hat{\G}}
= \bigcup_{g \in G[\ast,\iota_1(e)]} 
\bigl( \rho_{e[g]}^{\hat{\G}} \colon 
\hat{G}_{\iota_1(e)} \longrightarrow \hat{G}_{\iota_2(e)}\bigr),
\]
each $\hat{\rho}_e$ is a bijection between the partition sets
($\hat{\G}/\I$ is complete). Extending this analysis to walks 
$w = e_1\cdots e_n$ from $\iota_1(w)$ to $\iota_2(w)$
in $\I$, we write $w[g]$ for the unique lifting of the generator
sequence $e_1 \cdots e_n$ that starts with $e_1[g]$ (at any 
$\hat{g} \in \hat{G}[\ast,g]$ for $g \in G[\ast,\iota_1(e_1)]$) and,
writing $w_i := e_1 \cdots e_i$ for the prefixes, produces the walk  
\[
w[g] = e_1[g] \cdot e_2[gw_1^\G] \cdots e_n[gw_{n-1}^\G] \in \hat{\I}^\ast.
\]

We obtain $\hat{\rho}_w$ as the composition of the 
$\hat{\rho}_{e_i}$ along this walk, 
\[
\hat{\rho}_w= \bigcup_{g \in G[\ast,\iota_1(w)]} 
\bigl( \rho_{w[g]}^{\hat{\G}} \colon
\hat{G}_{\iota_1(w)} \longrightarrow \hat{G}_{\iota_2(w)} \bigr),
\]
which is again a bijection between the partition sets involved, and
maps every
$\hat{g} \in \hat{G}[\ast,g] \subset \hat{G}_s$ (for $g \in
G[\ast,s]$, $s = \iota_1(e_1)$) according
to
\[
\hat{\rho}_w \colon 
\hat{g} 
\longmapsto  \hat{g}\cdot e_1[g] \cdot e_2[gw_1^\G]
\cdots e_n[gw_{n-1}^\G].
\]

We let $\tilde{\G}/\I$ be the groupoid generated by this action of the
$\rho_e^{\hat{\G}}$ for  
$e \in E$, in the inverse semigroup $I(\hat{G})$. Note how this
groupoid structure is induced in terms of the derived groupoidal
action of $\I^\ast$ on the Cayley graph of $\hat{\G}$, which combines 
in parallel the effect of all the $\rho_{e[g]}$ for $g \in G[\ast,\iota_1(e)]$ 
according to the native groupoidal action 
of $\hat{\I}^\ast$ on $\hat{\G}$. 
So the elements of $\tilde{\G}$ are 
the partial bijections $\hat{\rho}_w$ of $\hat{G}$
induced by walks $w$ in $\I$.

For a walk $w$ from $s$ to $s$ in $\I$, $\hat{\rho}_w$ represents $1_s$ in the 
groupoid $\tilde{\G}$ induced by this action, if, and only if, 
$\hat{\rho}_w$ is the identity of $\hat{G}_s$, i.e.\ if, and only
if, $\rho^{\hat{G}}_{w[g]}= \mathrm{id}_{\hat{G}[\ast,g]}$ for all $g
\in G[\ast,s]$, if, and only if, right multiplication with   
$(w[g])^{\hat{\G}}$ is a local identity 
on $\hat{G}[\ast,g]$ for all $g
\in G[\ast,s]$. 
But any two
$g,g' \in G[\ast,s]$ are related by a
unique $\I$-rigid automorphism $\sigma$ of $\hat{\I}$ (i.e.\ the Cayley graph
of $\G$, not the groupoid) that maps $g$ to $g'$.
Due to the assumption of full symmetry of $\hat{\G}$ over $\hat{\I}$,
this automorphism extends to a symmetry $\hat{\sigma}$ of the groupoid 
$\hat{\G}$, and this extension automorphically relates the actions of 
$\rho_{w[g]}^{\hat{\G}}$ at $\hat{G}[\ast,g]$ and 
$\rho_{w[g']}^{\hat{\G}}$ at $\hat{G}[\ast,g']$,
as indicated in the following commuting diagram:
\[
\nt\hspace{-.5cm}
\xymatrix{
G[\ast,s] \ar[dd]|{\makebox(10,10){$\sc\sigma$}} & *++{g} \ar[dd]|{\makebox(10,10){$\sc\sigma$}}
& *++{e} \ar[dd]|{\makebox(10,10){$\sc\sigma$}} & *++{e[g]} \ar[dd]|{\makebox(10,10){$\sc \hat{\sigma}$}} 
& *++{w[g]} \ar[dd]|{\makebox(10,10){$\sc \hat{\sigma}$}} 
& *++{\hat{G}[\ast,g]} \ar[dd]|{\makebox(10,10){$\sc \hat{\sigma}$}}
\ar[rr]^{\hat{\rho}_w}_{\rho^{\hat{\G}}_{w[g]}}
&& 
*++{\hat{G}[\ast,g w^\G]}  \ar[dd]|{\makebox(10,10){$\sc \hat{\sigma}$}}
\\
\\
*++{G[\ast,s]}  & g' 
& *++{e} & *++{e[g']} 
& *++{w[g']} 
& *++{\hat{G}[\ast,g']} 
\ar[rr]^{\hat{\rho}_w}_{\rho^{\hat{\G}}_{w[g']}}
&& 
\hat{G}[\ast,g' w^\G] 
}
\]

So the following are equivalent:
\bre
\item the walk 
$w = e_1\cdots e_n$  in $\I$ from $s$ to $s$ 
generates the unit $1_s$ in the groupoid $\tilde{\G}/\I$:
$w^{\tilde{\G}} = 1_s$.
\item
the lift $w[g]$ of this walk to a walk at any $\hat{g} \in
\hat{G}[\ast,g] \subset \hat{G}_s$
in the Cayley graph of $\hat{\G}$ loops back to $\hat{g}$, i.e.\ 
$(w[g])^{\hat{\G}}$ acts as the identity of $\hat{G}[\ast,g]$, for all $g \in G[\ast,s]$.
\item
the lift $w[g]$ of this walk to a walk at $\hat{g}$
loops back to $\hat{g}$ for some $\hat{g}  \in
\hat{G}[\ast,g] \subset \hat{G}_s$ and some $g \in G[\ast,s]$.
\ere

For claim~(a) of the proposition it is immediate that $\tilde{\G}/\I$ is simple, by simplicity of $\G$;
compatibility with $\H/\I$ follows from compatibility of $\G$ with $\H/\I$:
any $\hat{\rho}_w$ that represents a unit $1_s$ in $\tilde{\G}$
must, by construction of $\tilde{\G}$, act as the identity 
of $\hat{G}_s$ in the amalgamation pattern
$\hat{\G}/\I$, so that in particular $w^\G = 1_s$ in $\G$. 

For~(b), full symmetry of $\tilde{\G}/\I$ over $\H/\I$ is immediate if 
$\G$ itself is fully symmetric over $\H/\I$: any automorphism of
$\H/\I$ lifts to a groupoid automorphism of $\G$, which in particular
induces an automorphism of its Cayley graph, and hence of the
incidence pattern $\hat{\I}$, which in turn lifts to an automorphism of
$\hat{\G}$ and thus of the groupoid $\tilde{\G}/\I$ (whose definition in
terms of $\hat{\G}$ refers just to the Cayley graph of $\hat{\G}$ and 
partitions relative to $\I$). 

For~(c), we argue that any coset cycle in $\tilde{\G}/\I$ induces a
coset cycle of the same length in $\hat{\G}/\hat{\I}$. The relevant
cosets of $\tilde{\G}$ are generated by sub-groupoids of the form 
$\tilde{\G}[\alpha]$ for generator sets $\alpha = \alpha^{-1} \subset E$;
with this sub-groupoid of $\tilde{\G}$ we can associate the sub-groupoid
$\hat{\G}[\hat{\alpha}]$ 
of $\hat{\G}$, generated by
$\hat{\alpha} := \pi^{-1}(\alpha) = \{ e[g] = (g,ge) \in \hat{E} \colon e \in
\alpha, g \in G[\ast,\iota_1(e)] \}$. It remains to check that a
cyclic tuple of pointed cosets  
$\bigl( \tilde{g}_i,  \tilde{g}_i \tilde{\G}[\alpha_i]
\bigr)_{i \in \Z_n}$ that forms a coset cycle in $\tilde{\G}$,
for suitable choices of the $\hat{g}_i$ in $\hat{\G}$
translates into a cyclic tuple 
$\bigl( \hat{g}_i, \hat{g}_i \hat{\G}[\hat{\alpha}_i]
\bigr)_{i \in \Z_n}$, 
that forms a coset cycle in $\hat{\G}/\hat{\I}$ (cf.\ conditions~(i)
and~(ii) from Definition~\ref{cosetcycdef}):

(i)~From $\tilde{g}_{i+1} \in g_i \tilde{\G}[\alpha_i]$ we find $w_i \in
\alpha_i^\ast$ such that $\tilde{g}_{i+1} = \tilde{g}_i
w_i^{\tilde{\G}}$ and put $\hat{g}_{i+1} := \hat{g}_i \cdot
(w_i[g_i])^{\hat{\G}}$. Here $w_i[g_i]$ is the unique lift of the 
walk $w_i$ in $\I$ to a walk in $\hat{\I}$ at $g_i$ if $\hat{g}_i \in
\hat{G}[\ast,g_i]$, and consists of links in $\hat{\alpha}_i$.
The choice of a starting element, say $\hat{g}_0$, is arbitrary, but 
it is important that tracing the sequence of the $w_i[g_i]$ in
$\hat{\G}$ takes us back to that starting point --- which it does because the 
product/concatenation of the $w_i$ represents a unit in $\tilde{\G}$,
which implies the same for their lifts in $\hat{\G}$.

(ii)~From $\tilde{g}_i \tilde{\G}[\alpha_i\cap \alpha_{i-1}] 
\cap \tilde{g}_{i+1} \tilde{G}[\alpha_i\cap
\alpha_{i+1}] = \emptyset$, we directly infer the corresponding
condition  for the translation to $\hat{\G}$ as follows. If $\hat{g}$ were in the
corresponding intersection in $\hat{\G}$, it could be represented as
$\hat{g} = \hat{g}_i \hat{u}^{\tilde{\G}} = \hat{g}_{i+1} \hat{w}^{\tilde{\G}}$ for
walks $\hat{u} \in (\hat{\alpha}_i\cap \hat{\alpha}_{i-1})^\ast$ and 
$\hat{w} \in (\hat{\alpha}_i\cap \hat{\alpha}_{i+1})^\ast$ in
$\hat{\I}$. These walks project to walks $u = \pi(\hat{u})$ and $w =
\pi(\hat{u})$ in $\I$, which
would show that also $\tilde{g}_i \tilde{\G}[\alpha_i\cap \alpha_{i-1}] 
\cap \tilde{g}_{i+1} \tilde{G}[\alpha_i\cap
\alpha_{i+1}] \not= \emptyset$.

For~(d), the direct product $\tilde{\H}' := \H\otimes \tilde{\G}$ is a
covering of $\H/\I$ by definition;
it remains to show that $\tilde{\H}'/\tilde{\I}'$ is simple and strongly
coherent as an amalgamation pattern over the incidence pattern
$\tilde{\I}'$ induced by the Cayley graph of $\tilde{\G}$.
Recall from Definition~\ref{coherencedef} that simplicity requires  
the partial isomorphisms $\rho_{\tilde{e}}$ in $\tilde{\H}'$ 
to be maximal among all 
$\rho_w$ in $\tilde{\H}'$ along walks $\tilde{w}$ from $\iota_1(\tilde{e})$ to
$\iota_2(\tilde{e})$ in $\tilde{\I}'$. 
The links $\tilde{e}$ of $\tilde{\I}'$ correspond to right
multiplication in $\tilde{\G}$ by generators $e^{\tilde{\G}}$ at
$\tilde{g} \in \tilde{G}$, which in 
$\tilde{\H}' = \H \otimes \tilde{\G}$ act as 
\[
\rho_{\tilde{e}} \colon
(h,\tilde{g}) \longmapsto 
(\rho_e(h), \tilde{g} e^{\tilde{\G}})
\] 
on $h \in A_{s}$ for $s = \iota_1(e)$. Any walk $\tilde{w}$ from
$\tilde{g} = \iota_1(\tilde{e})$ to $\tilde{g}'=\iota_2(\tilde{e})$ in 
$\tilde{\I}'$ is such that $\tilde{w}\nt^{\tilde{\G}} = e^{\tilde{\G}}$ in $\tilde{\G}$. Since 
$\tilde{\G}$ is a groupoid over $\I$, this implies in particular that, 
for the projection $w$ of $\tilde{w}$ to $\I$, 
$w^\G = e$  in $\G$, whence $\rho_w(h) = \rho_e(h)$ (where defined) is 
guaranteed by compatibility of $\G$ with $\H$.

Strong coherence, by Definition~\ref{coherencedef}, requires that, 
for any two walks $\tilde{w}_1, \tilde{w}_2$ 
from $\tilde{g}$ to $\tilde{g}'$ in $\tilde{\I}'$, there is a walk 
$\tilde{w}$ such that $\rho_{\tilde{w}}$ in $\tilde{\H}'$ is a common extension of
the $\rho_{\tilde{w}_i}$ for $i=1,2$. The reason essentially is that 
$\tilde{w}_1^{\tilde{\G}} = \tilde{w}_2^{\tilde{\G}}$ implies that 
the actions of the $\tilde{w}_i$ on $((\H \otimes \G) \otimes
\hat{\G})$ agree, where strong coherence yields a common 
extension. More specifically, consider the effect of the 
actions of $\tilde{e}$ (links in $\tilde{\I}'$) along the walks
$\tilde{w}_i$ in $\tilde{\I}'$ in the groupoid component $\tilde{\G}$
of $\tilde{\H}' = \H \otimes \tilde{\G}$.
The exploration of $\tilde{\G}$ given above shows that $\tilde{w}_1^{\tilde{\G}} =
\tilde{w}_2^{\tilde{\G}}$ implies that the projections of 
the walks $\tilde{w}_i$ to walks $w_i$ in $\I$ induce the same action
on $\hat{\G}$. In terms of matching walks $w_i[g]$ in the Cayley
graph of $\hat{\G}$, which are the lifts of the $w_i$ at appropriate
elements $\hat{g} \in \hat{\G}[\ast,g]$, this
implies that $(w_1[g])^{\hat{\G}} = (w_2[g])^{\hat{\G}}$.
Strong coherence of $((\H \otimes \G) \otimes \hat{\G})$ therefore
implies the existence of a walk $w[g]$ in the Cayley graph of
$\hat{\G}$ whose action on $((\H \otimes \G) \otimes \hat{\G})$
extends those of the $w_i[g]$. Now $(w[g])^{\hat{\G}} =
(w_i[g])^{\hat{\G}}$ for all appropriate $g$ 
implies that also
$w^{\tilde{\G}} = \tilde{w}^{\tilde{\G}}$, since these groupoid
elements are defined in term of the action on $\hat{\G}$.
Since their action in the $\H$-component of 
$\tilde{\H}' = \H \otimes \tilde{\G}$ is fully determined by the
projections $w$ and $w_i$ to $\I$ in both $((\H \otimes \G) \otimes
\hat{\G})$ and $\H \otimes \tilde{\G}$, $\rho_{w_i} \subset \rho_w$ in
$\H$ implies the same also over $\H \otimes \tilde{\G}$.

Claim~(e) follows from~(b) together with Lemma~\ref{coverinprodsymmlem}.
\eprf

\section{Generic realisations in reduced products} 
\label{redprodrealsec}

Recall the equivalence relation
$\approx$ defined in connection with Lemma~\ref{cohapproxlem}:
it is the equivalence relation induced on the universe $H =
\dot{\bigcup}_{s \in S} A_s$ of an
amalgamation pattern $\H = (\str{H},
(\str{A}_s)_{s \in S}, (\rho_e)_{e \in E})$ over $\I = (S,E)$ 
by regarding as equivalent any elements $a \in A_s$ and $\rho_e(a) \in
A_{s'}$ for links $e \in E$ with $\iota(e) = (s,s')$.  Also recall the
notion of an atlas from Definition~\ref{atlasdef}, and of a direct
product $\H \otimes \G$ from Definition~\ref{prodgroupdef}.  
The $\approx$-equivalence class of $a \in H$ is denoted as $[a]$.

\bD
\label{redproddef}
Let $\G/\I$ be compatible with $\H/\I$. The \emph{reduced product} 
$\H\otimes \G/{\approx}$ 
of an amalgamation pattern $\H/\I$ 
with a groupoid $\G/\I$ is based on the quotient structure of the 
relational structure $\hat{\str{H}}$ of the 
direct product $\hat{\H} = \H\otimes \G$ w.r.t.\ the equivalence
relation $\approx$ whenever this quotient is well-defined.%
\footnote{I.e.\ whenever the direct product $\hat{\H}$ is globally
  consistent, which is in particular the case if the direct product 
is strongly coherent, by 
Lemma~\ref{cohapproxlem}.}
We endow this relational structure $\hat{\str{H}}/{\approx}$ with 
an atlas induced by the families of subsets 
\[
U_s := \bigl\{ u[g] \colon g \in G[\ast,s] \bigr\}
\mbox{ where }
u[g] = [A_{\iota_2(g)} \times \{g\}] =  
\bigl\{ [(a,g)] \colon a \in A_{\iota_2(g)} \bigr\}
\]
for $s \in S$, with the natural isomorphisms 
\[
\barr{@{}rcl@{}}
\pi_{u,s} \colon \hat{\str{H}}/{\approx} &\rightarrow & \str{A}_s
\\
{[(a,g)]} &\longmapsto & a
\earr
\]
for $u = u[g] \in U_s$, which is well-defined due to
coherence of the product. With this atlas on 
$U = \bigcup_{s \in S} U_s$ over the relational 
structure $\hat{\str{H}}/{\approx}$, we define the reduced product as
\[
\H\otimes \G/{\approx} :=
\bigl( \hat{\str{H}}/{\approx},U,(U_s)_{s \in S}, (\pi_{u,s})_{u\in U_s}\bigr).
\]
\eD

It is clear from the format of these reduced products in relation 
to $\H/\I$ that they are natural candidates for realisations of
$\H/\I$ according to Definition~\ref{realdef}. And indeed, whenever 
the direct product $\H\otimes \G$ is simple and strongly coherent,
Lemma~\ref{coverreallem} shows that $\H\otimes \G/{\approx}$ is a
realisation of $\H \otimes \G$ which induces a realisation of $\H/\I$ 
by Remark~\ref{coverealrem}. 
So Corollary~\ref{singlecovercor} guarantees the existence of finite realisations 
of $\H/\I$ based on reduced products of $\H/\I$ with suitable finite
groupoids $\G/\I$, as available by Theorem~\ref{Nacycgroipoidthm}.

\bT
\label{realisethm}
Every finite amalgamation pattern $\H/\I$ admits finite realisations
by natural reduced products 
$\H \otimes \G/{\approx}$ 
based on a direct product $\H\otimes \G$ of $\H$ with a suitable
finite groupoid $\G/\I$, such that $\H \otimes \G$ itself is a simple,
strongly coherent covering of $\H$. Such realisations can be chosen to 
be fully symmetric over the given amalgamation pattern $\H/\I$, and 
$\G/\I$ can be chosen to be $N$-acyclic for any given $N \geq 2$.
\eT

\bR
\label{caninfasredprodrealrem}
The canonical  infinite realisation of $\H/\I$ of
Remark~\ref{caninfrealrem} can similarly be cast as a reduced product of
$\H$ with the free groupoid over $\I$ of 
Remark~\ref{freegroupoidrem}.
\eR

\subsection{Degrees of acyclicity in realisations} 

As seen above, Corollary~\ref{singlecovercor} together with 
Lemma~\ref{coverreallem} offer a route to finite realisations
with additional acyclicity properties, obtained as 
reduced products with groupoids that are $N$-acyclic 
for any chosen $N \geq 2$. We recall that  $2$-acyclicity 
is essential towards strong coherence. We now want to
analyse the effect that higher levels of coset acyclicity 
may have on the resulting realisation.
To this end we first review an established
notion of hypergraph acyclicity (\cite{Berge}, also known as
$\alpha$-acyclicity in the literature, e.g.~\cite{BeeriFaginetal}). Graded,
i.e.\ quantitatively localised or size-bounded  
versions of this natural
notion can be applied to the hypergraph structure 
of the atlas of a finite realisation.
Realisations obtained as reduced products with suitable $N$-acyclic
groupoids can then be shown to have an atlas that is $N$-acyclic
in the sense of hypergraph acyclicity. This means that the family 
of distinguished substructures induced by the co-ordinate domains,
or the isomorphic copies of the template structures $\str{A}_s$ 
that make up the realisation, locally overlap in a tree-like
(viz.\ tree-decomposable) fashion. This will further imply 
graded local universality properties for the realisation.

\subsection{Hypergraph acyclicity}

\bD
\label{hypgraphdef}
A \emph{hypergraph} $(A,U)$ is a structure consisting of a
\emph{vertex set} $A$ and a \emph{set of hyperedges} $U \subset 
\mathcal{P}_{\mathrm{fin}} (A)$ such that $A = \bigcup U$.%
\footnote{Hyperedges are to
  be (non-empty) finite subsets of the set of vertices even when the
  set of vertices is infinite; forbidding vertices that are outside all hyperedges is
  essentially w.l.o.g.\ (noting that singleton hyperedges are
  allowed) but helps to avoid trivial special provisions
  in some constructions and arguments.}
\eD

With a hypergraph $(A,U)$, with set of hyperedges $U \subset
\mathcal{P}_{\mathrm{fin}} (A)$ over vertex set $A$, we associate its \emph{Gaifman
  graph}, which is a simple undirected graph over the same vertex
set. Its edge relation $R$ links two 
distinct vertices $a$ and $a'$ if they occur in  a common
hyperedge:
\[
R = \{ (a,a') \in A^2 \colon a \not= a', a,a' \in u \mbox{ for some } u
\in U \}.
\]

In other words, the Gaifman graph consists of the union of cliques 
over the $u \in U$.

Recall that a \emph{chord} in a cycle $(a_i)_{i \in \Z_n}$ is an edge  
linking two vertices that are not next neighbours in the cycle (which
can only occur in cycles of length greater than~$3$).

\bD
\label{hypacycdef}
A hypergraph $(A,U)$ is \emph{acyclic} if it is both
\emph{chordal} and \emph{conformal}, i.e.\ if
\bre
\item
its Gaifman graph has no chordless cycles of length greater than~$3$ (chordality);
\item
every clique in its Gaifman graph is fully contained in some hyperedge (conformality).
\ere  
\eD

Useful equivalent characterisations in the case
of finite hypergraphs involve the notion of
\emph{tree-de\-com\-posa\-bi\-li\-ty} and a notion of
\emph{retractability} (cf.~\cite{BeeriFaginetal}).

\bD
\label{treedecompdef}
A \emph{tree decomposition} $(T,E,\lambda \colon T \rightarrow U)$  
of a finite hypergraph $(A,U)$ consists of a labelling of the hyperedges 
$u \in U$ by the
vertices of a finite (graph-theoretic) tree $(T, E)$ via a surjective function 
$\lambda\colon T \rightarrow U$ such that, for every $a \in A$, the subset 
$\{ t \in T \colon a \in \lambda(t) \} \subset T$ is connected in
$(T,E)$. 

A finite hypergraph $(A,U)$ is \emph{retractable} (in the sense of Graham's
algorithm) if it can be reduced to the empty structure be recursive
application of the following retraction steps:
\bre
\item[--]
deletion from $A$ and every $u \in U$ of 
a single vertex $a \in A$ that is an element  of just one hyperedge; 
\item[--]
deletion from $U$ of a single hyperedge 
that is a subset of some other hyperedge.
\ere
\eD

It is obvious that a hypergraph is acyclic if, and only if, 
all its finite induced sub-hypergraphs are acyclic. For finite hypergraphs 
it is also well known, e.g.\ from~\cite{BeeriFaginetal}, that 
acyclicity coincides with the existence of a tree decomposition as well
as with retractability. Rather than outright acyclicity, we shall
focus on a notion of $N$-acyclicity as a suitable quantitative
approximation to acyclicity, so that 
acyclicity becomes the limit across all
levels of $N$-acyclicity.

\bD
\label{hypNacycdef}
For $N \geq 3$, a hypergraph $(A,U)$ is \emph{$N$-acyclic} if 
every induced sub-hyper\-graph with up to $N$ many vertices
is acyclic.
\eD

It is easy to check that $N$-acyclicity is equivalent to the
combination of 
\bre
\item
$N$-chordality:
its Gaifman graph has no chordless cycles of lengths greater than~$3$ 
and up to $N$;
\item
$N$-conformality:
every clique of size up to $N$ in its Gaifman graph is contained in some hyperedge.
\ere  

In the following we also refer to \emph{$N$-acyclic atlases}
or to \emph{$N$-acyclic realisations}, meaning that the hypergraph
associated with the co-ordinate domains of the atlas is $N$-acyclic.

\bO
\label{atlastreedecompobs}
If $(A,U)$ is the hypergraph of co-ordinate
domains of an atlas of a finite relational structure $\str{A}$ on $A$
according to Definition~\ref{atlasdef}, then
a tree decomposition $(T,E,\lambda \colon T \rightarrow U)$ 
of $(A,U)$ is such that every
tuple in a relation of $\str{A}$ must be fully contained in one of the
$\lambda(t)$ and such that all elements in the overlap
$\lambda(t_1) \cap \lambda(t_2)$ of two charts must be represented in
every $\lambda(t)$ along the shortest connecting path between $t_1$
and $t_2$ in the tree $T$. Model-theoretically, such a tree
decomposition provides a representation of the underlying structure as
a free amalgam of its distinguished substructures $\str{A}\restr u$
for $u \in U$. The free amalgam
can be put together in a stepwise fashion by following in reverse
order a retraction as in Definition~\ref{treedecompdef}.
\eO

\subsection{Acyclicity in reduced products}

Compare Definition~\ref{redproddef} for reduced products and recall
that compatibility of the groupoid $\G/\I$ with the amalgamation
pattern $\H/\I$ is equivalent with the coherence of  
the direct product $\H \otimes \G$. The strictly stronger assumption of 
strong coherence of the direct product $\H \otimes \G$ implies that 
the reduced product $\H\otimes \G/{\approx}$ is well-defined. 
The next lemma refers to the hypergraph induced by the
natural atlas of this reduced product $\H\otimes \G/{\approx}$, whose
hyperedges are the co-ordinate domains of its charts onto the
$\str{A}_s$ of $\H/\I$,
\[
u[g] := [A_s \times \{g\}] = \{ [(a,g)] \colon a \in A_s \},
\]
for all $g \in G[\ast,s]$, $s \in S$ (note that $g \in \G$ fully
determines $u[g]$ since $s = \iota_2(g)$). 
Note that an element $[(a,g)]$ of the reduced product, which is represented
by a pair $(a,g)$ of the direct product with $a \in A_{s}$ and $g \in
G[\ast,s]$ for some $s \in S$ (which thus is an element of $u[g]$), 
is also an element of $u[h]$ for some $h \in \G[\ast,s']$ if, and only if,
it is represented by some pair $(a',h) \approx (a,g)$ for some 
$a' \in A_{s'}$. This is the case if, and only if, there is a walk 
$w$ from $s$ to $s'$ in $\I$ such that 
$h = g w^\G$ in $\G$ and $\rho_w(a) = a'$ in $\H$. We observe
that any such $w$ must be comprised of generators $e \in E$ for which 
$\rho_{w'}(a) \in \mathrm{dom}(\rho_e)$ for some $w'$ from $s$ to
$\iota_1(e)$ for which $a \in \mathrm{dom}(\rho_{w'})$, so that 
$(a,g) \approx (\rho_{w'}(a), gw'\nt^\G)
\approx (\rho_{w'e}(a), g(w'e)\nt^\G)$. It will therefore be
convenient to define these generator sets as
\[
\alpha_{a,g} = \{ e \in E \colon \rho_w(a) 
\in \mathrm{dom}(\rho_e) \mbox{ for some walk $w$ from $s$ to 
$\iota_1(e)$ in $\I$ } \}.
\]

Note that $\alpha_{a,g}^{-1} = \alpha_{a,g}$ is closed under edge reversal (inversion of
generators) since $\rho_w(a) \in \mathrm{dom}(\rho_e)$ if, and only if,
$\rho_{we}(a) \in \mathrm{dom}(\rho_{e}^{-1})=\mathrm{dom}(\rho_{e^{-1}})$.
The connection between cosets w.r.t.\ these
generator sets and the hypergraph structure of the atlas of
$\H\otimes\G/{\approx}$ is that
$g\G[\alpha_{a,g}] \subset G$ is
the set of groupoid elements
at which $[(a,g)]$ is represented.

\bO
\label{gensetobs}
In a reduced product $\H\otimes\G/{\approx}$, 
for $s \in S$, $a \in A_s \subset H$, $g \in G[\ast,s]$ and
$\alpha_{a,g} \subset E$ as above:
\[
\{ h \in \G \colon [(a,g)] \in u[h] \} = g\G[\alpha_{a,g}].
\]
\eO

For $N$-acyclicity of 
$\G/\I$ see Definition~\ref{Nacycdef} in connection with the following.

\bL
\label{nacycreallem}
For any amalgamation pattern $\H/\I$ and 
finite groupoid $\G/\I$ such that the direct product $\H\otimes\G$ is
strongly coherent: if $\G/\I$ is $N$-acyclic for some $N \geq 3$, then the 
hypergraph induced by the
co-ordinate domains $u[g]$ of the natural atlas of the reduced product
$\H\otimes\G/{\approx}$ 
is $N$-acyclic.
\eL

\prf
Denote as $(A,U)$ the 
hypergraph formed by the universe $A$ of the relational structure 
$\str{A}$ of the reduced product $\H \otimes \G/{\approx}$
together with the set of co-ordinate domains of its atlas. 
Its vertex set $A$ consists of
the $\approx$-equivalence classes $[(a,g)]$ of elements $(a,g)$ of the direct
product 
and its  hyperedges are the sets 
\[
u[g] := [A_{\iota_2(g)} \times \{g\}]  = \{ [(a,g)] \colon a \in
A_{\iota_2(g)} \}.
\]
Assume that $\G/\I$ is $N$-acyclic. We show that $(A,U)$ is $N$-acyclic.

\medskip\noindent\emph{$N$-chordality.}
For $N$-chordality, we need to show that any chordless cycle
in the Gaifman graph of $(A,U)$ must have length greater than $N$.
It suffices to argue that any chordless cycle in $(A,U)$ induces a
coset cycle of the same length in $\G$. Let $(c_i)_{i \in \Z_n}$ be a
chordless cycle in the Gaifman graph of $(A,U)$: for each $i \in \Z_n$
there is some hyperedge $u_i \in U$ such that 
$c_{i-1}, c_i \in u_i$, while there is no such $u \in U$ 
for the pair $c_i,c_j$ if $j \not= i\pm 1$ in $\Z_n$. Fix a sequence 
of such $u_i = u[g_i]$ so that $c_i \in u_{i},u_{i+1}$.
Since $c_i \in u_i$ we may choose representatives such that 
$c_i =[(a_i,g_i)]$. We 
associate with each $[(a_i,g_i)]$ the generator set $\alpha_i =
\alpha_{a_i,g_i} \subset E$ discussed in Observation~\ref{gensetobs}:
\[
\alpha_i = \{ e \in E \colon \rho_w(a_i) 
\in \mathrm{dom}(\rho_e) \mbox{ for some walk $w$ from $s_i$ to 
$\iota_1(e)$ in $\I$ } \}.
\]

With $i \in \Z_n$ we associate the 
sub-groupoid $\G[\alpha_i]$ and the pointed coset $(g_i, g_i \G[\alpha_i])$,
and show that 
\[
\bigl(g_i, g_i \G[\alpha_i]\bigr)_{i \in \Z_i}
\]
forms a coset cycle in $\G/\I$ by checking the two defining conditions of 
Definition~\ref{cosetcycdef}:
\bre
\item
$g_{i+1} \in g_i \G[\alpha_i]$;
\item
$g_i \G[\alpha_i\cap \alpha_{i-1}] \cap g_{i+1} \G[\alpha_i\cap
\alpha_{i+1}] = \emptyset$.
\ere

For (i), observe that $c_{i} = [(a_{i},g_{i})] \in u_{i+1} = u[g_{i+1}]$ 
implies that there is some walk $w$ from $s_i$ to $s_{i+1}$ 
in $\I$ such that $g_{i+1} = g_i w^\G$ in $\G$ and such that 
$a_i \in A_{s_i}$ is mapped to some element of 
$A_{s_{i+1}}$ by $\rho_w$ in $\H/\I$ (cf.\ Observation~\ref{gensetobs}).
It follows that this path is made up of links $e \in \alpha_i$ so that
$w^\G \in \G[\alpha_i]$ whence $g_{i+1} \in g_i \G[\alpha_i]$.

The argument for (ii) relies on the fact that the cycle 
formed by the $c_i$ in $(A,U)$ is chordless: 
suppose that, to the contrary of~(ii), there were some groupoid element 
$g \in g_i \G[\alpha_i\cap \alpha_{i-1}] \cap g_{i+1} \G[\alpha_i\cap
\alpha_{i+1}]$. 
By Observation~\ref{gensetobs}, $c_j = [(a_j,g_j)] \in u[g]$ if 
$g \in g_j \G[\alpha_j]$. So obviously $c_i,c_{i+1} \in u[g]$. 
But also $c_{i-1} \in u_i=u[g_i]$ by our initial choice of the $u_i$,
whence, by Observation~\ref{gensetobs}, $g_i \in g_{i-1}
\G[\alpha_{i-1}]$. Since 
$g \in g_i \G[\alpha_{i-1}]$, this implies that $g \in g_{i-1}\G[\alpha_{i-1}]$
too, whence, again by  Observation~\ref{gensetobs}, $c_{i-1} \in u[g]$
as well. So $u[g]$ would be a chord.

\medskip\noindent\emph{$N$-conformality.}
Towards an indirect proof of $N$-conformality consider a 
minimal clique $C \subset A$ in the Gaifman graph of $(A,U)$ 
that is not contained in one of the 
hyperedges $u \in U$; minimality here means that every proper subset of $C$ is
contained in some hyperedge. It suffices to show that any such clique 
induces a coset cycle of length $|C|$ in $\G$. Enumerate $C$ as 
$C = \{ c_i \colon i \in \Z_n \}$ for $n = |C| \geq 3$.
For $i \in \Z_n$ let $u_i$ be 
a hyperedge $u_i = u[k_i]$ that contains $C \setminus \{ c_i \}$, so that
$c_i \in u_{j}$ if, and only if, $j \not= i$.
Also fix a choice of representatives 
for the vertices $c_i \in C$ in the direct product $\H \otimes \G$,
according to $c_i = [(a_i,g_{i})]$ with $a_i \in A_{s_{i}}$ and 
$g_{i} \in G[\ast,s_{i}]$, and let as above 
\[
\alpha_i = \{ e \in E \colon \rho_w(a_{i}) 
\in \mathrm{dom}(\rho_e) \mbox{ for some walk $w$ from $s_{i}$ to 
$\iota_1(e)$ in $\I$ } \}.
\]
Further let, for $i \in \Z_n$, $\beta_i := \bigcap_{j \not=i,i+1}
\alpha_i$ and 
$h_i := k_{i}^{-1} k_{i+1}$.
We claim that 
\[
\bigl( h_i, h_i \G[\beta_i] \bigr)_{i \in \Z_n}
\] 
forms a coset cycle in $\G$, and need to establish conditions~(i)
and~(ii) as above.

Towards~(i) we note that
$c_{j} \in u_{i} = u[k_{i}]$ for $i \not=j$,
by Observation~\ref{gensetobs}
implies that $k_i \in g_j\G[\alpha_j]$ for all $i \not= j$ so that in particular
$k_{i}^{-1} k_{i+1} \in \G[\alpha_j]$ for all $j \not= i,i+1$. 
Since $\G$ is $2$-acyclic, Observation~\ref{twoacycobs} implies that 
$k_{i}^{-1} k_{i+1} \in \bigcap_{j \not= i,i+1} \G[\alpha_j] = \G[\beta_i]$,
or that $k_{i+1} \in k_i \G[\beta_i]$ as desired.

For~(ii), assume to the contrary that 
there were some 
$k \in k_i \G[\beta_i \cap \beta_{i-1}] \cap k_{i+1}
\G[\beta_i \cap \beta_{i+1}]$. Note that 
$\beta_i \cap \beta_{i-1} = \bigcap_{j \not= i} \alpha_j$
and $\beta_i \cap \beta_{i+1} = \bigcap_{j \not= i+1} \alpha_j$.
We claim that $u[k]$ would contain all of $C$,
contradicting assumptions. 

Indeed, for $j \not= i$, $c_j = [(a_j,g_j)] \in u_i = u[k_i]$ so that,
by Observation~\ref{gensetobs}, 
$k_i \in g_j \G[\alpha_j]$. Since also $k \in k_i \G[\beta_i \cap
\beta_{i-1}] \subset \G[\alpha_j]$ for all $j \not= i$, 
$k \in g_j\G[\alpha_j]$ and $c_j \in u[k]$ for $j \not= i$ follows,
again by Observation~\ref{gensetobs}. 
 A strictly analogous argument, based on $k \in k_{i+1} \G[\beta_i \cap
\beta_{i+1}]$ and $c_j \in u_{i+1} = u[k_{i+1}]$ for $j \not= i+1$,
shows that $c_j \in  u[k]$ for all $j \not= i+1$.
Overall, therefore $c_j \in u[k]$ for all $i$, the desired contradiction.
\eprf

\subsection{Acyclicity in realisations}

Together with Theorem~\ref{realisethm} we therefore find that not just 
fully symmetric realisations but fully symmetric realisations of any
degree of hypergraph acyclicity can be obtained as reduced products
with suitable finite groupoids.  Section~\ref{Nunivsec} below will show why 
$N$-acyclic realisations are particularly desirable.

\bC
\label{nacycrealcorr}
For any 
$N \geq 2$, every finite amalgamation pattern $\H/\I$ admits finite realisations
by natural reduced products 
$\H \otimes \G/{\approx}$ 
based on a direct product $\H\otimes \G$ of $\H$ with a suitable
finite $N$-acyclic groupoid $\G/\I$, such that $\H \otimes \G$ itself is a simple,
strongly coherent covering of $\H$. Such realisations can be chosen to 
be fully symmetric over the given amalgamation pattern $\H/\I$ and 
$N$-acyclic (in the sense that the hypergraph of co-ordinate domains,
or of the isomorphic copies of the $\str{A}_s$ 
that make up the realisation, is $N$-acyclic).
\eC

\bR
\label{caninfasredprodrealacycrem}
The canonical  infinite realisation of $\H/\I$ as discussed in 
Remarks~\ref{caninfrealrem} and~\ref{caninfasredprodrealrem} 
above is $N$-acyclic for all $N$, i.e.\ its
hypergraph of co-ordinate domains is acyclic. This 
can be seen directly or via its nature as a reduced product with the 
acyclic  free groupoid over $\I$ of
Remarks~\ref{freegroupoidrem} and~\ref{freegroupoidacycrem}.
\eR

\subsection{A universality property for N-acyclic realisations} 
\label{Nunivsec}

\bD
\label{universalrealdef}
A realisation $\A/\H$ of an amalgamation pattern $\H/\I$ is 
\emph{$N$-universal} if every induced substructure of size up to $N$
of $\str{A}$ can be homomorphically mapped into any other (finite or infinite) 
realisation of $\H/\I$.  
\eD

\bP
\label{universalrealprop}
$N$-acyclic realisations are $N$-universal.
\eP

Note that the canonical infinite realisation of $\H/\I$ 
(cf.\ Remarks~\ref{caninfrealrem} and~\ref{caninfasredprodrealrem}),
being $N$-acyclic for all $N$ by
Remark~\ref{caninfasredprodrealacycrem}, 
must be $N$-universal for all $N$. In fact, it does admit full 
homomorphisms into any other realisation, by a countable chain 
argument analogous to the following finite chain argument.

\prf
For given $\H/\I$,  
let $\A/\H = (\H/\I,\str{A},U,(U_s)_{s \in S}, (\pi_{u,s})_{u\in
  U_s})$ be an $N$-acyclic realisation, 
$\B/\H = (\H/\I,\str{B},V,(V_s)_{s \in S}, (\pi_{v,s})_{v\in  V_s})$ 
any other realisation.
By the characterisation of hypergraph acyclicity in
terms of tree decompositions from~\cite{BeeriFaginetal}, we know that 
any induced substructure $\str{C}  = \str{A}\restr C \subset \str{A}$ 
of size $|C| \leq N$ admits a tree decomposition 
$(T,E, \lambda \colon T \rightarrow U)$
according to Definition~\ref{treedecompdef} such that 
$C \subset \bigcup_{t \in T} \lambda(t)$ and 
for every $c \in C$ 
the set $\{ t \in T \colon c \in \lambda(t) \}$ is connected in
$T$. We orient the tree from a chosen root vertex $t_0$ and 
enumerate the vertices of $T$ as $T = \{ t_i \colon i \leq n \}$
in such a manner that each set $T_i = \{ t_j \colon j \leq i \}$ is
connected. Let $u_j = \lambda(t_j) \in U_{s_j}$ so that 
$\str{A}\restr u_j \simeq \str{A}_{s_j}$ via a chart of $\A$,
let $A_i = \bigcup_{j \leq i} u_i$ and $\str{A}_i = \str{A}\restr A_i$. Let
$\str{C}_i = \str{C}\restr A_i$, so that $\str{C}_n = \str{C}$. 
The desired homomorphism $h \colon \str{C} \rightarrow \str{B}$ is
constructed from an increasing chain of 
homomorphisms  $h_i \colon \str{C}_i \rightarrow \str{B}$ for
$i=1,\ldots, n$ such that for each individual $j \leq i$ every
$h_i\restr (C \cap u_j)$ extends to an isomorphism $h_{ij}^\ast$
between $\str{A}\restr u_j \simeq 
\str{A}_{s_j}$ and some $\str{B}\restr v_j$ for $v_j \in V_{s_j}$ that
is induced by corresponding charts of $\A$ and $\B$. Such 
$h_i$ can be found inductively as follows. For $i =0$, let $u_0 = \lambda(t_0) \in U_{s_0}$. We 
may then choose any $v_0 \in V_{s_0}$ in $\B$, and find an isomorphism 
between $\str{A}_0 = \str{A}\restr u_0 \simeq \str{A}_{s_0}$ and 
$\str{B}\restr v_0 \simeq \str{A}_{s_0}$ as induced by 
corresponding charts of $\A$ and $\B$; the restriction of this
isomorphism to $C \cap u_0$ can serve as $h_0$.
For the extension step from
$h_i \colon \str{A}_i \rightarrow \str{B}$ to  
$h_{i+1} \colon \str{A}_{i+1} \rightarrow \str{B}$, 
let $t = t_{i+1}$, $s= s_{i+1}$ and $u = u_{i+1}$ be the building blocks
associated with the new vertex $t_{i+1}$ in $T_{i+1}$. Consider the
immediate predecessor vertex $t_j$ of $t$ in the chosen
orientation of $T$ and let $h_{ij}^\ast(\str{A}\restr u_j) = \str{B}\restr
v_j \simeq \str{A}_{s_j}$ be the local extension of $h_i$ to an
isomorphism.
If $u_j \cap u =
\emptyset$, then $\str{A}_{i+1}$ can be handled like $\str{A}_0$ in
the base case.
Otherwise, since $\B$ also is a realisation of $\H/\I$, 
there is some $v \in V_{s}$ such that $\str{B}\restr v \simeq
\str{A}_s$ intersects with $\str{B}\restr v_j$ in 
$h_{ij}^\ast(u_j \cap u)$ in a manner that is compatible with the charts 
at $u$ and $v$ in $\A$ and $\B$, respectively. 
To see this, we can use 
the fact that $u_j \cap u$ in $\A$ is induced by some $\rho_w$ 
for some walk $w$ from $s_j$ to $s$ in $\I$; and the succession of 
the same sequence of links from $v_j$ in $\B$ leads to a target 
$v$ whose overlap with $v_j$ must contain all of $h_{ij}^\ast(u_j \cap u)$. 
It follows that the extension of $h_i$ by 
the isomorphism between $\str{B}\restr v \simeq
\str{A}_s$ and $\str{A}\restr u \simeq \str{A}_s$, induced by
corresponding charts,
is well-defined in restriction to $C$, 
and meets the requirements for $h_{i+1}$. For
well-definedness, as a homomorphism between relational structures, it is
essential that the $u_i$ form a tree decomposition of an atlas for $\str{C}$: by
Observation~\ref{atlastreedecompobs} there can be no
identities between elements of $\str{C}_i$ and $\str{C}\restr u$ other
than those in $u_j \cap u$, and no relational links between any
elements from $C \cap (u \setminus u_j)$ and from $C_i$. 
\eprf

\section{Two applications} 
\label{symmsec}

\subsection{Hypergraph coverings}
\label{hypcovsec}

As we saw in one of the first motivating examples,
Example~\ref{explodedviewhypex}, any (finite) hypergraph 
$(B,S)$, with $S \subset \mathcal{P}_{\mathrm{fin}}(B)$  
as its set of hyperedges over the vertex set $B$, 
gives rise to an amalgamation pattern that describes an 
exploded view of $(B,S)$. This amalgamation pattern has $(B,S)$ as
a natural albeit trivial realisation. Recall from Definition~\ref{hypgraphdef} that 
we always assume $B = \bigcup S$.
The associated incidence pattern is $\I(B,S)$
with $S$ as its set of sites and $E = \{ (s,s') \in S^2\colon s \not=
s', s \cap s' \not= \emptyset \}$ as the set of links. More formally,
$\I(B,S) = (S, E, \iota, \cdot^{-1})$ is the $2$-sorted encoding of 
the intersection graph of the set $S$ of hyperedges, 
with directed edges and edge reversal. The amalgamation pattern 
that represents the exploded view of $(B,S)$ over this incidence
pattern $\I = \I(B,S)$ is 
\[
\H(B,S)/\I
= \bigl(  \I;
H, \delta \colon H \rightarrow S; P, \eta \colon P \rightarrow E
\bigr)
\]
where 
\bre
\item[--]
$H = \bigcup_{s \in S} s \times \{ s\}$ is the disjoint union of
the hyperedges $s \in S$, partitioned by $\delta$ according to
$\delta(a,s) = s$ so that $\delta^{-1}(s) = s  \times \{ s\}
\subset H$ represents the hyperedge $s \in S$;
\item[--]
$P$ is the set of partial
bijections $p_e = p_{(s,s')} \colon (s \cap s') \times \{ s \}
\rightarrow (s \cap s') \times \{ s' \}$ 
representing the overlap between $s$ and $s'$ in $(B,S)$
by mapping $(b,s)$ to $(b,s')$ 
for pairs $(s,s') = e \in E$ of hyperedges with a
non-trivial intersection.
\ere

Clearly  $\H(B,S)$ is strongly coherent and simple; in fact the action of
$\I^\ast$ through compositions $\rho_w$ along walks in $\I = \I(B,S)$
is very simple: if $\iota(w) = (s,s')$, then $\rho_w$ is the natural bijection 
between $(s\cap s') \times \{ s\}$ and   $(s\cap s') \times \{ s' \}$. 
Also the symmetries of $\H(B,S)$ in the sense of Definition~\ref{amalgautdef}
are in bijective correspondence with the symmetries of $(B,S)$ as a hypergraph, 
viz.\ with permutations of $B$ that fix $S$ as a set of subsets of $B$
and thus also preserve overlaps. 
Obviously $(B,S)$ itself can be cast as a realisation of $\H(B,S)/\I$, 
by associating each hyperedge $s \in S$ with the natural chart into 
$\delta^{-1}(s) = s \times \{ s \}$ in $\H$. So $(B,S)$ becomes
the hypergraph of the co-ordinate domains of the atlas that connects
$(B,S)$ with its exploded view. In fact $(B,S)$ corresponds to the
quotient $\H(B,S)/{\approx}$ of the strongly 
coherent and simple $\H(B,S)/\I$ as in 
Observation~\ref{strongcoherelqutientlem}.
Proposition~\ref{hypcovprop} below shows 
that, quite naturally, every realisation of $\H(B,S)/\I$ 
is a hypergraph covering of $(B,S)$ in the sense of the
following definition.

\bD
\label{hypcovdef}
A hypergraph \emph{homomorphism} 
between hypergraphs $(A,U)$ and
$(B,S)$ is given by a map $\pi \colon A \rightarrow B$ whose 
restriction to each 
hyperedge $u \in U$ bijectively maps that $u$ onto an image hyperedge 
$\pi(u) \in S$. We naturally extend $\pi$ to the second sort (of
hyperedges) and denote homomorphisms as in $\pi \colon (A,U) \rightarrow (B,S)$.

A homomorphism $\pi \colon (A,U) \rightarrow (B,S)$ is a
\emph{covering} 
of the hypergraph $(B,S)$ by the hypergraph $(A,U)$ if it maps $U$
surjectively onto $S$ (and hence also $A$ onto $B$), and satisfies the following
lifting property for overlaps between hyperedges: 
for every pair of hyperedges $(s,s')$ of $(B,S)$  with nontrivial
overlap, and every $u \in U$ with $\pi(u) = s$ there is 
some $u' \in U$ s.t.\ $\pi(u') = s'$ and  $s \cap s' = \pi(u\cap u')$
(equivalently: s.t.\ $\pi$ restricts to a bijection of $u \cup u'$ onto
$s \cup s'$).

A covering $\pi \colon (A,U) \rightarrow (B,S)$ is \emph{fully
  symmetric} over $(B,S)$ if every symmetry of
$(B,S)$ (or of $\H(B,S)$) extends to a symmetry of 
$\pi \colon (A,U) \rightarrow (B,S)$ and if the 
subgroup of $B$-rigid symmetries (symmetries that fix  
$B$, and therefore $\H(B,S)$, pointwise) acts transitively on $\pi^{-1}(s) \subset U$ 
for every $s \in S$.
\eD

Note that our notion of a covering is what in other contexts is
called a \emph{branched covering}, since the multiplicities of
overlaps between hyperedges in $\pi^{-1}(s) \subset U$ and
$\pi^{-1}(s') \subset U$ may be greater than one. Simple examples show
that coverings of a given (finite) hypergraph by (finite or infinite)
hypergraphs that are acyclic, or even just $N$-acyclic for specific
$N$, may be necessarily branched. (The situation for graphs is quite
different: finite unbranched coverings by $N$-acyclic graphs are
always available, and can be obtained as direct products with finite
groups of large girth~\cite{OttoAPAL04}.)

\bE
\label{cartwheelex}
Consider the hypergraph with $4$ vertices and $3$ hyperedges of size
$3$ that corresponds to a triangulation of a triangle with one extra central
vertex. This hypergraph is not $3$-acyclic, and every non-trivial
covering hypergraph must have higher multiplicities at the central vertex
(infinite for a properly acyclic covering).
\eE

\bP
\label{hypcovprop}
Every realisation $\A$ of $\H(B,S)/\I$ for $\I = \I(B,S)$ 
induces a covering of the hypergraph $(B,S)$ by the hypergraph $(A,U)$ 
associated with the atlas of the realisation $\A$.
This covering is fully symmetric over $(B,S)$ if the realisation $\A$ is fully
symmetric over $\H(B,S)$, and the covering hypergraph is $N$-acyclic
if the realisation $\A$ is $N$-acyclic.
\eP

\prf
Let $\A$ be  a realisation of $\H(B,S)/\I$ with atlas $(A,U,(U_s),
(\pi_{u,s}))$ so that every $\pi_{u,s}$ for $u \in U_s$ bijectively
maps $u_s\subset A$ to $s \times \{ s \} \subset H$. Let $(\pi_{u,s})_1$
be the composition of $\pi_{u,s}$ with projection to the first
component of the image, which is a bijection from $u \in U_s$ onto~$s$.
If $a \in u \cap u'$ for $u \in U_s$ and $u' \in U_{s'}$, then the
properties of $\A$ as a realisation imply that, for some walk $w$ from 
$s$ to $s'$ in $\I$, $\pi_{u',s'}(a) = \rho_w(\pi_{u,s}(a))$ in
$\H(B,S)$, which implies that $(\pi_{u',s'})_1(a) = (\pi_{u,s})_1(a)$
by the special nature of the $\rho_w$ in $\H(B,S)$. It follows that
the global union of all the $(\pi_{u,s})_1$ is a well-defined
surjection $\pi \colon A \rightarrow B$ of $A$ onto $B = \bigcup S$. 
It remains to show that $\pi \colon (A,U) \rightarrow (B,S)$ is a hypergraph
covering. The claims regarding symmetry and $N$-acyclicity are
obvious, and  
clearly $\pi$ is a hypergraph homomorphism that is surjective also at
the level of hyperedges. The lifting property for non-trivial overlaps between
hyperedges $s,s' \in S$ at some $u \in U_s$ is obvious from the
corresponding condition on realisations: since $s \cap s' \not=
\emptyset$ and $s \not= s'$, $e=(s,s')$ is a link in $\I = \I(B,S)$,
which has to be realised  in $\A$ 
by the overlap with some $u' \in U_{s'}$, and compatibility with $\rho_e$ (in the sense of condition~(i) of
Definition~\ref{realdef}) means that $u \cap u'$ precisely matches
$s\cap s'$.
\eprf

\bC
\label{acyccovcorr}
For every $N$, every finite hypergraph admits fully symmetric finite coverings by
$N$-acyclic hypergraphs.
\eC

\bE
Consider the simple hypergraph consisting of
the facets of the $3$-simplex; this hypergraph is associated 
with the faces of the tetrahedron, or it may be seen as just the 
complete $3$-uniform hypergraph on four vertices. Its Gaifman
graph does not have any chordless cycles of length greater than~$3$,
but conformality is violated by the $4$-clique of all its vertices. 
A finite $N$-acyclic covering for 
$N \geq 4$ cannot be unbranched since it must unravel the local
$3$-cycles of hyperedges that meet in a single vertex into cycles of 
lengths $3n$ for some $n > 1$. While it is not hard to visualize a
natural fully symmetric covering based on a quotient of a 
hexagonal grid pattern (locally two-fold, for $n=2$), which is $5$-acyclic, 
there seems to be no similarly obvious ad-hoc construction for 
higher degrees of acyclicity as guaranteed by Corollary~\ref{acyccovcorr}.
\eE

\subsection{Lifting local to global symmetries}
\label{HLsec}

Recall that a partial isomorphism $p \in
\mathrm{part}(\str{A},\str{A}')$ is an isomorphism between induced
substructures of $\str{A}$ and $\str{A}'$. In case $\str{A}= \str{A}'$, we
speak of \emph{partial automorphisms}, which are partial or local
symmetries within structure $\str{A}$, i.e.\ elements of the 
inverse semigroup of partial bijections that respect the
interpretations of given relations of $\str{A}$. It is not hard to see that a 
collection of partial automorphisms in $\mathrm{part}(\str{A},\str{A})$  
can simultaneously be extended to global automorphisms in the
automorphism group $\mathrm{aut}(\str{B})$ of some infinite
extension $\str{B} \supset \str{A}$. A straightforward 
construction of such solutions to the \emph{extension problem
for partial automorphisms} could proceed along the lines of
Remark~\ref{caninfasredprodrealrem} for the construction of 
infinite realisations as reduced products. The problem becomes more
subtle if, for finite $\str{A}$ and a collection of
its partial automorphisms, one asks for \emph{finite} solutions or
for \emph{finite} solutions within a restricted class of
structures. Another important variant of the problem asks 
for \emph{finite} solutions within a restricted class $\mathcal{C}$ 
provided that at least infinite solutions exist in $\mathcal{C}$. Note
the conditional character of this variant, which is in the style
of a finite model property for the extension task.  Classes
$\mathcal{C}$ of relational structures that satisfy this version are said to have the 
\emph{extension property for partial automorphisms} (for short: EPPA)
in terminology introduced in~\cite{HerwigLascar}.

Classical results for the unconditional extension task 
in restriction to finite structures were first obtained for 
the class of all finite graphs by Hrushovski in~\cite{Hrushovski}, then 
for the class of all finite structures in any fixed finite relational signature 
by Herwig in~\cite{Herwig}. A strongest general result in this direction
is the result of Herwig and Lascar in~\cite{HerwigLascar},
Theorem~\ref{HLthm} below: any class of relational structures 
defined in terms of finitely many
finite \emph{forbidden homomorphisms} satisfies (EPPA). 
This result subsumes the earlier, unconditional results for classes of
finite structures, since infinite solutions are always available in
those cases. For the result for finite graphs, a very elegant, 
elementary and direct proof is presented in~\cite{HerwigLascar}, which
avoids both the substantial intricacies of the much stronger general
(EPPA) result in~\cite{HerwigLascar} and the non-trivial algebraic arguments 
of the original proof in~\cite{Hrushovski}. 

As it turns out, fully symmetric finite realisations of naturally
induced amalgamation patterns offer an alternative direct route to 
particularly generic finite solutions for the extension problem for
finite automorphisms. This connection is especially striking in the general case of 
the conditional Herwig--Lascar result, where the insights 
from Section~\ref{Nunivsec}
into the connection between local acyclicity and local genericity up
to homomorphisms play out in a particularly nice way.

\bD
\label{EPPAdef}
Let $\sigma$ be a finite relational signature, $\mathcal{C}$ a class of 
$\sigma$-structures. 
An \emph{instance} of the \emph{extension problem for partial
  automorphisms} over $\mathcal{C}$ is a pair 
$(\str{A}_0,P)$ consisting of a finite $\sigma$-structure $\str{A}_0
\in \mathcal{C}$ and a set $P$ of partial isomorphisms in $\str{A}_0$ 
(partial automorphisms, local symmetries of $\str{A}_0$); 
a \emph{solution} for this instance is any $\str{A}$ in $\mathcal{C}$ 
into which $\str{A}_0$ embeds isomorphically such that 
the image of each $p \in P$ extends to a full automorphism of $\str{A}$. 

The class $\mathcal{C}$ has the 
\emph{extension property for partial automorphisms}  (EPPA) if 
every instance $(\str{A}_0,P)$ of the extension problem over
$\mathcal{C}$ that admits any solution also admits a finite 
solution in $\mathcal{C}$.
\eD

In the following we shall always assume that the set $P$ 
of partial automorphisms of an instance is closed 
under inverses: for each $p \in P$, also $p^{-1} \in P$.
This is w.l.o.g.\ in the sense that $p^{-1}$ can be adjoined to $P$ 
without  changing the set of solutions.

For classes $\mathcal{C}$ like the class of all graphs or of all
$\sigma$-structures, the condition that an instance has a solution
in $\mathcal{C}$ is universally fulfilled, so that (EPPA) for such
classes boils down to the requirement that every instance over 
$\mathcal{C}_{\mathrm{fin}}$ admits a solution in
$\mathcal{C}_{\mathrm{fin}}$, where $\mathcal{C}_{\mathrm{fin}}$
stands for the subclass of finite structures in $\mathcal{C}$.

\bD
\label{forbidhomdef}
A class $\mathcal{C}$ of $\sigma$-structures is defined in
terms of \emph{forbidden homomorphisms} if there is some 
finite set $X$ of finite $\sigma$-structures such that 
$\mathcal{C}$ is the class of those $\sigma$-structures $\str{A}$ 
that admit \emph{no} homomorphism $h \colon \str{X} \rightarrow \str{A}$
from any $\str{X} \in X$. 
\eD

\bT[Herwig--Lascar]
\label{HLthm}
Any class $\mathcal{C}$ of $\sigma$-structures that is defined in
terms of forbidden homomorphisms 
has the extension property for partial automorphisms (EPPA).
\eT

Positive examples include the classes of all (finite) graphs~\cite{Hrushovski}, of 
(finite) $K_n$-free graphs~\cite{HerwigLascar}, of all (finite)
$\sigma$-structures~\cite{Herwig}. 
Examples of classes $\mathcal{C}$ with (EPPA) that are not
defined in terms of forbidden homomorphisms but covered by the 
immediate specialisation of Theorem~\ref{eppathm} below, are the classes of
(finite) conformal $\sigma$-structures, i.e.\ structures for which
all cliques in the Gaifman graph are generated by single relational
atoms, as treated in~\cite{HO}. It appears to be still open whether the class of all
tournaments or the class of all $4$-cycle-free graphs have (EPPA).

\bD
\label{symmsoldef}
We call a solution $\str{A}$ to the instance $(\str{A}_0,P)$ 
of the extension problem a \emph{fully symmetric solution} 
if it possesses an atlas $(\str{A},U,(\pi_u)_{u \in U})$ of co-ordinate maps
$\pi_u \colon \str{A}\restr u \simeq \str{A}_0$ onto the single co-ordinate
structure $\str{A}_0$ such that 
\bre
\item
the automorphism group of $\str{A}$ acts transitively on the set 
$U$ of the co-ordinate domains in a manner compatible with 
the charts (i.e.\ by symmetries of the
atlas that are rigid w.r.t.\ the single co-ordinate structure);
\item
any two co-ordinate domains $u,u' \in U$ with non-trivial intersection
are related by the action of a composition of partial isomorphisms
from $P$ in the sense that $p_{u'} \circ p_u^{-1} = p_n \circ \cdots
\circ p_1$ for a sequence $w = p_1\cdots p_n \in P^\ast$.%
\footnote{The equality of partial compositions includes the 
equality of the domains. Recall that we assume $P$ to be closed under
inverses; $P^\ast$ can be understood as either the set of finite
words over the alphabet $P$, or, equivalently as the set of 
walks in an incidence pattern $\I_P = ( \{ 0 \}, P)$ 
with a site $0$ and loops for $p \in P$.}
\ere
\eD

\bT
\label{eppathm}
For the class $\mathcal{C}$ of all $\sigma$-structures and any $N \in \N$:
every instance $(\str{A}_0,P)$ for $\str{A}_0 \in \mathcal{C}$ 
admits finite solutions $\str{A} \in \mathcal{C}_{\mathrm{fin}}$ which
are fully symmetric and whose atlas is $N$-acyclic.
\eT

\prf
The desired finite solution is obtained from a finite, fully symmetric and
$N$-acyclic realisation of an amalgamation pattern 
$\H(\str{A}_0,P)$ derived from the extension task $(\str{A}_0, P)$ as
follows. We use the incidence pattern  $\I =
\I(\str{A}_0,P) = \I_P$ with the singleton set $S = \{ 0 \}$ 
for its sort of sites, and edges $e_p \in E$ 
for each $p \in P$, with constant $\iota(p) = (0,0)$
and with edge reversal according to $e_p^{-1} =
e_{p^{-1}}$. We may identify the set $P$ with the sort $E$ of links of $\I_P$.
Formally, the amalgamation pattern $\H(\str{A}_0,P)$ then is 
\[
\H(\str{A}_0,P)/\I_P = \bigl( \I; \str{A}_0, 0 \colon A_0 \rightarrow \{ 0 \}; P,
\mathrm{id} \colon P \rightarrow P \bigr)
\]
with the constant function $0$ (as $\delta$) for the trivial $S$-partition of $H =
A_0$, and the identity function on $P$ (for $\eta$) that labels the partial
isomorphism $p \in \mathrm{part}(\str{A}_0,\str{A}_0)$ by its name,
which is the link $p = e_p \in E = P$. 

We claim that every realisation $\A/\H = (\H,\str{A},U, (\pi_{u})_{u\in U})$ of
$\H= \H(\str{A}_0,P)$ that is fully
symmetric over $\H(\str{A}_0,P)$ is a solution to the extension
problem for $(\str{A}_0,P)$. More specifically, the $\sigma$-structure $\str{A}$ of $\A$
embeds an isomorphic copy of $\str{A}_0$ as a substructure 
in such a way that the image of every 
$p \in \mathrm{part}(\str{A}_0,\str{A}_0)$
extends to some automorphism of $\str{A}$.

Any embedding of $\str{A}_0$ according to one of the charts
$\pi_u \colon \str{A}\restr u \simeq \str{A}_0$ can be used. We fix
some $u_0 \in U$ and regard $\str{A}\restr u_0$ as the distinguished 
copy of $\str{A}_0$ in $\str{A}$. Let its chart be $\pi_0 := \pi_{u_0} \colon
\str{A}\restr u_0 \simeq \str{A}_0$. The incarnation of 
$p \in \mathrm{part}(\str{A}_0,\str{A}_0)$ in this distinguished
copy of $\str{A}_0$ is $\pi_0^{-1} \circ p \circ \pi_0\in 
\mathrm{part}(\str{A}\restr u_0,\str{A}\restr u_0)$.

For $p \in P$ choose $u \in U$ for $e_p$ as in
condition~(i) of Definition~\ref{realdef} for realisations, such that 
\begin{equation}
\label{(ast)} 
\pi_u \circ \pi_{0}^{-1} = p \in
\mathrm{part}(\str{A}_0,\str{A}_0),
\end{equation}
where the composition on the left is a partial composition of the
restrictions that fit together precisely in $\pi_{0}^{-1}(\mathrm{dom}(p)) =
\pi_u^{-1}(\mathrm{image}(p))$.

As $\A$ is fully symmetric over $\H$, it possesses an $\H$-rigid
symmetry $\pi$ that maps $u$ to $u_0$ in a manner that is compatible with 
$\pi_u$ and $\pi_0$: 
\begin{equation}
\label{(dagger)} 
\pi_{0} \circ \pi = \pi_{u}, 
\end{equation}
where the
composition on the left is defined on all of $u$.

Let $a \in u_0$ be in the domain of $\pi_0^{-1} \circ p \circ \pi_0$
(i.e.\ of $p$ in the distinguished copy of $\str{A}_0$ in $\str{A}$). 
Then $p(\pi_{0}(a)) = \pi_u(a)$ by~(\ref{(ast)}), so that 
$\pi_{0}^{-1} (p(\pi_{0}(a))) = \pi_0^{-1}(\pi_u(a)) = \pi(a)$
by~(\ref{(dagger)}). So the symmetry $\pi$ of $\A$ extends 
$\pi_0^{-1} \circ p \circ \pi_0\in \mathrm{part}(\str{A}\restr u_0,\str{A}\restr u_0)$
to an automorphism of $\str{A}$.

The extra symmetry and acyclicity properties of the solution $\str{A}$
with its natural atlas are immediate from corresponding properties of $\A$.
\eprf

\bO
\label{groupobs}
Considering the manner in which suitable finite realisations 
of $\H/\I$ were obtained as reduced products with suitable groupoids
$\G/\I$ in Section~\ref{redprodrealsec}, it is worth noting that 
the incidence pattern $\I_P$ for $\H(\str{A}_0,P)/\I_P$ is such that 
any groupoid $\G/\I_P$ is in fact a Cayley \emph{group}.
\eO

The following strengthening of the Herwig--Lascar result of
Theorem~\ref{HLthm} is a consequence of Theorem~\ref{eppathm}
and our considerations about $N$-acyclicity and $N$-universality in 
Section~\ref{Nunivsec}. In analogy with Definition~\ref{universalrealdef} we
say that a $\sigma$-structure 
$\str{A}$ is \emph{$N$-universal} within a class $\mathcal{C}_0$ of
$\sigma$-structures (in our case, a class of solutions to an instance
of an extension problem) 
if every substructure of $\str{A}$ of size up to $N$ admits homomorphisms
to every structure in $\mathcal{C}_0$.

\bC
\label{eppacorr}
Consider the instance $(\str{A}_0,P)$ of the extension problem in 
the class of all $\sigma$-structures and let $\str{A}$ be any 
finite, $N$-acyclic and fully symmetric solution as in Theorem~\ref{eppathm}.
\bae
\item
Then $\str{A}$ is $N$-universal within the class
$\mathcal{C}_0 = \mathcal{C}_0(\str{A}_0,P)$ 
of all solutions to the instance $(\str{A}_0,P)$.
\item
Let the class $\mathcal{C}$ be defined in terms of forbidden
homomorphisms and $N$ sufficiently large, viz.\ an upper bound for the sizes of the forbidden homomorphisms
that define $\mathcal{C}$. Then $\str{A} \in
\mathcal{C}_{\mathrm{fin}}$ if the instance $(\str{A}_0,P)$ has any 
finite or infinite solution in $\mathcal{C}$. This implies 
(EPPA) for the class $\mathcal{C}$.
\eae
\eC

\prf
Claim~(b) is an immediate consequence of~(a). 
We deal with the more specific claim~(a). 
 
Let $\str{A}$ be a finite, $N$-acyclic, fully symmetric 
solution as in the proof of Theorem~\ref{eppathm}
with the atlas $(\str{A},U,(U_s), (\pi_{u})_{u\in U})$ 
of charts $\pi_u \colon \str{A}\restr u \simeq \str{A}_0$.
Let $\str{B}$ be any other solution for the instance $(\str{A}_0,P)$ of the
extension problem over the class of all $\sigma$-structures, with
$\str{A}_0 \simeq \str{B}_0 \subset \str{B}$ for the embedded copy of
$\str{A}_0$. For $p \in P$, let $\pi_p \in \mathrm{Aut}(\str{B})$ 
be an extension of the image of $p \in
\mathrm{part}(\str{A}_0,\str{A}_0)$ in $\str{B}_0 \subset \str{B}$.
We need to construct, for a given $\str{D} \subset \str{A}$ of
size up to $N$, a homomorphism $h \colon \str{D} \rightarrow \str{B}$. 

As in the proof of Proposition~\ref{universalrealprop}, we can
use a tree decomposition of $\str{D}$ that is induced by the
$N$-acyclic atlas of $\str{A}$. This tree decomposition will govern 
an inductive choice of a sequence of homomorphisms
$h_i$ for substructures $\str{D}_i \subset \str{D}$. 
In more detail: By
$N$-acyclicity of the atlas,  
$\str{D}  = \str{A}\restr D \subset \str{A}$ for $|D| \leq N$,  
$\str{D}$ admits a tree decomposition  
$(T,E, \lambda \colon T \rightarrow U)$ such that 
$D \subset \bigcup_{t \in T} \lambda(t)$ and such that 
the sets $\{ t \in T \colon d \in \lambda(t) \} \subset T$ are
connected for all $d \in D$. Let $t_0$ be the root of $T$ w.r.t.\ a
chosen orientation and enumerate $T$ as $T = \{ t_i \colon i \leq n \}$
such that all initial subsets $T_i = \{ t_j \colon j \leq i \} \subset T$
are connected. Let $\lambda(t_j) = u_j \in U$ so that 
$\pi_j := \pi_{u_j} \colon \str{A}\restr u_j \simeq \str{A}_0$.
Let $A_i = \bigcup_{j \leq i} u_i$, $\str{A}_i = \str{A}\restr A_i$,
and $\str{D}_i = \str{D}\restr A_i$.

The homomorphism $h \colon \str{D} \rightarrow \str{B}$ is the final member
of an increasing sequence of homomorphisms $h_i \colon \str{D}_i
\rightarrow \str{B}$ for $i=1,\ldots, n$ with this additional property:
for each $j \leq i$, 
$h_i\restr (D \cap u_j)$ extends to an isomorphism $h_{ij}^\ast \colon 
\str{A}\restr u_j \simeq \str{B}_j \subset \str{B}$ where $\str{B}_j
\simeq \str{A}_0$ is an image of $\str{A}_0 \simeq \str{B}_0 \subset \str{B}$ under
some composition $\pi_{w_j}$ of the designated extensions $\pi_p$ of the $p \in P$
to automorphisms of $\str{B}$, corresponding to some word $w_j = p_1
\cdots p_m\in P^\ast$, according to 
\[
\str{B}_j = \pi_{w_j}(\str{A}_0) = (\pi_{p_m} \circ \cdots \circ
\pi_{p_1}  ) (\str{A}_0) \subset \str{B}.
\]

Such $h_i$ can be chosen inductively as follows. 
For $i =0$, let $u_0 = \lambda(t_0) \in U$. We may combine 
the isomorphic embedding of $\str{A}_0 \simeq \str{B}_0 \subset \B$ 
and the isomorphism between $\pi_{u_0} \colon 
\str{A}\restr u_0 \simeq \str{A}_0$ to obtain the desired isomorphism
$h_{00}^\ast$, and use its restriction to $D \cap u_0$ as
$h_0$.

For the extension step from
$h_i \colon \str{A}_i \rightarrow \str{B}$ to  
$h_{i+1} \colon \str{A}_{i+1} \rightarrow \str{B}$, let 
$u = \lambda(t_{i+1})$ and consider the
immediate predecessor $t_j$ of $t_{i+1}$ in $T$. Let 
$h_{ij}^\ast(\str{A}\restr u_j) = \str{B}_j \simeq \str{A}_0$ 
be the local extension of $h_i$ to an isomorphism, where 
$\str{B}_j = \pi_{w_j}(\str{A}_0)\subset \str{B}$.

Since $\str{A}$ is a fully symmetric solution, 
$u_j \cap u$ in $\A$ is induced by some $\rho_{w'}$ 
for a word $w' \in P^\ast$. The corresponding operation of 
the composition $\pi_{w'}$ acts as an automorphism on $\str{B}$ and  
thus produces an overlap between $\str{B}_j = \pi_{w_j}(\str{A}_0)$
and $\str{B}_{i+1} := \pi_{w_jw'} (\str{A}_0)  = \pi_{w'}(\str{B}_j)$ 
that contains the $h_{ij}^\ast$-image of $\str{A}\restr (u_j \cap u)$.
We can use this $\str{B}_{i+1}$, the natural
isomorphism induced by the chart of $\str{A}\restr u$ in the atlas of
$\str{A}$ and the automorphism $\pi_{w_jw'}^{-1}$ that maps
$\str{B}_{i+1}$ onto $\str{B}_0 \subset \str{B}$ as $h^\ast_{i+1,i+1}$
and join its restriction to $u \cap D$ to $h_i$ in order to obtain $h_{i+1}$.
The argument for well-definedness as a homomorphism between relational
structures is strictly analogous to that in the proof of 
Proposition~\ref{universalrealprop}, based on the properties of a tree decomposition.
\eprf

\bibliographystyle{plain}

\end{document}